\documentclass{article}

\usepackage{setspace}

\usepackage{amsmath}

\usepackage{mathtools}

\begin{document}

                                          \title  {On non-periodic and non-dense billiard trajectories Part 1}

\author{George William Tokarsky\\Math Department, University of Alberta}\maketitle
\singlespace\

In a paper [1] from 1983, G.A. Galperin examined the existance of non-periodic and not everywhere dense billiard trajectories in
various polygons and in particular in a right triangle. The particular example that he gave involving a right triangle was based on a sequence of mirror images of an isosceles triangle $ABC$ with $<A=<B=x$ which in this paper we will measure in degrees as in Fig. 1 and which had been accepted [1,3 to 9] for the last thirty years as the definitive example of the existance of a non-periodic and not everywhere dense billiard trajectory in a triangle. Unfortunately that result turned out to be incorrect. This was proved by the author in [2]. The significance of this is that the question of whether every triangle, indeed every isosceles triangle has only periodic or dense trajectories or trajectories that end at a vertex is still an open problem.

\doublespace
Figure 1 

\singlespace
 We first repeat here the short proof in [2]  using a similar but not the same notation as in Galperin's original paper.  
Using the well known process of straightening out a billiard trajectory, Galperin's triangle example was constructed from 44 copies of isosceles triangle $ABC$ using a sequence of reflections as in Fig. 1 in which the first and last reflection is in side $AC$ and all alternate reflections are in side $AB$ and where $A_{22}B_{22}$ and $A_{16}B_{16}$ are parallel to the base $A_{0}B_{0}$. Points labelled $A_{i}$, $B_{i}$, $C_{i}$, $L_{i}$, $M_{i}$, $P_{i}$, $S_{i}$ all represent the same points $A,B,C,L,M,P,S$ on triangle $ABC$. Two parallel corridors through the interior of this tower of triangles were then formed, first using the line $C_{1}C_{22}$ and then introducing a parallel through $C_{16}$ and another parallel line through $C_{17}$ as shown. This construction sets up an exchange of segments on $PL$ in which $PS$ $\rightarrow$ $ML$ and $SL$ $\rightarrow$ $PM$. If the ratio of $PS$ to $SL$ which is the same as the ratio of the widths of the two corridors is irrational then any billiard trajectory within and parallel to the corridors would be non-periodic and everywhere dense in the corridors but not everywhere dense in the given isosceles triangle. Since this ratio is a continuously differentiable function of the acute angle $x$ of the isosceles triangle, there would be an angle $x$ for which this ratio is irrational provided this ratio is \textbf{non-constant}.

\doublespace 
\singlespace However as it happens this ratio turns out to be a constant and in fact equals $1/2$ for every value of the acute angle $x$ for which the corridors can be formed. To calculate this ratio note that the angle opposite the base $AB$ and all its mirror images are labelled as "C" vertices. We can then calculate the coordinates of every "C" vertex using  a little trigonometry together with the algorithm one from section 6 and under the assumption that the distance between any two successive "C" vertices is one unit.

If we let the coordinates of $C_{1}$ be (0,1) then the coordinates of $C_{16}$ are $(sin4x,5+8cos2x+3cos4x)$ of $C_{17}$ are 
$(sin4x,6+8cos2x+3cos4x)$ and the coordinates of $C_{22}$ are $(sin2x+sin4x,8+11cos2x+3cos4x)$. Now the line through $S_{0}$ and $L_{16}$  is exactly the line through $C_{1}$ and $C_{22}$ and has slope
$m=(7+11cos2x+3cos4x)/(sin2x+sin4x)$ and equation $y=mx+1$. This means that the other two parallel lines have the same slope and then the ratio of the widths of the two corridors is the same as the ratio that is cut off by any transversal of the two corridors.

In particular if we extend the line segment $C_{17}$$C_{16}$ and find its point of intersection D with the line $C_{1}$$C_{22}$, we find that its coordinates are $(sin4x,1+msin4x)$ and the distance from $C_{16}$ to D is $4+8cos2x+3cos4x-msin4x$ which equals 1 since $(3+8cos2x+3cos4x)(sin2x+sin4x) = (7+11cos2x+3cos4x)sin4x$ which can be verifed by multiplying out and using the trig identity $2sinAcosB=sin(A+B)+sin(A-B)$.

\doublespace
Figure 2 
\singlespace

This means that the ratio of the widths of the two corridors is exactly 2 to 1 for every value of the acute angle $x$, which is a constant rational number and thus Galperin's example fails. Note we don't have to worry about the case where $sin2x+sin4x =0$ as then the lines would be vertical and if the corridors could be formed (which they can't) the ratio would be 1 to 1 which is again rational.

One might then hope to salvage Galperin's method by finding a different sequence of mirror images which produces an irrational ratio. It turns out that this is also doomed to failure as all corridors that satisfy the conditions as in Galperin's original paper produce a constant rational ratio as it is the purpose of this paper to eventually show. 

\subsection*{1. Tower of mirror images of an isosceles triangle}

Let triangle ABC be isosceles with $<A=<B=x$ and by convention \textbf{oriented counterclockwise} from A to B to C. A sequence of mirror images in the sides of triangle ABC will be called a \textbf{tower} if 

1. every alternate mirror image is in side AB and

2. its  first and last mirror image is in side AC or BC.  

\noindent Any tower must then have an even number of triangles in it.  If we number the triangles consecutively the first or starting triangle will always be oriented counterclockwise and so will every odd numbered triangle while every even numbered triangle will be oriented clockwise including the last triangle. All successive mirror images of A,B and C will be called A,B, and C points. \textbf{Successive "C" points will always be assumed to be one unit distance apart} and are successively labelled $C_{1}$ (we start at $C_{1}$ instead of $C_{0}$) to $C_{n}$ while the base AB and its mirror images are successively labelled 
$A_{0}$$B_{0}$ to $A_{n}$$B_{n}$. This means the orginal triangle is labelled $A_{0}$$B_{0}$$C_{1}$. Note that this indexing differs from that used by Galperin. Observe that each "C" point has an unique label while the "A" and "B" points can have multiple labels. For example if the first mirror image is in side AC then that side has the labels $A_{0}$$C_{1}$ as well as the label $A_{1}$$C_{1}$. Any C point corresponding to a mirror image in AC will be called a \textbf{black point} and if it corresponds to a mirror image in BC will be called a \textbf{blue point}. Each C point in the tower will have an unique color.  Further the C points are in fact ordered in an increasing order that follows the ordering of the formation of the sequence of mirror images of the tower. So $i<j$ if and only if $C_{i}$ was formed before $C_{j}$ in the sequence of mirror images.

\doublespace
Figure 3 
\singlespace

\subsection*{2. Rhombus poolshots and rhombus towers}

Again let triangle ABC be isosceles with $<A=<B=x$ and oriented counterclockwise from A to B to C.
Given a finite billiard trajectory  or poolshot, we will call it a \textbf{rhombus pool shot} or \textbf{rhombus billiard trajectory} if 

1. every alternate reflection is in the base AB 

2. its first and last reflection is in side AC or BC and

3. it starts and finishes at side AB and doesn't hit a vertex.

\noindent If we \textbf{straighten} out this poolshot then we get a corresponding \textbf{rhombus tower} of mirror images with all A points occuring on one side and all B points occuring on the other side  of the straightened poolshot (called the \textbf{associated rhombus poolshot}). 
Now observe that the blue C points and the black C points are also on opposite sides of the straightened  poolshot with all the blue points on the same side as all the A points and all the black points on the same side as all the B points and it follows that the convex hull of the blue C points and the convex hull of the black 
C points are disjoint. This further means we cannot have a $blue-black-blue$ collinear situation where a black point is between  two blue points. Similarly we cannot have a $black-blue-black$ collinear situation. 

Caution: A rhombus tower is a tower but a tower with property 3 above need not be a rhombus tower for example if the convex hulls of the blue and black points intersect. This prevents any poolshot being associated with the tower.

Notice since the poolshot starts at AB, if we straighten it out upwards in \textbf{standard position} with the base AB placed horizontal with A to the left of B and C above the base, then all A points and blue C points are on the $left$ $side$ and all B points and black C points are on the $right$ $side$.

\doublespace
Figure 4
\singlespace

\noindent
\textbf{Convention:} Any straightened rhombus poolshot can be viewed as forming the positive Y coordinate axis in an XY coordinate system by introducing a perpendicular X axis through the starting point of the poolshot on side AB. With this convention all A points will be on the left side and all B points will be on the right side of the poolshot.

\subsection*{3. The Non-Rhombus Tower Collinear Test} As previously noted given a finite tower if we have a blue-black-blue collinear situation with the black point between the two blue points (or a black-blue-black collinear situation) then no rhombus poolshot can create this tower of mirror images as the convex hulls of the blue and black points would not be disjoint. Vectorially if $v=(a,b)$ is a non-zero vector from the first blue point to the black point and $w=(c,d)$ is a non-zero vector from the black point to the last blue point, the three points are collinear if and only if $ad=bc$. Hence we get the 

\textbf{The Non-Rhombus Tower Collinear Test:} A tower is not a rhombus tower if it has three collinear blue-black-blue (or black-blue-black) C points in that betweeness order or equivalently if $ad=bc$ where $v=(a,b)$ is a vector from the first C point to the second and $w=(c,d)$ is a vector from the second to the third C point as there would be no associated rhombus pool shot.

\doublespace
\singlespace
Note that if the three C points are in the correct betweeness order then a and c have the same sign and b and d have the same sign. Also observe that a=0 if and only if c=0 and similarly for b and d and that one of a and b is non-zero and one of c and d is non-zero.

\subsection*{4. Side Sequences and Codes}

Since we will be dealing almost exclusively with rhombus billiard trajectories or rhombus poolshots in an isosceles triangle ABC, we need a compact notation to describe the successive sides that are hit.

The usual notation is to label the sides of triangle ABC as follows AB=1, AC=3 and BC=2. Then for example 13131212131 represents a \textbf{side sequence} of reflections in those respective sides. We will use the convention that the first 1 tells us that the poolshot starts at side AB and the last 1 tells us that it ends at side AB and no further reflections are considered. Observe that no two consecutive integers are the same and since every alternate reflection is in the base AB then every alternate integer will be a 1. 

(We will also use the same side sequence notation to describe the sequence of reflected images of triangle ABC that form a tower which may or may not be a rhombus tower according as to whether it is associated with a poolshot or not. For example the tower in Fig. 5 with x=17 degrees is described by the side sequence 1313131313131 and isn't a rhombus tower. We make a similar  convention as above that the first 1 just represents that first appearance of triangle ABC (and not a reflection) and the last 1 represents the last appearance of triangle ABC and not a reflection.)

\doublespace
Figure 5
\singlespace

The notation can become cumbersome to read if we are dealing with a long side sequence of reflections for example 131312121213131312121313121312131. Hence we introduce its corresponding \textbf{code sequence} 4 6 6 4 4 2 2 2 2 or 4
$6^2$ $4^2$ $2^4$ in exponential code. The first integer 4 represents how many times 1 and  3 are switched at the start, the second integer 6 then represents the following number of succesive switches of 1 and 2, the third integer 6 then represents the following number of succesive switches of 1 and 3 and so forth.
Alternately the code notation  can be gotten from the side sequence notation by successively counting the number of 3's then the number of 2's then the number of 3's and so forth and then multiplying by two to get the code sequence. 

\doublespace
\singlespace
\noindent
\textbf{Important comment:} The sequence of code numbers divided by two represents the successive number of black and blue points (or blue and black points) as they alternate in the tower. For the example above there are 2 black points followed by 3 blue points, 3 black points and so forth.

\doublespace
\singlespace
Convention: If we use exponential code and if the exponent is 0, this means that the code number does not appear. For example 4 $6^0$ $4^2$ $2^4$ is just $4^3$ $2^4$ in exponential code. If the base is zero this also means the code number does not appear. For example 0 $6^0$ $0^2$ $2^4$ is just $2^4$.

\doublespace
\singlespace
Observe that if the side sequence starts 13 then given its code sequence, we can recover the original side sequence by starting the sequence with 13... If we choose to recover the side sequence by starting it with 12... then we will recover it with the 3's and 2's interchanged. This is not a problem as it just represents a relabeling of the sides of triangle ABC. Also observe that since a rhombus poolshot starts and ends with AB then its code sequence consists of all even integers. As with side sequences we can also use a code sequence to describe a tower of mirror images of triangle ABC be it a rhombus tower or not.

We will make the convention that if there are three dots in front of a code sequence (or following a code sequence) then this means there is at least one \textbf{code number} preceding (or following) that code sequence in which case we will call it a \textbf{subcode} For example ... 2 4 ... is a subcode of the code sequence 2 2 4 4. Corresponding to any subcode, there is a sequence of mirror images of triangle ABC which is a tower in its own right and which we will call a \textbf{subtower}.

\subsection*{5. Corridor Towers}

Now take a rhombus tower of mirror images of triangle ABC that starts with half a rhombus (the isosceles triangle ABC) and finishes with half a rhombus and which has $n-1$ rhombi in between and which involves 2n copies of triangle ABC. In order to form the corridors as in Galperin's paper, the following three conditions must also be satisfied.

 1.  The first reflection is in side $A_{1}$$C_{1}$=$A_{0}$$C_{1}$ and the last reflection is in side $A_{n}$$C_{n}$.

 2. $A_{n}$$B_{n}$ must be parallel to $A_{0}$$B_{0}$ and that there must exist a \textbf{special integer} "m" with $0<m<n$ such that $A_{m}$$B_{m}$ is also parallel to $A_{0}$$B_{0}$ and such that the (2m-1)st reflection is in side $C_{m}$$B_{m}$ and the (2m+1)st reflection is in side $C_{m+1}$$B_{m+1}$.

We can then form corridors as follows. The \textbf{first or right corridor} is formed by taking the line $C_{1}$$C_{n}$ and extending it to hit $A_{0}$$B_{0}$ at $S_{0}$, $A_{m}$$B_{m}$ at $L_{m}$ and $A_{n}$$B_{n}$ at $M_{n}$. The line $S_{0}$$M_{n}$ then forms a \textbf{(right) boundary line of the first corridor}.
The other \textbf{(left) boundary line of the first corridor} $P_{0}$$Q_{n}$ is formed by introducing a parallel line through $C_{m}$ and extending it to hit $A_{0}$$B_{0}$ at $P_{0}$, $A_{m}$$B_{m}$ at $M_{m}$ and $A_{n}$$B_{n}$ at  $Q_{n}$. Now introduce a third parallel through $C_{m+1}$ hitting $A_{m}$$B_{m}$ at $S_{m}$ and $A_{n}$$B_{n}$ at $P_{n}$ and form the \textbf{second wider corridor} with boundary  lines $S_{m}$$P_{n}$ and $L_{m}$$M_{n}$. We will also call the corridor with boundary lines $S_{m}$$P_{n}$ and $M_{m}$$Q_{n}$ the \textbf{left corridor}.

3. Finally the first and second corridors must have positive width with the width of the second corridor larger than the first and the corridors must go through the interiors of the sequence of mirror images of triangle ABC and have no vertices between their boundary lines although there can be and are vertices on their boundary lines.

If a rhombus tower satisfies all three conditions, we will call it a \textbf{corridor rhombus tower} or \textbf{corridor tower}.

It is worth noting that any \textbf{associated}  rhombus poolshot that creates the rhombus tower need not be parallel to the boundary lines or sides of the corridors if they exist but can always be taken to be so and henceforth we will assume that is the case. Also note that if the corridor tower is placed in standard position with base $A_{0}$$B_{0}$  horizontal and the straightened pool shot going upwards, then it must have positive slope (leans to the right) otherwise we would not have been able to form the second corridor. 

Also observe that $C_{1}$ and $C_{n}$ are black points while the \textbf{special C points} $C_{m}$ and $C_{m+1}$ are blue points and if the corridor tower is in standard position then all the C points to the right of the associated (parallel) rhombus poolshot are black points and all the C points to the left are blue points.
The original Galperin example composed of 44 copies of triangle ABC can be described by the side sequence $131212121313131312121212131313121213131212131$ or 2 6 8 8 6 4 4 4 2 in code. Observe that the sum of the code
numbers gives us the number of triangles involved. The jth integer in a code  sequence is called the \textbf{jth code number}.
Below Fig. 6 is another example of a corridor tower 131212121313131212131 or 2 6 6 4 2 composed of 20 triangles which will turn out to be the shortest corridor tower and in which the ratio of the widths is also 2 to 1 for every angle x for which the corridors can be formed. This will eventually be proved in the Classification Theorem.

\doublespace
Figure 6
\singlespace 

The subtower containing the corridor from the segment $P_{0}$$S_{0}$ to the segment $M_{m}$$L_{m}$ where $m$ is the \textbf{special integer} will be called the \textbf{first level} and subtower containing the wider corridor from the segment $S_{m}$$L_{m}$ to the segment $P_{n}$$M_{n}$ will be called the \textbf{second level}. Note that if we consider each level separately then each level starts with AB and ends with AB and is a tower in its own right and we can use our code notation to describe each level. The first level is a tower whose first reflection is in side AC and last reflection is in side BC while the second level is a tower whose first reflection is in side BC and last reflection is in side AC.

\doublespace
Figure 7a and 7b
\singlespace

\noindent \textbf{Fact 1:} If a  corridor tower or rhombus tower has "2n" in its code where $x=<A=<B$, then $0<x<90/n$. 
\singlespace
\noindent Proof: Since then the straightened pool shot  must cross an angle of  2nx which means that $2nx <180$. QED

\doublespace
Figure 8
\singlespace

Now observe that the first level of a corridor tower is \textbf{antisymmetric} (which here means that if we look at its side sequence, and interchange the 2's and the 3's and reverse its order then we get exactly the same side sequence) since the ray $M_{m}$$C_{m}$ leaves $A_{m}$$B_{m}$ at the same angle as the ray $S_{0}$$C_{1}$ leaves $A_{0}$$B_{0}$. Similarily the second level is antisymmetric. As an example the first level of the original corridor is 131212121313131312121212131313121
 (2 6 8 8 6 2 in code) and its second level is 1213131212131 (2 4 4 2 in code) both of which read the same forwards and                                                                    
backwards if we interchange the 2's and 3's. It follows that the code sequence of the first level and of the second level are both \textbf {symmetric} and each has an \textbf{even} number of even code numbers.

\doublespace
\singlespace
\noindent \textbf{Fact 2:} $A_{i}$$B_{j}$ is parallel to $A_{k}$$B_{l}$ if there is the same number of blue points and black points between the two line segments. In particular this is the case if the sequence of code numbers between the two line segments is even and symmetric.

\subsection*{6. Algorithm One}

\singlespace
There is a nice algorithm that can be used to calculate the coordinates of the C points in a corridor tower (or more generally any arbitrary tower) in \textbf{standard position} assuming the coordinates of the first C point are (0,1) and that the distance between successive C points is exactly one unit.

Starting at the base (which has Y coordinate $\frac{1}{2}$) of the corridor tower and proceeding upwards count the number of black points then the number of blue points and continue alternating in this manner. As an example for the original corridor this gives us the sequence 1 3 4 4 3 2 2 2 1 which is exactly half of its code sequence 2 6 8 8 6 4 4 4 2.  We can now form the corresponding \textbf{vertical array} of integers  as follows and illustrated below.

1. The first line in the vertical array which we call a \textbf{black line} because it consists wholly of black points  is the integer 0 which represents the first black point $C_1$ and we assign it the coordinates (0,1) or (sin0, cos0)

2. The color of the points in a line alternately changes from black to blue or from blue to black as we go from one line in the array to the next line. For black lines consecutive integers decrease by 2 (and for blue lines they increase by 2)  throughout that line and this continues on to the first integer on the next line after which it switches from increasing to decreasing or vice versa

3. The jth code number divided by 2 indicates the number of integers on the jth line which we call its \textbf{length} which for this example gives us the successive lengths 1 3 4 4 3 2 2 2 1

4. The ith integer in this array reading from up to down and from left to right represents the point $C_i$ and has coordinates  (X,Y) where X is a sum of the form ${a_1sin2x+a_2sin4x+...+a_tsin2tx+...+a_ksin2kx}$ = 

\noindent ${a_0sin0x+a_1sin2x+a_2sin4x+...+a_tsin2tx+...+a_ksin2kx}$ where the coefficient of sin2tx is calculated by counting up the number of times the integer 2t appears with a plus sign minus the number of times it appears with a negative sign in the array up to and including the integer representing $C_i$. Similarily Y is a sum of the form  ${a_0+a_1cos2x+a_2cos4x+...+a_tcos2tx+...+a_kcos2kx}$ = 

\noindent ${a_0cos0x+a_1cos2x+a_2cos4x+...+a_tcos2tx+...+a_kcos2kx}$ where the coefficient of cos2tx is calculated by counting up the number of times the integer 2t appears in the array irregardless of its sign up to and including the integer representing $C_i$.

\doublespace

\noindent In our example we get the following vertical array.

\singlespace
\underline{0*}  black line

\underline {-2 0 2} blue line

4 2 0 -2 black line

\underline {-4 -2 0 2} blue line

\underline {4 2 0} black line

-2* 0* blue line

2 0 black line

\underline{-2 0} blue line

2* black line

\doublespace

\singlespace

\textbf{Note 1:} We have starred the integers corresponding to $C_{1}$,$C_{16}$,$C_{17}$ and $C_{22}$ and we add a horizontal line whenever a line changes its length or whenever its length is unknown. We used this algorithm to calculate the coordinates of these points at the beginning of this paper. For example the coordinates of $C_{22}$ were found to be $(sin2x+sin4x,8+11cos2x+3cos4x)$.

\textbf{Note 2:} A vector between any two C points can easily be found by looking at only the integers in the array between the two C points (excluding the first C point and including the last C point) and by counting as in step four above. For example the vector from the third starred point to the fourth is given by (sin2x,2+3cos2x).

\textbf{Note 3:} For any \textbf{arbitrary tower}, the rules are the same as enumerated above with modifications as set out below. 

a. If it is in standard position the first integer 0 in the array can be a black or a blue point.

b. We may also choose to start the vertical array with any even integer 2n, positive or negative in which case the coordinates of the first C point which can be black or blue will be considered to be set at (sin2nx, cos2nx) and this would correspond to a tower which is not in standard position if 2nx is not a multiple of 360. Below is an example of the tower above where the first C point has coordinates (-sin2x, cos2x). A   vector from the third starred point to the fourth is now given by (-2sin2x-sin4x,2+2cos2x+cos4x).

\doublespace

\singlespace
\underline{-2*}  black line

\underline {-4 -2 0} blue line

2 0 -2 -4 black line

\underline {-6 -4 -2 0} 

\underline {2 0 -2} 

-4* -2* 

0 -2

\underline{-4 -2} blue line

0* black line

\doublespace

\singlespace

c. We can also recoordinatize the C points by subtracting a constant vector (u,v) from the coordinates of every C point. In this way the first C point can be considered to start with any pair of coordinates. It is important to observe however that this will not change any relationships between the vectors for example parallel vectors will still be parallel.

\textbf{Note 4:} Given isosceles triangle ABC and a fixed angle x then there is a one to one correspondence between code sequences of even positive integers and towers of triangle ABC whose first reflection is in side AC (and similarly if its  first reflection is in side BC). In other words the code sequence determines an unique tower and the code sequence can be taken as its name. Further the code sequence determines an unique vertical array once the first even integer and its color has been chosen.

\subsection*{7. Tests}

We need to develop a sequence of useful tests and start with the following.

\textbf{Important Fact:} If a vertical poolshot enters a rhombus A$C_{1}$B$C_{2}$ say first at A$C_{1}$ and doesn't hit a vertex and crosses the diagonal AB then it leaves the rhombus at either A$C_{2}$ or B$C_{2}$ such that  the Y coordinate of vertex $C_{2}$ is greater than the Y coordinate of  vertex $C_{1}$.

Proof: Assume that the vertical  pool shot goes through A$C_{1}$ as shown on the diagram and that AB has positive slope $m_{1}$. Now since the diagonals are perpendicular the slope of $C_{1}$$C_{2}$ is the negative reciprocal and the Y coordinate must increase from $C_{1}$ to $C_{2}$. Similarily if AB has negative or zero slope. QED

\doublespace
Figure 9
\singlespace

\textbf{Consequence 1:} Given a straightened vertical  rhombus pool shot involving isosceles triangle ABC, then the Y coordinates of the "C" points keep increasing as the pool shot moves upwards.

Note if the rhombus poolshot as above is not vertical and we drop perpendiculars from the "C" points onto the straightened poolshot viewed as a positive coordinate axis, then the location  of the feet of the perpendiculars (which we will still call the \textbf{Y coordinate of "C"}) increase as the pool shot moves forward. See Fig. 10. Also note that the pool shot can be considered as a vector and as it goes through the sequence of mirror images of triangle ABC, it induces an increasing order on them and also an increasing order on the "C" points. This means that the "C" points are  in fact labelled such that the Y coordinate of $C_{i}$ is less than the Y coordinate of $C_{j}$ if and only if $i$ is less than $j$. We will refer to $C_{i}$ as the \textbf{lower} C point and $C_{j}$ as the \textbf{higher} C point.

Caution: Given a tower in standard position which is not a rhombus tower and with the usual ordering of the C points in the same order that they are formed by the  sequence of reflections, the Y coordinates need not increase as the index $i$ in $C_{i}$ increases.

\doublespace
Figure 10
\singlespace

 Now observe that if we are given the two corridors in a corridor rhombus tower which we can assume is in standard position and if we form a vector v=(a,b) from the special blue point $C_{m+1}$ to a \textbf{higher} blue C point and a vector w=(c,d) from the  special blue point $C_{m}$ to a \textbf{higher} black C point (note this means that b and d are both positive), then the two vectors cannot be parallel (see Fig. 11).  Now the two vectors are parallel if w=kv where k=d/b and hence ka=c which is equivalent to ad=bc.

\doublespace
Figure 11
\singlespace

This leads to the following test.

\textbf{Non-Corridor Test 1:} Given a rhombus tower satisfying the first two corridor conditions and a vector v=(a,b) from the special blue point $C_{m+1}$ to a higher blue C point which is parallel (ad=bc) to a vector w=(c,d) from the  special blue point $C_{m}$ to a higher black C point then the two corridors cannot be formed with all three requisite properties and it is not a corridor rhombus tower.

\textbf{Consequence 2:} Given a corridor rhombus tower and if we join the first black point $C_1$ (or the last black point $C_n$) to any black point then no blue point can be collinear with those two in any order. It follows that if we join $C_1$ (or $C_n$) to any blue point, no black point can be collinear with these two.

Proof: Since the black and blue points are separated by the right corridor and since $C_{1}$ has the smallest Y coordinate and $C_{n}$ has the largest and the Y coordinates of the C points increase as its subscript increases. QED

\doublespace
Figure 12
\singlespace

This gives the following test.

\textbf{Non-Corridor Test 2:} Given a prospective corridor rhombus tower satisfying the first condition in which there is a blue point collinear in any order with either the first black point $C_1$ (or the last black point $C_n$) and another black point, then the corridors cannot be formed with all three properties and it is not a corridor rhombus tower. 

\textbf{Consequence 3:} Given a rhombus tower where $C_{i}$ is a blue point, $C_{j}$ is a black point and $C_{k}$ is a blue point with $i<j<k$  then the pool shot must go through the segment $C_{i}$$C_{j}$ first and then go through the segment $C_{j}$$C_{k}$ second. 

Proof: Since the Y coordinates of the C points increase as the index increases, it follows that the Y coordinate of every point on $C_{i}$$C_{j}$ is less than the Y coordinate of every point on $C_{j}$$C_{k}$ and the result follows. QED

\doublespace
Figure 13
\singlespace

This means that the orientation of the triangle $C_{i}$$C_{j}$$C_{k}$ in that order must be counterclockwise. 
Similarly if $C_{i}$ is a black point, $C_{j}$ is a blue point and $C_{k}$ is a black point with $i<j<k$  in which case the orientation of $C_{i}$$C_{j}$$C_{k}$ is clockwise.

\textbf{Non-Rhombus Tower Test 3:} Given a tower where $C_{i}$ is a blue point, $C_{j}$ is a black point and $C_{k}$ is a blue point with $i<j<k$ and if the orientation of triangle $C_{i}$$C_{j}$$C_{k}$ is clockwise or the three points are collinear, then there is no rhombus poolshot and hence it is not a rhombus tower which further means it is not a corridor rhombus tower.

Similarly if $C_{i}$ is a black point, $C_{j}$ is a blue point and $C_{k}$ is a black point with $i<j<k$  and the orientation of $C_{i}$$C_{j}$$C_{k}$ is counterclockwise or the three points are collinear then it is not a rhombus tower.

\doublespace
\singlespace
This is an important test for the non-existance of the corridors and we can rephrase it using the well known algorithm that if the coordinates of $C_{i}$ are ($x_{0}$,$y_{0}$), of $C_{j}$ are ($x_{1}$,$y_{1}$) and of $C_{k}$ are ($x_{2}$,$y_{2}$) then the orientation of triangle $C_{i}$$C_{j}$$C_{k}$ is counterclockwise if ($x_{1}$-$x_{0}$)($y_{2}$-$y_{0}$)-($x_{2}$-$x_{0}$)($y_{1}$-$y_{0}$) is positive and clockwise if negative. If the difference is zero then the three points are collinear. Writing this in vector form we get the following.

\textbf{Non-Rhombus Tower  Test 4:} Given a tower and

a. if we form a vector v=(a,b) from any blue point C to a higher black point C' and another vector w=(c,d) from C' to a higher blue point C" and if the orientation CC'C'' is clockwise or the points are collinear then there is no rhombus poolshot.
Note this is the case if $a(b+d) \leq b(a+c)$ which is equivalent to $ad \leq bc$.

b. alternately if we form a vector v=(a,b) from C to C' and another vector u=(a+c,b+d)  from C to C"  and  if the orientation CC'C'' is clockwise or the points are collinear then there is no rhombus poolshot if $ad \leq bc$.

In either case it is not a rhombus or corridor rhombus tower if $ad \leq bc$

c. Similarly if we form a vector w=(c,d) from any black point C to a higher blue point C' and another vector v=(a,b) from C' to a higher black point C" and if the orientation CC'C'' is counterclockwise or the points are collinear  which is then still equivalent to $ad\leq bc$ then there is no rhombus poolshot.

\doublespace
Figure 14
\singlespace

Indeed Test 4 can be extended further to get

\textbf{Non-Rhombus Tower  Test 5:} Given a tower and if we take a vector v=(a,b) from a blue point C to a higher black point C' and a vector w=(c,d) from a possibly different black point C'' to a higher blue point C''', and if the orientation of the triangle with coordinates (0,0), (a,b), (a+c,b+d) is clockwise (or the triangle with coordinates (0,0), (c,d), (a+c,b+d) is counterclockwise) then there is no rhombus poolshot. In either case there is no rhombus poolshot if $ad < bc$ which means the tower is not a rhombus or corridor rhombus tower. 

\doublespace
\singlespace
The version of the above where the points are collinear is treated separately below and can be extended to allow them to form parallel lines in the following way.

\doublespace
\singlespace
\textbf{Non-Rhombus Tower  Test 6:} Given a code sequence of a tower and its corresponding vertical array which produces two \textbf{parallel} vectors  pointing in the increasing direction such that one v=(a,b) is from a blue (black) point $C_{i}$ to a higher black (blue) point $C_{j}$ and the other w=(c,d) is from a black (blue) point $C_{k}$ to a higher blue (black) point $C_{p}$, then there is no rhombus pool shot and hence the tower is not a rhombus or corridor rhombus tower. Note the vectors are parallel if $ad = bc$.

Proof: As the pool shot would have to go through one of them first say $C_{i}$$C_{j}$ and then $C_{k}$$C_{p}$ second which would cause the two blue C's (and the two black C's) to be on opposite sides of the straightened poolshot which is impossible. QED

\doublespace
Figure 15
\singlespace

Combining the two tests, we get 

\textbf{Non-Rhombus Tower  Test 7:} Given a tower and if we take a vector v=(a,b) from a blue point C to a higher black point C' and a vector w=(c,d) from a possibly different black point C'' to a higher blue point C''', and if $ad \leq bc$ then there is no rhombus poolshot which means the tower is not a rhombus or corridor rhombus tower. 

\doublespace
\singlespace

As an example we can use test 6 to prove that the subcode $6^2$ 4 2 $4^2$... cannot occur in the code sequence of a rhombus poolshot.

Proof: Consider the two starred points, the first of which is  blue and the second one black and the two double starred points the first of which is black and the second one blue. Observe that since the integers between them are exactly the same (excluding the initial points and including the terminal points), the corresponding vectors must be parallel and also satisfy the conditions set out in test 6. Hence no rhombus poolshot contains this code sequence. Observe that this array starts at -2.

\doublespace
\singlespace

-2* 0 2  blue

\underline {4 2 0*}  black

\underline{-2 0}

\underline 2**  black

0 2

\underline {4 2}

0**... blue

\doublespace

\noindent QED

\singlespace

\textbf{Non-Corridor Test 8:} Given a rhombus tower satisfying the first two corridor conditions and a vector v=(a,b) from a lower blue point to the special blue point $C_{m}$ and  a vector w=(c,d) from a lower black point to the last black point $C_{n}$   then the two corridors cannot be formed with all three requisite properties and it is not a corridor rhombus tower if $ad<bc$. Observe that if we use the vectors v and w or the vectors -v and -w, the test is the same.

Proof: Since the triangle with coordinates P(0,0), Q(a,b), R(a+c,b+d) cannot have clockwise orientation from P to Q to R which is the case if $ad<bc$.

\doublespace
\noindent QED

\doublespace
Figure 16

\singlespace
\textbf{Non-Corridor Test 9:} Given a rhombus tower satisfying the first corridor condition and a vector v=(a,b) from any black point to a lower blue point  and  a vector w=(c,d) from the last black point $C_{n}$ to a lower black point     then the two corridors cannot be formed with all three requisite properties and it is not a corridor rhombus tower if $ad\leq bc$.  Similarly if we use  a vector  v=(a,b) from the first black point $C_{1}$ to a higher black point and a vector  w=(c,d) from any black point to a higher blue point then it is not a corridor rhombus tower if $ad\leq bc$. Again if we use the vectors v and w or the vectors -v and -w, the test is the same.

Proof: Since the points with coordinates P(0,0), Q(a,b), R(a+c,b+d) cannot have clockwise orientation from P to Q to R or be collinear which is the case if $ad\leq bc$.

\doublespace
\noindent QED

\doublespace
Figure 17

\singlespace
\noindent Finally observe that if in a corridor rhombus tower we take a vector v=(a,b) from a lower blue point to the special blue point $C_{m}$ and  a vector w=(c,d) from a lower black point to the last black point $C_{n}$ and if these two vectors are parallel then the line through $C_{m}$ with vector v and the line through $C_{n}$ with vector w are in fact along the boundary lines of the first corridor. Hence we get the following test.

\doublespace
\singlespace
\noindent
\textbf{Locked In Test:} Given a fixed rhombus tower satisfying the first two corridor conditions and a vector v=(a,b) from a lower blue point $C_{k}$ to the special blue point $C_{m}$ and  a vector w=(c,d) from a lower black point $C_{k'}$ to the last black point $C_{n}$ and if v and w are parallel (ad=bc) then all corridor rhombus towers ending with this same fixed rhombus tower all have the line through $C_{k}$ and $C_{m}$ and the line through $C_{k'}$ and $C_{n}$ as the boundary lines of the first corridor. We can equally use v and w or -v and -w.

Note this means that in all such corridor rhombus towers, the ratio of the widths of the two corridors is exactly the same.

\doublespace
Figure 18

\subsection*{8. Algorithm Two}

\singlespace
To apply these tests is going to require that we check various trig identities or trig inequalities many of which will involve certain kinds of patterns of finite sums of sines and cosines. To handle these we first reduce them to a shorter sum by the following algorithm.

\doublespace

\noindent Given the sum $a_{0}$+$a_{1}$cos2x+$a_{2}$cos4x+$a_{3}$cos6x+...+$a_{k}$cos2kx

\singlespace

{\setlength\parindent{0pt} 
 Step 1: Multiply by 2sinx and use the trig identity $2sinxcos2nx=sin(2n+1)x-sin(2n-1)x$. This results in the sum

\noindent (2$a_{0}$-$a_{1})$sinx+($a_{1}$-$a_{2}$)sin3x+($a_{2}$-$a_{3}$)sin5x+ ... +($a_{k-1}$-$a_{k}$)sin(2k-1)x+$a_{k}$sin(2k+1)x
which is of the form $c_{1}$sinx+$c_{2}$sin3x+$c_{3}$sin5x+...+$c_{k}$sin(2k+1)x}.

\doublespace
\singlespace

{\setlength\parindent{0pt} Step 2: Multiply by 2sinx again and use the trig identity $2sinxsin(2n+1)x=cos(2n)x-cos(2n+2)x$ to get

\noindent $c_{1}$+($c_{2}$-$c_{1}$)cos2x+($c_{3}$-$c_{2}$)cos4x+...+($c_{k}$-$c_{k-1}$)cos(2k)x-$c_{k}$cos(2k+2)x}.

\doublespace
\singlespace
\noindent To illustrate with an example if we start with the sum of k+1 terms in the pattern below

{\setlength\parindent{0pt} $(2nk-2n+1)+(4nk-4n+2)cos2x+(4nk-6n+1)cos4x+(4nk-10n+1)cos6x+...+(6n+1)cos(2k-2)x+(2n+1)cos2kx$ 

it becomes at step one

(2n+1)sin3x+ 4nsin5x+4nsin7x+...+4nsin(2k-1)x+(2n+1)sin(2k+1)x 

and then becomes at step two just a sum of four terms

(2n+1)cos2x+(2n-1)cos4x+(-2n+1)cos2kx+(-2n-1)cos(2k+2)x} 

\noindent which is exactly the original sum times $4sin^2x$.

\doublespace
\noindent Similarily given the sum $d_{1}sin2x+d_{2}sin4x+d_{3}sin6x+...+d_{k}sin2kx$

\doublespace
\singlespace

{\setlength\parindent{0pt} Step 1: Multiply by 2sinx and use the trig identity $2sinxsin2nx=cos(2n-1)x-cos(2n+1)x$ to get

\noindent $d_{1}cosx+(d_{2}-d_{1})cos3x+(d_{3}-d_{2})cos5x+...+((d_{k}-d_{k-1}))cos(2k-1)x-d_{k}cos(2k+1)x$
which is of the form

\noindent $e_{1}cosx+e_{2}cos3x+e_{3}cos5x+...+e_{k}cos(2k-1)x+e_{k+1}cos(2k+1)x$}.

\doublespace
\singlespace

{\setlength\parindent{0pt} Step 2: Multiply be 2sinx again and use the trig identity $2sinxcos(2n+1)x=sin(2n+2)x-sin2nx$ to get

\noindent $(e_{1}-e_{2})sin2x+(e_{2}-e_{3})sin4x+...+(e_{k}-e_{k+1})sin(2k)x+e_{k+1}sin(2k+2)x$}.

\doublespace
\noindent Example:
\singlespace

{\setlength\parindent{0pt} $(2nk-2n+1) + (4nk-4n+2)sin2x  + (4nk-6n+1)sin4x  + (4nk-10n+1)sin6x + ...  +(6n+1)sin(2k-2)x  + (2n+1)sin2kx$ 

becomes at step one

$(2nk -2n + 1)cosx+(-2n -1)cos3x+(-4n)cos5x+(-4n)cos7x+ ... +(-4n)cos(2k-1)x+(-2n -1)cos(2k+1)x$ 

and then becomes at step two just a sum of four terms

$(2nk+2)sin2x+(2n -1)sin4x+(-2n+1)sin2kx+(-2n -1)sin(2k+2)x$} which is exactly the original sum times $4sin^2x$.

\section*{9. Useful Trig Identities}

\textbf{The Main Trig Identity:} $16sin^4x$$[sin2kx+sin(2k+2)x+ ... +sin(2k+2s)x]=sin(2k-4)x-3sin(2k-2)x+3sin2kx-sin(2k+2)x-sin(2s+2k-2)x+3sin(2s+2k)x-3sin(2s+2k+2)x+sin(2s+2k+4)x$

\doublespace
\singlespace

\noindent which can be proved using the identities $2sinxsin2nx=cos(2n-1)x-cos(2n+1)x$ and $2sinxcos(2n+1)x=sin(2n+2)x-sin2nx$.

\doublespace
\singlespace
\noindent Special Cases

1. s=0 then $16sin^4x$$sin2kx=sin(2k-4)x-4sin(2k-2)x+6sin2kx-4sin(2k+2)x+sin(2k+4)x$

2. s=1 then $16sin^4x$$[sin2kx+sin(2k+2)x]=sin(2k-4)x-3sin(2k-2)x+2sin2kx+2sin(2k+2)x-3sin(2k+4)x+sin(2k+6)x$

3. s=2 then $16sin^4x$$[sin2kx+sin(2k+2)x+sin(2k+4)x]=sin(2k-4)x-3sin(2k-2)x+3sin2kx-2sin(2k+2)x+3sin(2k+4)x-3sin(2k+6)x+sin(2k+8)x$

If $sinx$$\not= 0$, another useful identity is 

\noindent
$sinz+sin(z+2x)+sin(z+4x)+...+sin(z+2nx)=sin((n+1)x)sin(z+nx)/sinx $

As an example if we let z=4x and replace n by 2n this then becomes
$sin4x+sin6x+...+sin(4n+2)x+sin(4n+4)x=sin((2n+1)x)sin(2n+4)x/sinx$ This would be greater than zero if $0<(2n+4)x < 180$ and $n\geq0$.

\section*{10. Classification of all Corridor Towers}

\singlespace

In order to show that Galperin's method can never work, we must first classify all corridor towers. This will take some work. We start by first characterizing the second level of a  corridor tower which is a rhombus tower all by itself.

Recall as previously noted that \textbf{any corridor in standard position} must lean to the right as otherwise we would not be able to form the second level. This means that the X coordinate of any black point (other than the first which is set at zero) must be greater than zero. But this entails that the code sequence must start at 2 6 or higher. If it started 2 2..., the second black point would have X coordinate -sin2x a negative value (keeping in mind that x is an acute angle between 0 and 90) and if it started 2 4... the second black point would have X coordinate 0. This means the corridor must start 131212121... and then by antisymmetry the first level must end ...131313121.

\doublespace
Figure 19
\singlespace

We can now use this to help us show that the second level is always of the code sequence form 2 $4^{2k}$ 2 for some integer $k\geq0$. We also make use of  the following lemmas. Any rule prefaced  by the word corridor indicates that it is about a corridor tower.

\textbf{Corridor Lemma 1:} The second level must start 1213... (starts 2 ... in code) and by antisymmetry the second level ends in ...2131 (ends ... 2 in code).

Proof: Since  the last reflection by definition is in side AC  the second level ends ...31  and by antisymmetry must then start 12.... This means the second level starts 12$(12)^r$13 for some integer $r\geq0$. If $r>0$ then the line through $C_{m+1}$ and $C_{m+r+1}$ which are blue points is parallel to the line through the blue point $C_{m}$ and the black point $C_{m+r+2}$ since both are perpendicular to the bisector of the angle $C_{m}$$B_{m}$$C_{m+r+2}$ of size 2(r+2)x at $B_{m}$ Hence the corridor can't exist by non-corridor test 1 and the second level must start 1213.... QED

\doublespace
Figure 20
\singlespace

\textbf{Cor 1:} This means the right or first corridor as it passes from the first level to the second level is of the form ...131313121213...
( ... 4 ... in code)

The next rule is a rule about rhombus poolshots. Any rule prefaced  by the word rhombus indicates that it is about a rhombus poolshot or rhombus tower.

\doublespace
\singlespace
\textbf{Rhombus Rule A:} The subcode  ...2n 2n+2k where 2n is not the first code number (or 2n+2k 2n... where 2n is not the last code number) $n\geq1$, $k\geq2$ never appears in the code sequence of any rhombus pool shot. Caution: This rule does not hold if 2n is the first or last code number as can be seen from the existance of the corridor 2 6 6 4 2 in Fig. 6.

Proof: Suppose it does then since we are free to orient the corresponding sequence of mirror images of triangle ABC however we want and since there is a code number preceding 2n, we can start with triangle ABC in standard position and then follow that by a sequence of mirror images corresponding to the code 2n 2n+2k. Now observe that the three "C" points which we have starred below are collinear (shown below) and since this is a black-blue-black collinear situation, it is impossible by the Non-Rhombus Tower Collinear Test. Note that the second starred point lies between the other two since the Y coordinates of the C points are increasing.

\doublespace
\singlespace
\underline {..0*}  black

\underline{-2 0 2 4 ... 2n-6 (2n-4)*}   blue

2n-2 2n-4 2n-6 ... 2 0 -2 (-4)* ...  black

\doublespace

\singlespace

Using this vertical array, a vector v=(a,b) from the first black starred point to the second starred blue point is $(sin4x+sin6x+...+sin(2n-4)x,1+2cos2x+cos4x+cos6x+...+cos(2n-4)x)$ and a vector w=(c,d) from second starred blue point to the third starred black point is
 
\noindent $(sin6x+...+sin(2n-4)x+sin(2n-2)x,1+2cos2x+2cos4x+cos6x+...+cos(2n-4)x+cos(2n-2)x)$ and the three points are collinear since ad=bc.

There are special cases which we will do here but for other rules we may leave it to the reader.

\noindent If n=1, v=$(-sin2x,cos2x)$ and w=$(-sin2x-sin4x,1+cos2x+cos4x)$ which are collinear since ad=bc.

\noindent If n=2, v=$(-sin2x,1+cos2x)$ and w=$(-sin4x,1+2cos2x+cos4x)$ which are collinear since ad=bc.

\noindent If n=3, v=$(0,1+2cos2x)$ and w=$(0,1+2cos2x+2cos4x)$ which are collinear since ad=bc or equivalently since the three points lie on the vertical line x=0.

If n greater than 3, we can use the general case. QED

\doublespace
Figure 21a and 21b
\singlespace

This means that if a code of a rhombus tower contains the integer 2n other than at the start or at the end, then any adjacent integer  must be either 2n-2, 2n or 2n+2. In other words it can change by no more than 2 if at all. We will also say that other than at the end or the start, 2n \textbf{forces} the next or preceding code number in a rhombus tower to differ by no more than two.

\textbf{Cor 2:} It follows that the right corridor as it passes from the first level to the second level is of the form ...21313131212131... (or in code by ... 6 4 ...) which means that the first level ends ...2131313121 (... 6 2 in code) and hence by antisymmetry the first level  starts exactly 131212121... or (2 6 ... in code).

Proof: Since by Cor 1, it is of the form ...131313121213... and can't be of the form ...2k 4...  for $k\geq4$ by Rhombus rule A.  QED

\doublespace
\singlespace

\textbf{Rhombus Rule B:} The subcodes

\noindent a. ... 2n $(2n+2)^{2k+1}$ 2n ... with $n\geq1$, $k\geq0$ where 2n is not the first or the last code number never appears in the code sequence of any rhombus pool shot and  

\noindent b. 2n $(2n-2)^{2k+1}$ 2n with $n\geq2$, $k\geq0$ never appears in the code sequence of any rhombus pool shot. Note for part b there need not be any other  code numbers.

Proof: 

a. The black-blue vector between the two starred points is parallel to the blue-black vector between the double starred points since the integers between them (excluding the first and including the last) are the same and hence by the Non-Rhombus Tower Test 6 there is no poolshot.

\doublespace

\singlespace

\noindent
\underline {..0*} black

\noindent
\underline{-2 0 2 4 ... 2n-6 (2n-4)*} blue

\noindent
$\left.
\begin{aligned}
&\text{2n-2 2n-4 2n-6 ... 0 -2} \\
&\text{-4 -2 0 2 ... 2n-6 2n-4} \\
&\vdots\quad\vdots\\
&\underline {\text{2n-2 2n-4 2n-6 ... 0 -2  }} \\
\end{aligned} 
\;\right\}$ 2k+1 lines

\noindent
\underline{-4** -2 0 2 ... 2n-8 2n-6}  blue

\noindent
(2n-4)** ...  black

\doublespace
\singlespace

\doublespace
Figure 22
\singlespace

b. The blue-black vector between the two starred points is parallel to the black-blue vector between the double starred points since the integers between them (excluding the first and including the last) are the same and hence by the Non-Rhombus Tower Test 6 there is no poolshot. Note we have oriented the sequence of mirror images so that the first starred integer is a blue -2*.

\doublespace

\singlespace

\noindent
\underline{-2* 0 2 4 ... 2n-6 2n-4} 

\noindent
$\left.
\begin{aligned}
&\text{(2n-2)* 2n-4 ... 4 2} \\ 
&\text{0 2 ... 2n-6 2n-4} \\
&\vdots\quad\vdots\\
&\underline {\text{2n-2 2n-4 ... 4 2**  }} \\
\end{aligned} 
\;\right\}$ 2k+1 lines

\noindent
0 2 ... 2n-4 (2n-2)**

\doublespace
\singlespace

\doublespace
Figure 23
\singlespace

\textbf{COR 3:} It follows that the right corridor as it passes from the first level to the second level can't be of the form ... 6 $4^{2k+1}$ 6 ... with $k\geq0$.

\noindent Proof: Use Rhombus rule B, part b with n=3.
QED

\doublespace
\singlespace

\textbf{Rhombus Rule C:} The subcode ... 2n $(2n+2)^{2k}$ 2n+4 for $n\geq1$, $k\geq0$ never appears in the code of any rhombus pool shot where 2n is not the first code number. Similarily for  2n+4 $(2n+2)^{2k}$ 2n ... where $n\geq1$, $k\geq0$ where 2n is not the last code number.

Proof: Note this is true for k=0 by Rhombus rule A so we can assume that $k>0$. Now the black-blue vector between the two starred points below is parallel to the blue-black vector between the double starred points since the integers between them (excluding the first and including the last) are the same and hence by the Non-Rhombus Tower Test 6 there is no poolshot.

\doublespace

\singlespace

\noindent
\underline {..0*} black

\noindent
\underline{-2 0 2 4 ... 2n-6 (2n-4)} blue

\noindent
$\left.
\begin{aligned}
&\text{2n-2 2n-4 2n-6 ... 0 -2} \\
&\text{-4* -2 0 2 ... 2n-6 2n-4} \\
&\vdots\quad\vdots\\
&\underline {\text{-4** -2 0 2 ... 2n-6 2n-4  }} \\
\end{aligned} 
\;\right\}$ 2k lines

\noindent
2n-2 2n-4 2n-6 ... -2 -4**  black

\doublespace
\singlespace

\doublespace
Figure 24
\singlespace

\textbf{COR 4:} The right corridor as it passes from the first level to the second level can't contain the subcode ... 6 $4^{2k}$ 2 with $k\geq0$.

Proof: First observe that k=0 is impossible by Cor 2 as it must at least be of the form ... 6 4 ... Now  observe that  a 2 cannot be the last code number since the first 4 corresponds to the two special blue points (and a blue line) which would mean that this 2 would also correspond to a blue point (since the color of the lines alternate) which contradicts the definition of a corridor that the last C point is black. On the other hand  if 2 is not the last code number then 6 $4^{2k}$ 2 ...  is impossible by the second part of Rhombus rule C using n=1 . QED 

\

\textbf{Corridor Lemma 2:} The right corridor as it passes from the first level to the second level can't be of the form  
...6 $4^{2k+1}$ 2... with $k\geq0$ where 2 is not the last code number.

Proof: Consider the two special blue C points $C_{m}$, $C_{m+1}$ (the first two points starred below) and the C points $C_{m+4k+2}$ and $C_{m+4k+3}$ The vector v from the blue point $C_{m}$ to the black point $C_{m+4k+2}$ is exactly equal (and hence parallel) to the vector w from the blue point $C_{m+1}$ to the blue point $C_{m+4k+3}$ since both equal $(sin2x,2k+1+(2k+1)cos2x)$ and hence by the Non-Corridor Test 1 there is no corridor. 

\doublespace

\singlespace

\noindent
$\left.
\begin{aligned}
&\text{-2** 0*} \\
&\text{2 0} \\
&\vdots\quad\vdots\\
&\underline {\text{-2 0 }} \\
\end{aligned} 
\;\right\}$ 2k+1 lines

\noindent
\underline{2**}  black

\noindent
0* ... blue

\noindent QED

\doublespace
\singlespace

\textbf{COR 5:}  It follows that the right corridor as it passes from the first level to the second level ends at ... 6 4 2 or continues on in the form  ... 6 4 4 ...

Proof: If it doesn't end at ... 6 4 2 and continues on then since the second level must end in a 2, the only choices are ... 6 4 2 ... which is impossible by Corridor Lemma 2 since 2 is not the last code number or ... 6 4 6 ... which is impossible by Rhombus Rule B or ... 6 4 4 ... which must be the case.
QED

\

\textbf{Rhombus Rule D:} The following subcodes never appear in the code of any rhombus poolshot

\noindent a. $(2n)^{2}$ $(2n+2)^{3}$  for $n>1$ Note: n=1 or $2^{2}$ $4^{3}$ is possible (see Fig. 25)

\noindent b. ... $(2n)^{3}$ $(2n+2)^{2}$ ... for $n>0$ if there is a previous and following code number. 

\noindent c. $(2n)^{k}$ $(2n+2)^{2}$ ... for $n>0$, $k>3$ if there is a following code number.

Note: c follows immediately from b.

Proof: a. Let the blue-black vector between the first two starred points below be $v=(a,b)$ where $a=2sin2x+2sin4x+...+2 sin(2n-4)x+sin(2n-2)x$ and $b=2+2cos2x+2cos4x+...+2cos(2n-4)x+cos(2n-2)x$ and let the black-blue vector between the second two starred points be $w=(c,d)$ where $c=sin2x+3sin4x+...+3 sin(2n-2)x+sin(2n)x$ and 

\noindent $d=3+5cos2x+3cos4x+...+3cos(2n-2)x+cos(2n)x$. Then  $ad<bc$ $\leftrightarrow$ $0<(bc-ad)$ $\leftrightarrow$ $0<16sin^4x(bc-ad)$ which by algorithm two reduces to $-3sin2x+3sin4x-sin6x-sin(2n-4)x+3sin(2n-2)x-3sin2nx+sin(2n+2)x>0$. By the Main Trig Identity with k=2 and s=n-3, this is the same as  
$16sin^4x(sin4x+sin6x+.. +sin(2n-2)x )> 0$. This last condition holds since $(2n+2)x<180$ because 2n+2 appears in the code and hence by the Non-Rhombus Tower Test 4 there is no rhombus poolshot. Observe we have oriented the sequence of mirror images so that the first starred integer is a blue -2*.

\doublespace
\singlespace

-2* 0 2 4 ... 2n-4 blue

\underline{2n-2 2n-4 ... 2 0*} black

-2 0 2 4 ... 2n-4 2n-2

2n 2n-2 2n-4 ... 2 0

{-2 0 2 4 ... 2n-4 2n-2*} blue

\doublespace

\singlespace

Note: In the special case $n=2$, we get $v=(sin2x,2+cos2x), w=(sin2x+sin4x,3+5cos2x+cos4x)$ and $ad=bc$ which still means that there is no rhombus pool shot by the Non-Rhombus Tower Collinear Test.

\doublespace

Figure 25

\singlespace
b. Let the black-blue vector between the first two starred points be $w=(c,d)$ where $c=sin2x+3sin4x+...+3sin(2n-4)x+sin(2n-2)x$ and $d=3+5cos2x+3cos4x+...+3cos(2n-4)x+cos(2n-2)x$ and let the blue-black vector between the last two starred points be $v=(a,b)$ where $a=sin4x+2sin6x+...+2sin(2n-2)x$ and 

\noindent $b=2+4cos2x+3cos4x+2cos6x+...+2cos(2n-2)x$. Then $ad<bc$  reduces by using algorithm two to $-3sin2x+3sin4x-sin6x-sin(2n)x+3sin(2n+2)x-3sin(2n+4)x+sin(2n+6)x>0$ and by the Main Trig Identity with k=2 and s=n-1 this becomes equivalent to 
$sin4x+sin6x+.. +sin(2n+2)x > 0$. This holds since $(2n+2)x<180$ and hence by the Non-Rhombus Tower Test 4 part c there is no rhombus poolshot. Here we have oriented the sequence of mirror images so that the first starred integer is a black 0*.

\doublespace

\singlespace

\underline{... 0*} black

-2 0 2 4 ... 2n-4 

2n-2 2n-4 ... 2 0 

\underline{-2 0 2 4 ... (2n-4)*} blue

2n-2 2n-4 ... 2 0 -2

\underline{-4 -2 0 2 4 ... 2n-4} 

(2n-2)* ... black

\doublespace

\singlespace

QED

Note: In the special case n=1, we get $w=(-2sin2x,1+2cos2x), v=(-2sin2x-sin4x,2+2cos2x+cos4x)$ and hence $ad<bc$ is equivalent to $sin4x>0$ which holds since $4x<180$ here.

In the special case n=2, we get $w=(-sin2x,3+3cos2x), v=(-sin4x,2+4cos2x+cos4x)$ and hence $ad<bc$ is equivalent to $sin4x+sin6x>0$ which holds since $6x<180$ here.

In the special case n=3, we get $w=(sin2x+sin4x,3+5cos2x+cos4x), v=(sin4x,2+4cos2x+3cos4x)$ and hence $ad<bc$ is equivalent to $sin4x+sin6x+sin8x>0$ which holds since $8x<180$ here.

\doublespace
\singlespace

\textbf{Rhombus Rule E:} Let $n>0$ then the subcodes

\noindent a. ... $(2n)^{2k+1}$ 2n+2 $(2n+4)^{2t+1}$ 2n+6 where $k>t\geq0$ never appears in the code of any rhombus poolshot.

\noindent b. ... $(2n)^{2k+1}$ 2n+2 $(2n+4)^{2t}$ 2n+2 ... where $k>t\geq0$ never appears in the code of any rhombus poolshot. (Note: t=0 follows from Rhombus Rule D part b)

Proof: a. It is enough to show that ... $(2n)^{2t+1}$ 2n+2 $(2n+4)^{2t-1}$ 2n+6 never appears for $t>0$. Observe that the three starred points form a black-blue-black collinear situation. Let v=(a,b) be a vector between the first two starred points where $a=tsin2x+(2t+1)sin4x+... +(2t+1)sin(2n-4)x+tsin(2n-2)x$ and 

\noindent
$b=2t+1+(3t+2)cos2x+(2t+1)cos4x+...+(2t+1)cos(2n-4)x+tcos(2n-2)x$ and w=(c,d) be a vector between the last two starred points where $c=tsin4x+(2t+1)sin6x+...+(2t+1)sin(2n-2)x+tsin2nx$ and $d=2t+1+(4t+2)cos2x+(3t+2)cos4x+(2t+1)cos6x+...+(2t+1)cos(2n-2)x+tcos2nx$. Then ad=bc and the points are collinear, hence no poolshot by the  Non-Rhombus Tower Collinear Test.

\doublespace
\singlespace

\noindent
\underline{... 0*}

\noindent
$\left.
\begin{aligned}
&\text{-2 0 2 4 ... 2n-4} \\
&\text{2n-2 2n-4 ... 2 0} \\
&\vdots\quad\vdots\\
&\underline {\text{-2 0 2 4 ... (2n-4)*}} \\
\end{aligned} 
\;\right\}$ 2t+1 lines

\noindent
\underline {2n-2 2n-4 ... 2 0 -2}

\noindent
$\left.
\begin{aligned}
&\text{-4 -2 0 2 4 ... 2n-4 2n-2} \\
&\text{2n 2n-2 2n-4 ... 2 0 -2} \\
&\vdots\quad\vdots\\
&\underline {\text{-4 -2 0 2 4 ... 2n-4 2n-2}} \\
\end{aligned} 
\;\right\}$ 2t-1 lines

\noindent
\ {2n 2n-2 2n-4 ... 2 0 -2 -4*}

Note: 1. In the special case n=1, we get $v=(-(t+1)sin2x,t+(t+1)cos2x)$, $w=(-(t+1)sin2x-(t+1)sin4x,2t+1+(3t+1)cos2x+(t+1)cos4x)$ and $ad=bc$.

\noindent 2. In the special case n=2, we get $v=(-sin2x,2t+1+(2t+1)cos2x), w=(-sin4x,2t+1+(4t+2)cos2x+(2t+1)cos4x)$ and $ad=bc$.

\noindent 3. In the special case n=3, we get $v=(tsin2x+tsin4x,2t+1+(3t+2)cos2x+tcos4x)$, $w=(tsin4x+tsin6x,2t+1+(4t+2)cos2x+(3t+2)cos4x+tcos6x)$ and $ad=bc$.

b. It is enough to show that ... $(2n)^{2t+1}$ 2n+2 $(2n+4)^{2t-2}$ 2n+2 ... never appears for $t>0$. Let $w=(c,d)$ be the black-blue vector between the first two starred points where $c=tsin2x+(2t+1)sin4x+... +(2t+1)sin(2n-4)x+tsin(2n-2)x$ and 

\noindent
$d=2t+1+(3t+2)cos2x+(2t+1)cos4x+...+(2t+1)cos(2n-4)x+tcos(2n-2)x$ and $v=(a,b)$ be the blue-black vector between the last two starred points where $a=tsin4x+2tsin6x+...+2tsin(2n-2)x+(t-1)sin2nx$ and $b=2t+4tcos2x+3tcos4x+2tcos6x+...+2tcos(2n-2)x+(t-1)cos2nx$. Then $ad<bc$ by using algorithm two reduces to $-3sin2x+3sin4x-sin6x-sin(2n)x+3sin(2n+2)x-3sin(2n+4)x+sin(2n+6)x>0$ and by the Main Trig Identity with k=2 and s=n-1 this becomes equivalent to $sin4x+sin6x+.. +sin(2n+2)x > 0$ which holds since $(2n+4)x<180$. Hence by the Non-Rhombus Tower Test 4 part c there is no rhombus poolshot.

\doublespace
\singlespace

\noindent
\underline{... 0*}

\noindent
$\left.
\begin{aligned}
&\text{-2 0 2 4 ... 2n-4} \\
&\text{2n-2 2n-4 ... 2 0} \\
&\vdots\quad\vdots\\
&\underline {\text{-2 0 2 4 ... (2n-4)*}} \\
\end{aligned} 
\;\right\}$ 2t+1 lines

\noindent
\underline {2n-2 2n-4 ... 2 0 -2}

\noindent
$\left.
\begin{aligned}
&\text{-4 -2 0 2 4 ... 2n-4 2n-2} \\
&\text{2n 2n-2 2n-4 ... 2 0 -2} \\
&\vdots\quad\vdots\\
&\underline {\text{2n 2n-2 2n-4 ... 2 0 -2}} \\
\end{aligned} 
\;\right\}$ 2t-2 lines

\noindent
\underline{-4 -2 0 2 4 ... 2n-4}

\noindent (2n-2)* ...

Note: 1. In the special case n=1, we get $w=(-(t+1)sin2x,t+(t+1)cos2x)$, $v=(-(t+1)sin2x-tsin4x,2t+(3t-1)cos2x+tcos4x)$ and $ad<bc$ is then equivalent to $sin4x>0$ which is true since $6x<180$ .

\noindent 2. In the special case n=2, we get $w=(-sin2x,2t+1+(2t+1)cos2x), v=(-sin4x,2t+4tcos2x+(2t-1)cos4x)$ and then $ad<bc$ is equivalent to $sin4x+sin6x>0$ which is true since $8x<180$ .

\noindent 3. In the special case n=3, we get $w=(tsin2x+tsin4x,2t+1+(3t+2)cos2x+tcos4x)$, $v=(tsin4x+(t-1)sin6x,2t+4tcos2x+3tcos4x+(t-1)cos6x)$ and then $ad<bc$ is equivalent to $sin4x+sin6x+sin8x>0$ which is true since $10x<180$ .

QED

\doublespace
\singlespace

\textbf{Rhombus Rule F:} If a rhombus poolshot has the subcode ... $(2n)^{2k+1}$ 2n+2 ...  where $n>0$, $k\geq0$ and if the code keeps continuing on to the right for at least 2k+1  spots, then it must be of the form ... $(2n)^{2k+1}$ 2n+2 $(2n+4)^{2k}$ ... 

\noindent Hence we write

 ... $(2n)^{2k+1}$ 2n+2 ... \textbf{forces}
... $(2n)^{2k+1}$ 2n+2 $(2n+4)^{2k}$ ... 

\noindent and similarly if the code keeps continuing on to the left for at least 2k+1 spots

... 2n+2 $(2n)^{2k+1}$  ... \textbf{forces}
... $(2n+4)^{2k}$ 2n+2 $(2n)^{2k+1}$ ... 

\noindent Note that if k=0, this means the subcode does not change.

Proof: Since the rhombus poolshot is at least of the code form ... $(2n)^{3}$ 2n+2 ... and keeps continuing on to the right, the only possibilities are

a. ... $(2n)^{2k+1}$ 2n+2 2n ... which is impossible by Rhombus rule B(a) since 2n is not the last code number or 

b. ... $(2n)^{2k+1}$ $(2n+2)^{2}$ ... which is impossible by Rhombus rule D(b) since $2k+1\geq3$ or 

c. ... $(2n)^{2k+1}$ 2n+2 2n+4 ... which must be the case. 

\noindent Now assume the poolshot is of the form ... $(2n)^{2k+1}$ 2n+2 $(2n+4)^{s}$ ... where $s\geq1$ is odd, then if s=2t+1 with $k>t\geq0$ and since the code continues to the right then it must be of the form 

a. ... $(2n)^{2k+1}$ 2n+2 $(2n+4)^{2t+1}$ 2n+6 which can never happen by Rhombus rule E(a) or

b. ... $(2n)^{2k+1}$ 2n+2 $(2n+4)^{2t+1}$ 2n+2 ... which is impossible by Rhombus rule B(a) or 

c. ... $(2n)^{2k+1}$ 2n+2 $(2n+4)^{2t+2}$ ... which must be the case. 

\noindent Now assume $2t+2<2k$ and since the code continues to the right, then it must be of the form 

a. ... $(2n)^{2k+1}$ 2n+2 $(2n+4)^{2t+2}$ 2n+6  which can never happen by Rhombus rule C or

b. ... $(2n)^{2k+1}$ 2n+2 $(2n+4)^{2t+2}$ 2n+2 ... which is impossible by Rhombus rule E(b) or 

c. ... $(2n)^{2k+1}$ 2n+2 $(2n+4)^{2t+3}$ ... which then must be the case. 

\noindent We can now keep repeating the process until we get ... $(2n)^{2k+1}$ 2n+2 $(2n+4)^{2k}$...

QED

\doublespace
\singlespace

\noindent Note: It follows from a further examination of the proof and the rhombus rules that if the subcode ... $(2n)^{2k+1}$ 2n+2 ... doesn't keep continuing to the right for at least 2k+1 spots, then it would be of one of the following forms for $k>t\geq0$.

1. ... $(2n)^{2k+1}$ 2n+2 2n

2. ... $(2n)^{2k+1}$ 2n+2 $(2n+4)^{2t}$ 2n+2

3. ... $(2n)^{2k+1}$ 2n+2 $(2n+4)^{2t+1}$ 2n+2

4.  ... $(2n)^{2k+1}$ 2n+2 $(2n+4)^{s}$ with $2k+1>s\geq0$

\doublespace

\singlespace
\textbf{Mid Corridor Growth Rule:} The right corridor as it passes from the first level to the second level can't be of the form

\doublespace
\singlespace
\noindent \textbf{a.} ... 6 $4^{2k}$ 6 $8^{2k-1}$ 10 $(12)^{2k-1}$ 14 ... ... ... 4n-6 $(4n-4)^{2k-1}$ 4n-2 $(4n)^{2k-2}$ 4n-2 ... for $n\geq2$, $k\geq1$

Note: The first three cases are

 \noindent ... 6 $4^{2k}$ 6 $8^{2k-2}$ 6 ... (n=2 case)

\noindent ... 6 $4^{2k}$ 6 $8^{2k-1}$ 10 $(12)^{2k-2}$ 10 ... (n=3 case)

\noindent ... 6 $4^{2k}$ 6 $8^{2k-1}$ 10 $(12)^{2k-1}$ 14 $(16)^{2k-2}$ 14 ... (n=4 case)

\doublespace
 
\noindent and it can't be of the form

\singlespace

\noindent \textbf{b.} ... 6 $4^{2k}$ 6 $8^{2k-1}$ 10 $(12)^{2k-1}$ 14 ... $(4n-4)^{2k-1}$ 4n-2 $(4n)^{2k}$ ... for $n\geq2$, $k\geq1$

Note: The first three cases are

 \noindent ... 6 $4^{2k}$ 6 $8^{2k}$ ... (n=2 case)

\noindent ... 6 $4^{2k}$ 6 $8^{2k-1}$ 10 $(12)^{2k}$  ... (n=3 case)

\noindent ... 6 $4^{2k}$ 6 $8^{2k-1}$ 10 $(12)^{2k-1}$ 14 $(16)^{2k}$ ... (n=4 case)

\doublespace
\singlespace
Proof: First recall that the transition from the first to second level starts ... 6 4 ... and note that the last 4n-2 in part a and the 4n in part b can't be the last code numbers in the second level since the second level ends with the code number 2.

a. Let w=(c,d) be the blue-black vector between the two single starred points where the first starred point is the special blue point $C_{m}$ and where $c= sin2x+sin4x+...+sin(2n-2)x$ and $d=2nk-2k+(4nk-6k-1)cos2x+(4nk-10k-1)cos4x...+(6k-1)cos(2n-4)x+(2k-1)cos(2n-2)x$. Let v=(a,b) be the blue-blue vector between the two double starred points the first of which is the special blue point $C_{m+1}$ and where $a=sin2x+sin4x+...+sin2nx$ and $b=2nk-1+(4nk-2k-1)cos2x+(4nk-6k-1)cos4x+...+(6k-1)cos(2n-2)x+(2k-1)cos(2n)x$. Now ad=bc and hence the two vectors are parallel and by Non-Corridor Test 1 the right corridor as it passes from the first level to the second level can't be of the form above. Observe that for k=1 this says that ... 6 $4^{2}$ 6 8 10 12  14 ... 4n-4 $(4n-2)^{2}$ ... is impossible for $n\geq2$. Note we have started the vertical array with the two special blue points.

\doublespace
\singlespace

\noindent
$\left.
\begin{aligned}
&\text{-2* 0**} \\
&\vdots\quad\vdots\\
&\underline {\text{2 0}} \\
\end{aligned} 
\;\right\}$ 2k lines

\noindent
\underline {-2 0 2}

\noindent
$\left.
\begin{aligned}
&\text{4 2 0 -2} \\
&\text{-4 -2 0 2} \\
&\vdots\quad\vdots\\
&\underline {\text{4 2 0 -2}} \\
\end{aligned} 
\;\right\}$ 2k-1 lines

\noindent
\underline {-4 -2 0 2 4}

\noindent
$\left.
\begin{aligned}
&\text{6 4 2 0 -2 -4} \\
&\text{-6 -4 -2 0 2 4} \\
&\vdots\quad\vdots\\
&\underline {\text{6 4 2 0 -2 -4}} \\
\end{aligned} 
\;\right\}$ 2k-1 lines

\noindent
\underline {-6 -4 -2 0 2 4 6}

\noindent
$\left.
\begin{aligned}
&\text{8 6 4 2 0 -2 -4 -6} \\
&\text{-8 -6 -4 -2 0 2 4 6} \\
&\vdots\quad\vdots\\
&\underline {\text{8 6 4 2 0 -2 -4 -6}} \\
\end{aligned} 
\;\right\}$ 2k-1 lines

\vdots

\vdots

\vdots

\noindent
\underline{2n-2 2n-4 ... 2 0 -2 ... -(2n-4)*}

\noindent
\underline {-(2n-2) -(2n-4) ... -2 0 2 ... (2n-2)}

\noindent
$\left.
\begin{aligned}
&\text{2n 2n-2 2n-4 ... 2 0 -2 ... -(2n-2)} \\
&\vdots\quad\vdots\\
&\underline {\text{-2n -(2n-2) -(2n-4) ... -2 0 2 ... (2n-2)}} \\
\end{aligned} 
\;\right\}$ 2k-2 lines

\noindent
\underline{2n 2n-2 2n-4 ...2 0 -2 ... -(2n-4)}

\noindent
-(2n-2)**

\doublespace
\singlespace
b. Let w=(c,d) be the blue-black vector between the two single starred points where the first starred point is the special blue point $C_{m}$ and where $c= sin2x+sin4x+...+sin2nx$ and $d=2nk+1+(4nk-2k+1)cos2x+(4nk-6k+1)cos4x+...+(6k+1)cos(2n-2)x+(2k+1)cos(2n)x$. Let v=(a,b) be the blue-blue vector between the two double starred points the first of which is the special blue $C_{m+1}$ and where $a=sin2x+sin4x+...+sin(2n-2)x$ and $b=2nk-2k+(4nk-6k+1)cos2x+(4nk-10k+1)cos4x...+(6k+1)cos(2n-4)x+(2k+1)cos(2n-2)x$. Now ad=bc and hence the two vectors are parallel and by Non-Corridor Test 1 the right corridor as it passes from the first level to the second level can't be of the form above.

\doublespace
\singlespace

\noindent
$\left.
\begin{aligned}
&\text{-2* 0**} \\
&\vdots\quad\vdots\\
&\underline {\text{2 0}} \\
\end{aligned} 
\;\right\}$ 2k lines

\noindent
\underline {-2 0 2}

\noindent
$\left.
\begin{aligned}
&\text{4 2 0 -2} \\
&\text{-4 -2 0 2} \\
&\vdots\quad\vdots\\
&\underline {\text{4 2 0 -2}} \\
\end{aligned} 
\;\right\}$ 2k-1 lines

\noindent
\underline {-4 -2 0 2 4}

\noindent
$\left.
\begin{aligned}
&\text{6 4 2 0 -2 -4} \\
&\text{-6 -4 -2 0 2 4} \\
&\vdots\quad\vdots\\
&\underline {\text{6 4 2 0 -2 -4}} \\
\end{aligned} 
\;\right\}$ 2k-1 lines

\noindent
\underline {-6 -4 -2 0 2 4 6}

\noindent
$\left.
\begin{aligned}
&\text{8 6 4 2 0 -2 -4 -6} \\
&\text{-8 -6 -4 -2 0 2 4 6} \\
&\vdots\quad\vdots\\
&\underline {\text{8 6 4 2 0 -2 -4 -6}} \\
\end{aligned} 
\;\right\}$ 2k-1 lines

\vdots

\vdots

\vdots

\noindent
\underline {2n-2 2n-4 ... 2 0 -2 ... -(2n-4)}

\noindent
\underline {-(2n-2) -(2n-4) ... -2 0 2 ... (2n-2)**}

\noindent
$\left.
\begin{aligned}
&\text{2n 2n-2 2n-4 ... 2 0 -2 ... -(2n-2)} \\
&\vdots\quad\vdots\\
&\underline {\text{-2n -(2n-2) -(2n-4) ... -2 0 2 ... (2n-2)}} \\
\end{aligned} 
\;\right\}$ 2k lines

\noindent
(2n)* ...

QED

\doublespace
\singlespace

Now i) if the right corridor as it passes from the first level to the second level is of the form 

 ... 6 $4^{2k}$ 6 $8^{2k-1}$ 10 $12^{2k-1}$ 14 ... 4n-2 $(4n)^{2k-2}$ ...  for $n\geq2$, $k\geq1$

\noindent then it is \textbf{forced} to be of the form 

... 6 $4^{2k}$ 6 $8^{2k-1}$ 10 $12^{2k-1}$ 14 ... 4n-2 $(4n)^{2k-1}$ ... 

\noindent as the only other two choices ... 6 $4^{2k}$ 6 $8^{2k-1}$ 10 $(12)^{2k-1}$ 14 ... 4n-2 $(4n)^{2k-2}$ 4n-2 ... and  ... 6 $4^{2k}$ 6 $8^{2k-1}$ 10 $(12)^{2k-1}$ 14 ... 4n-2 $(4n)^{2k-2}$ 4n+2 ... are  impossible by Mid Corridor Growth Rule part a and Rhombus Rule C. Note if k=1 then ... 6 $4^{2}$ 8 10 12 14 ... 4n-4 4n-2 ... forces ... 6 $4^{2}$ 8 10 12 14 ... 4n-4 4n-2 4n ... as ... 6 $4^{2}$ 8 10 12 14 ... 4n-4 $(4n-2)^2$ ...  and 4n-4 4n-2 4n-4 ... are impossible by Mid Corridor Growth Rule  part a and Rhombus Rule B.

\doublespace
\singlespace

ii) If the right corridor as it passes from the first level to the second level is of the form  

... 6 $4^{2k}$ 6 $8^{2k-1}$ 10 $(12)^{2k-1}$ 14 ... 4n-2 $(4n)^{2k-1}$ ... for $n\geq2$, $k\geq1$

\noindent then it is \textbf{forced} to be of the form 

... 6 $4^{2k}$ 6 $8^{2k-1}$ 10 $(12)^{2k-1}$ 14 ... 4n-2 $(4n)^{2k-1}$ 4n+2 ... 

\noindent as ... 6 $4^{2k}$ 6 $8^{2k-1}$ 10 $(12)^{2k-1}$ 14 ... 4n-2 $(4n)^{2k}$ ...  and ... 6 $4^{2k}$ 6 $8^{2k-1}$ 10 $(12)^{2k-1}$ 14 ... 4n-2 $(4n)^{2k-1}$ 4n-2 ... are  impossible by Mid Corridor Growth Rule part b and Rhombus Rule B.

\doublespace
\singlespace
iii) If the right corridor as it passes from the first level to the second level is of the form 

... 6 $4^{2k}$ 6 $8^{2k-1}$ 10 $(12)^{2k-1}$ 14 ... 4n-6 $(4n-4)^{2k-1}$ 4n-2 ...  for $n\geq2$, $k\geq2$

\noindent then by Rhombus rule F it is forced to be of the form 

... 6 $4^{2k}$ 6 $8^{2k-1}$ 10 $(12)^{2k-1}$ 14 ... 4n-6 $(4n-4)^{2k-1}$  4n-2 $(4n)^{2k-2}$ ... 

Note 1. The first two cases are 

\noindent ... 6 $4^{2k}$ 6 ... forces ... 6 $4^{2k}$ 6 $8^{2k-2}$ ...  (n=2 case)

\noindent ... 6 $4^{2k}$ 6 $8^{2k-1}$ 10 ... forces ... 6 $4^{2k}$ 6 $8^{2k-1}$ 10 $(12)^{2k-2}$ ... (n=3 case)

Note 2. By the note to Rhombus Rule F if this subcode doesn't keep continuing to the right for at least 2k-1 spots, then it would have to end in 4n-4, 4n-2 or 4n but since $n\geq2$, none of these values is 2. Hence as the right corridor must end in the code 2 this subcode must keep continuing to the right for at least 2k-1 spots and Rhombus Rule F applies.

\doublespace
\singlespace
\textbf{Corridor Lemma 3:} The right corridor as it passes from the first level to the second level can't be of the form  ... 6 $4^{2k}$ 6 ...  with $k\geq0$.

Proof: Since then by the above results starting with $n=2$ the code sequence ... 6 $4^{2k}$ 6 ... keeps growing arbitrarily large and hence never ends in 2 which is impossible. 

More specifically if $k\geq2$ then ... 6 $4^{2k}$ 6 ...

\noindent forces ... 6 $4^{2k}$ 6 $8^{2k-2}$ ... using iii) with n=2 which

\noindent forces ... 6 $4^{2k}$ 6 $8^{2k-1}$      ... using i) with n=2 which

\noindent forces ... 6 $4^{2k}$ 6 $8^{2k-1}$ 10 ... using ii) with n=2 which

\noindent forces ... 6 $4^{2k}$ 6 $8^{2k-1}$ 10 $12^{2k-2}$ ... using iii) with n=3 and the code numbers keep growing as we keep using i), ii) and iii).

If k=1, we cannot use Rhombus Rule F and instead use Mid Corridor Growth Rule a) and b) and Rhombus Rule B as follows.

Starting with ... 6 $4^{2}$ 6 ..., this

\noindent forces ... 6 $4^{2}$  6 8 ... since ... 4 6 4 ... by B and ... 6 $4^{2}$  $6^{2}$ ... by a) with k=1 are impossible, which

\noindent forces ... 6 $4^{2}$  6 8 10 ... since ... 6 8 6 ... by B and ... 6 $4^{2}$  6  $8^{2}$ ... by b) with k=1 are impossible which 

\noindent forces .. 6 $4^{2}$  6 8 10 12 ... since ... 8 10 8 ... by B and ... 6 $4^{2}$  6 8  $10^{2}$ ... by a) with k=1 are impossible
and the code numbers keep growing as we keep using a), b) and Rhombus Rule B.

Note k=0 is impossible as the code sequence must be of the form ... 6 4 ... by Cor 2.

QED

\doublespace
\singlespace
\textbf{Corridor Lemma 4:} The right corridor as it passes from the first to the second level is exactly of the form ... 6 $4^{2k+1}$ 2 for $k\geq0$ (which can also be written in the form ... 6 $4^{2s-2}$ 4 2 starting with $s\geq1$).

Proof: We have showed that the code sequence of the right corridor as it passes from the first to the second level

1. is of the form ... 6 4 .. by Cor 2

2. is not of the form ... 6 $4^{2k+1}$ 2 ... with $k\geq0$ by Corridor Lemma 2 where 2 is not the last code number

3. is not of the form ... 6 $4^{2k}$ 2 ... or  ... 6 $4^{2k}$ 2 with $k\geq0$ by Cor 4 

4. is not of the form ... 6 $4^{2k}$ 6 ... with with $k\geq0$ by Corridor Lemma 3 

5. is not of the form ... 6 $4^{2k+1}$ 6 ... with $k\geq0$ by Cor 3

6. is not of the form ... 6 $4^{2k+1}$ 6 or ... 6 $4^{2k}$ 6 with $k\geq0$ since it must end in the code 2 by Corridor Lemma 1

Since there is only one remaining choice, this means that it must be exactly of the form ... 6 $4^{2k+1}$ 2 where 2 is the last code number. QED

Note this means the second level is exactly described by the code 2 $4^{2k}$ 2 for some $k\geq0$.

\subsection*{11. The Classification Theorem}

We are now in a position to classify all corridor towers.

\doublespace
\singlespace
\textbf{The Classification Theorem:} The only corridor rhombus towers are of the code sequence form 2 $ J^N$ 4 2 for $N\geq1$ where J=6 $8^{2k-2}$ 6 $4^{2k-2}$ with $k\geq1$.

\doublespace
\singlespace
 2 (6 $8^{2k-2}$ 6 $4^{2k-2}$)  (6 $8^{2k-2}$ 6 $4^{2k-2}$) ... (6 $8^{2k-2}$ 6 $4^{2k-2}$) 4 2

\doublespace

\noindent or an extended version of the above with code sequence form

2 $K^s$ $ J^N$ 4 2 where K= $ J^N$ 4 $L^N$ 4 where L = $2^2$ $4^{2k-2}$ for $s\geq0$.

\singlespace
\noindent Note: It follows that the first level is exactly of the form 2 $ J^{N-1}$ 6 $8^{2k-2}$ 6 2 or in the extended version of the form 2 $K^s$ $ J^{N-1}$ 6 $8^{2k-2}$ 6 2 while the second level is exactly of the form 2 $4^{2k-2}$ 2.

\doublespace
\singlespace
To prove this we need to develop many more rules. In the meantime, it is easy to prove the following.

\noindent
\textbf{Fact:} The ratio of the width of second corridor to that of the first corridor is N+1 to 1 or equivalently the ratio of the width of the left corridor to the right corridor is N to 1.

\doublespace
\singlespace
\noindent Proof: Given the corridor rhombus tower 2 $ J^N$ 4 2 in standard position, the coordinates of $C_1$ are (0,1), of the special point $C_m$ are (u,v) where $u=(N-1)sin2x+Nsin4x$ and $v=N(4k-2)-2k+3+[N(6k-3)-2k+3]cos2x+N(2k-1)cos4x$ and of $C_n$ are (e,f) where $e=Nsin2x+Nsin4x$ and $f=N(4k-2)+2+[N(6k-3)+2]cos2x+N(2k-1)cos4x$. This means that the slope of the boundary line from $C_1$ to $C_n$ is given by m=b/a where a=Nsin2x+Nsin4x and b=N(4k-2)+1+[N(6k-3)+2]cos2x+N(2k-1)cos4x. Now if D is the intersection of this boundary line with the line through $C_{m+1}$ and $C_m$ then the coordinates of D are (u,1+mu). But then the distance from $C_m$ to D is v-1-mu=1/N and the result follows.

\doublespace

\singlespace
\noindent
\underline 0*

\noindent
$\left.
\parbox{0.5\linewidth}
{\underline {-2 0 2}

\noindent
$\left.
\begin{aligned}
&\text{4 2 0 -2}\\
&\quad\vdots\quad\vdots\\
&\underline {\text{-4 -2 0 2}}
\end{aligned} 
\;\right\}$ $2k-2$ \text{ times}

\underline {4 2 0}

\noindent
$\left.
\begin{aligned}
&\text{-2* 0}\\
&\quad \vdots \quad \vdots\\
&\underline {\text{2 0}}
\end{aligned} 
\;\right\}$ $2k-2$ \text{ times}

}
\right\}$ N times

\noindent
\underline{-2 0}

\noindent 2*

QED.

\doublespace

Figure 26

\singlespace

\noindent
\textbf{Consequence 1:} It follows that any poolshot leaving the base AB within and parallel to the right corridor of the form 2 $ J^N$ 4 2 can be considered to enter and re-enter the left corridor N times before returning to the right corridor and  becoming periodic after 32Nk-16N+8 reflections. If we then straighten out this poolshot and form the corresponding rhombus tower (not corridor tower), it will be described by the code sequence 2 $ J^N$ 4 2 $ (2\; 4^{2k-2}\; 2)^N$=2 $ J^N$ 4 $ (2^2\; 4^{2k-2})^N$ 2=2 $ J^N$ 4 $L^N$ 2. Observe that the first and last 2 are separated by an odd number of code numbers which means that the corresponding side sequence is of the form 1312...  ...2131. In other words, the first and last reflections are in side AC.

\doublespace
Figure 27

\singlespace
\noindent
\textbf{Consequence 2:} It further follows that we can extend the corridor rhombus tower backwards to add in multiples of this periodic path to create longer corridor rhombus towers. These will have the form (2 $ J^N$ 4 $L^N$ 2) extended to (2 $ J^N$ 4 $L^N$ 2) extended to ...(2 $ J^N$ 4 $L^N$ 2) extended to (2 $ J^N$ 4 2) which is the same as 2 $ J^N$ 4 $L^N$ 4 $ J^N$ 4 $L^N$ 4 $ J^N$ 4 $L^N$  ...
4 $ J^N$ 4  2 because of the preceding observation. These are just longer versions of 2 $ J^N$ 4 2 with code sequences
2 $K^s$ $ J^N$ 4 2 where K= $ J^N$ 4 $L^N$ 4 and L = $2^2$ $4^{2k-2}$ and with the same ratio N to 1 of the two corridors. We will call these the \textbf{periodic corridor extensions} of 2 $ J^N$ 4 2.

\doublespace

Figure 28

\singlespace

\subsection*{12. Two Subcodes Rules}

\textbf{Rhombus Rule G} The subcodes 6 $4^{2k}$ 6 and ... $4^{2k+2}$ ... with  $k\geq0$ cannot occur together in a code sequence of a rhombus pool shot if the corresponding vertical array of the code sequence produces a blue-black vector v and a black-blue vector w which have the same set of integers between them (excluding the initial points and including the terminal points). 

Proof: The two vectors would then be parallel and by the Non-Rhombus Tower Test 6 there is no rhombus pool shot. QED 

Note: This can happen if the vertical array is as follows where the integers between the starred points and the double starred points are exactly the same.

\doublespace
\singlespace

\underline {-2* 0 2}  blue

$\left.
\begin{aligned}
&\text{4 2} \\
&\text{0 2} \\
&\vdots\quad\vdots\\
&\underline {\text{0 2 }} \\
\end{aligned} 
\;\right\}$ 2k lines

\underline {4 2 0*}  black

\doublespace
\singlespace

\quad\vdots\quad\vdots\\

\underline{...2**} black

0 2 blue

4 2 black

$\left.
\begin{aligned}
&\text{0 2} \\
&\text{4 2} \\
&\vdots\quad\vdots\\
&\underline {\text{4 2}} \\
\end{aligned} 
\;\right\}$ 2k lines

0**... blue

Note: As a special case if k=0 this means $6^2$ and ... $4^2$ ... cannot occur in the code sequence of a rhombus pool shot if they produce two parallel vectors as above. An example where this happens was the $6^2$ 4 2 $4^2$... example previously given.

\doublespace
\singlespace

\noindent
\textbf{Rhombus Rule H:} The subcodes ... 6 $8^{2k'}$ 6 ... and ... $8^{2k}$ ... with $k>k'\geq0$ are not part of the code sequence of a rhombus pool shot if they produce two vectors (for example as shown below) one w=(c,d) which is black-blue and the other v=(a,b)  which is blue-black and for which $ad< bc$.

Proof: We can suppose that ... 6 $8^{2k'}$ 6 ... produces the black-blue vector 

\noindent $w=(c,d)=(sin4x,(2k'+2)+(4k'+4)cos2x+(2k'+1)cos4x)$ going between the two starred points below.

\singlespace
\noindent
\underline{...0*} black

\noindent
\underline {-2 0 2}

\noindent
$\left.
\begin{aligned}
&\text{4 2 0 -2} \\
&\text{-4 -2 0 2} \\
&\vdots\quad\vdots\\
&\underline {\text{-4 -2 0 2}} \\
\end{aligned} 
\;\right\}$ 2k' lines

\noindent
\underline {4 2 0}

\noindent
-2*... blue

\doublespace
\singlespace
\noindent On the other hand we can further suppose that  ... $8^{2k}$ ... produces the blue-black vector 

\noindent $v=(a,b)=(sin4x,2k+4kcos2x+(2k+1)cos4x)$ going  between the following two starred points.

\singlespace
\noindent
\underline {... 2*} blue

\noindent
$\left.
\begin{aligned}
&\text{4 2 0 -2} \\
&\text{-4 -2 0 2} \\
&\vdots\quad\vdots\\
&\underline {\text{-4 -2 0 2}} \\
\end{aligned} 
\;\right\}$ 2k lines

\noindent
4* ... black

\doublespace
\singlespace
\noindent Notice that $ad<bc$ is equivalent to $d<b$  since a=c=sin4x and $sin4x>0$ since $8x<180$. But then if $ k'<k$, $d=2k'+2+(4k'+4)cos2x+(2k'+1)cos4x<b=2k+4kcos2x+(2k+1)cos4x$ observing that $cos4x>0$ since $8x<180$. Hence by Non-Rhombus Tower Test 5 the two subcodes cannot be part of the same rhombus poolshot.

QED

\noindent Examples: By Rhombus Rule H, we get the following examples.

\noindent 1: As a special case if $k'=0$, ...$6^2$... and ... $8^2$ ... cannot occur in the code sequence of a rhombus pool shot if they produce two parallel vectors as above. For example ...$6^2$ 4 4 6 $8^2$... is impossible.

\noindent 2: The subcode  ... 6 $8^{2s-4}$ 6 $4^{2s-2}$ 6 $8^{2s-2}$ ... never appears in a rhombus pool shot  using $k'=s-2$ and $k=s-1$ for $s\geq2$.

\doublespace
\singlespace

\noindent
\textbf{Rhombus Rule I:} The subcodes  ... 6 $8^{2k'+1}$ 10  and ... $4^{2k+1}$ with $k>k'\geq0$ are not part of the code sequence of a rhombus pool shot if they produce two vectors one w=(c,d) which is black-blue and the other v=(a,b)  which is blue-black and for which ${ad}\leq{bc}$.

Proof: We can suppose that ...$4^{2k+1}$ produces the blue-black  vector $v=(a,b)=(sin2x,(2k+1)+(2k+1)cos2x)$ going between the following two starred points 

\singlespace
\noindent
\underline{...0*}

\noindent
$\left.
\begin{aligned}
&\text{2 0} \\
&\text{-2 0} \\
&\vdots\quad\vdots\\
&\underline {\text{2 0*}} \\
\end{aligned} 
\;\right\}$ 2k+1 lines

\doublespace
\singlespace
\noindent and that ... 6 $8^{2k'+1}$ 10 produces the  black-blue vector $w=(c,d)=
(sin4x,2k'+3+(4k'+6)cos2x+(2k'+3)cos4x)$ going between the below two starred points.

\noindent
\underline 0*

\noindent
\underline {-2 0 2}

\noindent
$\left.
\begin{aligned}
&\text{4 2 0 -2} \\
&\text{-4 -2 0 2} \\
&\vdots\quad\vdots\\
&\underline {\text{4 2 0 -2}} \\
\end{aligned} 
\;\right\}$ 2k'+1 lines

\noindent
-4 -2 0 2 4*

\doublespace
\singlespace
\noindent
If $k'=k-1$ then $v=(a,b)=(sin2x,2k+1+(2k+1)cos2x)$, $w=(c,d)=(sin4x,2k+1+(4k+2)cos2x+(2k+1)cos4x)$ and $ad=bc$. If $k'<k-1$ then $d<2k+1+(4k+2)cos2x+(2k+1)cos4x$ and $ad<bc$ observing that cos2x and cos 4x are positive since 8 occurs in the code. In all cases $ad\leq bc$ and there is no rhombus poolshot by the Non-Rhombus Tower Test 7.

QED

\noindent Examples: By Rhombus Rule I the following subcodes (reversing the order of the code numbers) produce exactly the situation above and hence are impossible in a rhombus poolshot where J= 6 $8^{2s-2}$ 6 $4^{2s-2}$.

\noindent 1. 10 $8$ 6 $4^{2}$ 4 ...  using $k=1, k'=0$

\noindent 2. 10 $8$ 6 $4^{2}$ $J$ 4 ...  using $k=1, k'=0, s=2$ 

\noindent 3. 10 $8$ 6 $4^{2}$ $ J^i$ 4 ...  using $k=1, k'=0, s=2$ for $i\geq0$

\noindent 4. 10 $8^{2s-3}$ 6 $4^{2s-2}$ 6 $8^{2s-2}$ 6 $4^{2s-1}$ ... using $k=s-1$ and $k'=s-2$, $s\geq2$

\noindent 5. 10 $8^{2s-3}$ 6 $4^{2s-2}$ $ J^i$ 4... using $k=s-1$ and $k'=s-2$ for $i\geq0$, $s\geq2$

\doublespace
\singlespace
\noindent \textbf{Rhombus Rule J:} 2n+2 $(2n)^{2k'}$ 2n+2 and ... $(2n)^{2k}$ ... with $0\leq{k'}\leq(k-1)$ and $1\leq{n}$ are not part of the code sequence of a rhombus pool shot if they produce two \textbf{identical} vectors one w=(c,d) which is black-blue and the other v=(a,b)  which is blue-black (for example as below). 

Proof: Consider the black-blue vector w between the two starred C points corresponding to the code ... $(2n)^{2k}$ ...

\doublespace
\singlespace

\noindent
\underline {...2*}

\noindent
$\left.
\begin{aligned}
&\text{0 2 4 ... 2n-2} \\
&\text{2n 2n-2 ...   4 2} \\
&\vdots\quad\vdots\\
&\text{0 2 4 ... 2n-2} \\
&\underline {\text{2n 2n-2 ...   4 2}} \\
\end{aligned} 
\;\right\}$ 2k lines

\noindent 0* ...

\doublespace
\singlespace
\noindent and the blue-black vector v between the two starred C points corresponding to the code 2n+2 $(2n)^{2k'}$ 2n+2

\noindent
\underline {-2* 0 2 ... 2n-2}

\noindent
$\left.
\begin{aligned}
&\text{2n 2n-2 ...   4 2} \\
&\text{0 2 4 ... 2n-2} \\
&\vdots\quad\vdots\\
&\text{2n 2n-2 ...   4 2} \\
&\underline {\text{0 2 4 ... 2n-2}} \\
\end{aligned} 
\;\right\}$ $2k'$ lines

\noindent 2n 2n-2 ... 2 0*

\doublespace
\singlespace
\noindent These two vectors are exactly the same if $2k'+2=2k$ which is equivalent to $k'=k-1$.   It follows that by Non-Rhombus Tower Test 6, the two subcodes cannot be part of the code sequence of a rhombus pool shot. A similar argument is made if $0\leq{k'}\leq(k-1)$.

\noindent QED

\doublespace
\singlespace

\noindent \textbf{COR:} The subcode 2n+2 $(2n)^{2k'}$ 2n+2 $(2n+4)^{2k-2}$ 2n+2 $(2n)^{2k-2}$ 2n+2 $(2n+4)^{2k-2}$ 2n+2 $(2n)^{2k-1}$ ... for $0\leq{k'}\leq(k-2)$, $n\geq1$ is impossible in a rhombus poolshot.

 Proof: Assuming $0\leq{k'}\leq(k-2)$, we get the vertical array below. Now use Rhombus Rule J.

\doublespace
\singlespace

\noindent
\underline {-2* 0 2 ... 2n-2}  

\noindent
$\left.
\begin{aligned}
&\text{2n 2n-2 ...   4 2} \\
&\vdots\quad\vdots\\
&\underline {\text{0 2 4 ... 2n-2}} \\
\end{aligned} 
\;\right\}$ $2k'$ lines

\noindent 
\underline {2n 2n-2 ... 2 0*}  

\noindent
$\left.
\begin{aligned}
&\text{-2 0 2 ... 2n-2 2n} \\
&\vdots\quad\vdots\\
&\underline {\text{2n+2 2n ... 2 0}} \\
\end{aligned} 
\;\right\}$ $2k-2$ lines

\noindent 
\underline {-2 0 2 ... 2n-2}

\noindent
$\left.
\begin{aligned}
&\text{2n 2n-2 ... 4 2} \\
&\vdots\quad\vdots\\
&\underline {\text{0 2 4 ... 2n-2}} \\
\end{aligned} 
\;\right\}$ $2k-2$ lines

\noindent 
\underline {2n 2n-2 ... 2 0}

\noindent
$\left.
\begin{aligned}
&\text{-2 0 2 ... 2n-2 2n} \\
&\vdots\quad\vdots\\
&\underline {\text{2n+2 2n ... 2 0}} \\
\end{aligned} 
\;\right\}$ $2k-2$ lines

\noindent 
\underline {-2 0 2 ... 2n-2}

\noindent
$\left.
\begin{aligned}
&\text{2n 2n-2 ... 4 2*} \\
&\text{0 2 4 ... 2n-2} \\
&\vdots\quad\vdots\\
&\underline {\text{2n 2n-2 ... 4 2}} \\
\end{aligned} 
\;\right\}$ $2k-1$ lines

\noindent 
0* ...

\noindent QED

\doublespace
\singlespace
If we let  R= 2n+2 $(2n)^{2k-2}$ 2n+2 $(2n+4)^{2k-2}$ this corollary reads 2n+2 $(2n)^{2k'}$ 2n+2 $(2n+4)^{2k-2}$ R 2n+2 $(2n)^{2k-1}$ ... for $0\leq{k'}\leq(k-2)$ is impossible in a rhombus poolshot. 

Similarly for 2n+2 $(2n)^{2k'}$ 2n+2 $(2n+4)^{2k-2}$ $R^i$ 2n+2 $(2n)^{2k-1}$ ...
for $i\geq{0}$.

Observing that for i=0 the above rule still applies assuming $0\leq{k'}\leq(k-2)$.

\noindent
\underline {-2* 0 2 ... 2n-2}   

\noindent
$\left.
\begin{aligned}
&\text{2n 2n-2 ...   4 2} \\
&\vdots\quad\vdots\\
&\underline {\text{0 2 4 ... 2n-2}} \\
\end{aligned} 
\;\right\}$ $2k'$ lines

\noindent 
\underline {2n 2n-2 ... 2 0*}  

\noindent
$\left.
\begin{aligned}
&\text{-2 0 2 ... 2n-2 2n} \\
&\vdots\quad\vdots\\
&\underline {\text{2n+2 2n ... 2 0}} \\
\end{aligned} 
\;\right\}$ $2k-2$ lines

\noindent 
\underline {-2 0 2 ... 2n-2}

\noindent
$\left.
\begin{aligned}
&\text{2n 2n-2 ... 4 2*} \\
&\vdots\quad\vdots\\
&\underline {\text{2n 2n-2 ... 4 2}} \\
\end{aligned} 
\;\right\}$ $2k-1$ lines

\noindent 
0* ...

\doublespace
\singlespace
\noindent Note: This leads to 

\noindent \textbf{Forcing Rule 1}

... $(2n)^{2k'}$ 2n+2 $(2n+4)^{2k-2}$ $R^i$ 2n+2 $(2n)^{2k-1}$ ... 

\noindent forces

 ...$(2n)^{2k'+1}$ 2n+2 $(2n+4)^{2k-2}$ $R^i$ 2n+2 $(2n)^{2k-1}$ ...

\noindent for $i\geq{0}$, $0\leq{k'}\leq{k-2}$, $n\geq1$ (assuming that the subcode extends to the left for at least 2 spots) 

Proof: If $n>1$ then ... 2n-2 $(2n)^{2k'}$ 2n+2 is impossible by Rhombus Rule C and  2n+2 $(2n)^{2k'}$ 2n+2 $(2n+4)^{2k-2}$ $R^i$ 2n+2 $(2n)^{2k-1}$ ... is impossible by this corollary. If n=1, this says ... $2^{2k'}$ 4 $6^{2k-2}$ $R^i$ 4 $2^{2k-1}$ ... forces ... $2^{2k'+1}$ 4 $6^{2k-2}$ $R^i$ 4 $2^{2k-1}$ ... as
4 $2^{2k'}$ 4 $6^{2k-2}$ $R^i$ 4 $2^{2k-1}$ ... is impossible by this corollary and since 2 is the smallest code number, there is no other choice.

QED

\doublespace
\singlespace
\noindent \textbf{Rhombus Rule K:} The subcodes ...2n-2 $(2n)^{2k'+1}$ 2n+2 and $(2n)^{2k-1}$ ... with ${k-2}\geq{k'}\geq0$ and $n>1$  are not part of the code sequence of a rhombus pool shot if they produce two identical vectors one which is black-blue and the other which is blue-black
(for example as below). In the special case  with n=2 this means ... 2 $4^{2k'+1}$ 6 and $4^{2k-1}$ ... with ${k-2}\geq{k'}\geq0$ are impossible together if they produce the two identical vectors.

Proof: Consider the black-blue vector between the two starred C points corresponding to the code ... 2n-2 $(2n)^{2k'+1}$ 2n+2.

\doublespace
\singlespace

\noindent
\underline {...0*}

\noindent
\underline {-2 0 2 4 ... 2n-6}

\noindent
$\left.
\begin{aligned}
&\text{2n-4 2n-6 ...   4 2 0 -2} \\
&\text{-4 -2 0 2 ... 2n-8 2n-6} \\
&\vdots\quad\vdots\\
&\underline {\text{2n-4 2n-6 ...   4 2 0 -2}} \\
\end{aligned} 
\;\right\}$ $2k'+1$ lines

\noindent -4 -2 0 2 ... 2n-6 2n-4*

\doublespace
\singlespace
\noindent and the blue-black vector between the two starred C points corresponding to the code $(2n)^{2k-1}$ ...

\noindent
$\left.
\begin{aligned}
&\text{-4* -2 0 2 ... 2n-8 2n-6} \\
&\text{2n-4 2n-6 ...   4 2 0 -2} \\
&\vdots\quad\vdots\\
&\text{2n-4 2n-6 ...   4 2 0 -2} \\
&\underline {\text{-4 -2 0 2 ... 2n-8 2n-6}} \\
\end{aligned} 
\;\right\}$ 2k-1 lines

\noindent (2n-4)* ...

\doublespace
\singlespace
\noindent These two vectors are exactly the same if $2k'+3 = 2k-1$ which is equivalent to $k'=k-2$. A similar argument is made if $0\leq{k'}\leq{k-2}$. It follows that by Non-Rhombus Tower Test 6, the two subcodes cannot be part of the code sequence of a rhombus pool shot.

\noindent QED

\doublespace
\singlespace
\noindent \textbf{COR:} The subcode ... 2n-2 $(2n)^{2k'+1}$ 2n+2 $(2n+4)^{2k-2}$ 2n+2 $(2n)^{2k-2}$ 2n+2 $(2n+4)^{2k-2}$ 2n+2 $(2n)^{2k-1}$ ... is impossible in a rhombus poolshot for $0\leq{k'}\leq{k-2}$, $n>1$.

\doublespace
\singlespace
\noindent Proof: This is impossible by applying the Rule above to

\doublespace
\singlespace

\noindent
\underline {...  0*}

\noindent
\underline {-2 0 2 4 ... 2n-6}

\noindent
$\left.
\begin{aligned}
&\text{2n-4 2n-6 ...   4 2 0 -2} \\
&\text{-4 -2 0 2 ... 2n-8 2n-6} \\
&\vdots\quad\vdots\\
&\underline {\text{2n-4 2n-6 ...   4 2 0 -2}} \\
\end{aligned} 
\;\right\}$ $2k'+1$ lines

\noindent
\underline {-4 -2 0 2 ... 2n-6 2n-4*}

\noindent
$\left.
\begin{aligned}
&\text{2n-2 2n-4  ...   2 0 -2 -4} \\
&\text{-6 -4 -2 0 2 ... 2n-6 2n-4} \\
&\vdots\quad\vdots\\
&\underline {\text{-6 -4 -2 0 2 ... 2n-6 2n-4}} \\
\end{aligned} 
\;\right\}$ $2k-2$ lines

\noindent
\underline {2n-2 2n-4  ...  2 0 -2}

\noindent
$\left.
\begin{aligned}
&\text{-4 -2 0 2 ... 2n-8 2n-6} \\
&\text{2n-4 2n-6 ...   4 2 0 -2} \\
&\vdots\quad\vdots\\
&\underline {\text{2n-4 2n-6 ...   4 2 0 -2}} \\
\end{aligned} 
\;\right\}$ $2k-2$ lines

\noindent
\underline {-4 -2 0 2 ... 2n-8 2n-6 2n-4}

\noindent
$\left.
\begin{aligned}
&\text{2n-2 2n-4  ...  2 0 -2 -4} \\
&\text{-6 -4 -2 0 2 ... 2n-6 2n-4} \\
&\vdots\quad\vdots\\
&\underline {\text{-6 -4 -2 0 2 ... 2n-6 2n-4}} \\
\end{aligned} 
\;\right\}$ $2k-2$ line

\noindent
\underline {2n-2 2n-4  ...  2 0 -2}

\noindent
$\left.
\begin{aligned}
&\text{-4* -2 0 2 ... 2n-8 2n-6} \\
&\text{2n-4 2n-6 ...   4 2 0 -2} \\
&\vdots\quad\vdots\\
&\text{2n-4 2n-6 ...   4 2 0 -2} \\
&\underline {\text{-4 -2 0 2 ... 2n-8 2n-6}} \\
\end{aligned} 
\;\right\}$ 2k-1 lines

\noindent (2n-4)* ...

\doublespace
\singlespace
\noindent QED

If we let R= 2n+2 $(2n)^{2k-2}$ 2n+2 $(2n+4)^{2k-2}$ this reads as ... 2n-2 $(2n)^{2k'+1}$ 2n+2 $(2n+4)^{2k-2}$ R 2n+2 $(2n)^{2k-1}$ ... is impossible  in a rhombus poolshot for $0\leq{k'}\leq{k-2}$, $n>1$. 

Similarly for ... 2n-2 $(2n)^{2k'+1}$ 2n+2 $(2n+4)^{2k-2}$ $R^i$ 2n+2 $(2n)^{2k-1}$ ...
for $i\geq{0}$. 

Observing that for i=0 the above rule still applies assuming $0\leq{k'}\leq{k-2}$.

\singlespace

\noindent
\underline {...  0*}

\noindent
\underline {-2 0 2 4 ... 2n-6}

\noindent
$\left.
\begin{aligned}
&\text{2n-4 2n-6 ...   4 2 0 -2} \\
&\text{-4 -2 0 2 ... 2n-8 2n-6} \\
&\vdots\quad\vdots\\
&\underline {\text{2n-4 2n-6 ...   4 2 0 -2}} \\
\end{aligned} 
\;\right\}$ $2k'+1$ lines

\noindent
\underline {-4 -2 0 2 ... 2n-6 2n-4*}

\noindent
$\left.
\begin{aligned}
&\text{2n-2 2n-4  ...   2 0 -2 -4} \\
&\text{-6 -4 -2 0 2 ... 2n-6 2n-4} \\
&\vdots\quad\vdots\\
&\underline {\text{-6 -4 -2 0 2 ... 2n-6 2n-4}} \\
\end{aligned} 
\;\right\}$ $2k-2$ lines

\noindent
\underline {2n-2 2n-4  ...  2 0 -2}

\noindent
$\left.
\begin{aligned}
&\text{-4* -2 0 2 ... 2n-8 2n-6} \\
&\text{2n-4 2n-6 ...   4 2 0 -2} \\
&\vdots\quad\vdots\\
&\text{2n-4 2n-6 ...   4 2 0 -2} \\
&\underline {\text{-4 -2 0 2 ... 2n-8 2n-6}} \\
\end{aligned} 
\;\right\}$ 2k-1 lines

\noindent (2n-4)* ...

\doublespace
\singlespace
\noindent Note: This leads to

\noindent
\textbf{Forcing Rule 2}

 ... $(2n)^{2k'+1}$ 2n+2 $(2n+4)^{2k-2}$ $R^i$ 2n+2 $(2n)^{2k-1}$ ... 

\noindent forces

 ...    $(2n)^{2k'+2}$ 2n+2 $(2n+4)^{2k-2}$ $R^i$ 2n+2 $(2n)^{2k-1}$ ... 

\noindent for $i\geq{0}$, $0\leq{k'}\leq{k-2}$, $n\geq1$ (assuming that the subcode extends to the left for at least 2 spots).

Proof:  2n+2 $(2n)^{2k'+1}$ 2n+2 is impossible by Rhombus Rule B and if $n>1$ ... 2n-2 $(2n)^{2k'+1}$ 2n+2 $(2n+4)^{2k-2}$ $R^i$ 2n+2 $(2n)^{2k-1}$ ... is impossible by this corollary. If n=1, this says that ... $2^{2k'+1}$ 4 $6^{2k-2}$ $R^i$ 4 $2^{2k-1}$ ... forces ... $2^{2k'+2}$ 4 $6^{2k-2}$ $R^i$ 4 $2^{2k-1}$ ... as 4 $2^{2k'+1}$ 4 is impossible by Rhombus Rule B and since 2 is the smallest code number, there is no other choice.

QED

\doublespace
\noindent
\textbf{Rhombus Rule L:}

\singlespace
The subcodes ... 2n $(2n+2)^{2k}$ 2n ... and ... $(2n+2)^{2k+2}$ ... for $k\geq0$, $n>0$ never appear together in a rhombus poolshot if  they create a black-blue vector w=(c,d) and a blue-black vector v=(a,b) for which $ad\leq{bc}$.

(It further follows that ... 2n $(2n+2)^{2k}$ 2n ... and ... $(2n+2)^{2s}$ ... are impossible in a rhombus poolshot if $k<s$ and the conditions above are satisfied.)

Proof:  Apply the Non-Rhombus Tower Test 7.

\noindent QED

\doublespace
\singlespace
\noindent Example:

\singlespace
First suppose that ... 2n $(2n+2)^{2k}$ 2n ... creates the vertical array below and the black-blue vector $w=(c,d)=
((k+2)sin4x+(2k+2)sin6x+...+
(2k+2)sin(2n-4)x+(k+1)sin(2n-2)x,2k+2+(4k+4)cos2x+(3k+2)cos4x+(2k+2)cos6x+...+
(2k+2)cos(2n-4)x+(k+1)cos(2n-2)x)$ between the two starred points.

\doublespace
\singlespace
\noindent
\underline {...0*} black
\singlespace
\noindent
\underline {-2 0 2 4 ... 2n-4}

\noindent
$\left.
\begin{aligned}
&\text{2n-2 2n-4 ... 2 0 -2}\\
&\quad\vdots\quad\vdots\\
&\underline {\text{-4 -2 0 2 ... 2n-6 2n-4}}
\end{aligned} 
\;\right\}$ 2k times

\noindent
\underline {2n-2 2n-4 ... 2 0}

\noindent
-2*... blue

\doublespace
\singlespace
\noindent and second suppose that  ... $(2n+2)^{2k+2}$ ... produces the vertical array below and the blue-black vector $v=(a,b)=((k+1)sin4x+(2k+2)sin6x+...
+(2k+2)sin(2n-4)x+(k+2)sin(2n-2)x,2k+2+(4k+4)cos2x+(3k+3)cos4x+(2k+2)cos6x+...+(2k+2)cos(2n-4)x+(k+2)cos(2n-2)x)$ between the two starred points.

\doublespace
\singlespace
\noindent
\underline{...(2n-4)*} blue

\noindent
$\left.
\begin{aligned}
&\text{2n-2 2n-4 ... 2 0 -2}\\
&\quad\vdots\quad\vdots\\
&\underline {\text{-4 -2 0 2 ... 2n-6 2n-4}}
\end{aligned} 
\;\right\}$ 2k+2 times

\noindent
{(2n-2)*...} black

Note these two vertical arrays above cannot be part of a vertical array of a rhombus poolshot by Rhombus Rule L since $ad<{bc}$ (by using Algorithm Two and special case 1 of the Main Trig Identity) would be  equivalent to $sin(2n+2)x>0$ and this would hold since 2n+2 is part of the subcode.

\doublespace

\singlespace
\textbf{COR 1:} The subcode ... 2n $(2n+2)^{2k}$ 2n $(2n-2)^{2k+2}$ 2n $(2n+2)^{2k+2}$ ... for $k\geq0$, $n\geq1$ never appears in a rhombus poolshot. 

\noindent (Note for k=0, ...$(2n)^2$ $(2n-2)^{2}$ 2n $(2n+2)^{2}$... is impossible and for n=1, ...2 $4^{2k}$ $2^2$ $4^{2k+2}$ ... is impossible and for n=1,k=0, ...$2^3$ $4^2$... is impossible which also follows from Rhombus Rule D(b))

Proof: Since it produces exactly the example above.

\noindent
QED 

Special Cases:

\noindent 1. n=1, $w=(-(2k+2)sin2x-ksin4x, k+1+(2k+2)cos2x+kcos4x)$, $v=(-(2k+2) sin2x-(k+1)sin4x,k+2+(2k+2)cos2x+(k+1)cos4x)$ and $ad=bc$.

\noindent 2. n=2, $w=(-(k+1)sin2x-ksin4x,2k+2+(3k+3)cos2x+kcos4x)$, $v=(-ksin2x-(k+1)sin4x,2k+2+(3k+4)cos2x+(k+1)cos4x)$ and $ad<bc$ is equivalent to $sin6x>0$.

\noindent 3. n=3, $w=(sin4x,2k+2+(4k+4)cos2x+(2k+1)cos4x)$, $v=(sin4x,2k+2+(4k+4)cos2x+(2k+3)cos4x)$ and $ad<bc$ is equivalent to $sin8x>0$.

\doublespace
\singlespace
\textbf{COR 2:} \textbf{(Forcing Rule 3)} In a rhombus poolshot the subcode below (if it extends at least two spots to the left) for $k\geq0$, $n\geq1$

... $(2n+2)^{2k}$ 2n $(2n-2)^{2k+2}$ 2n $(2n+2)^{2k+2}$ ...   

\noindent forces 

... $(2n+2)^{2k+1}$ 2n $(2n-2)^{2k+2}$ 2n $(2n+2)^{2k+2}$ ...

\noindent Proof: Since ... 2n $(2n+2)^{2k}$ 2n $(2n-2)^{2k+2}$ 2n $(2n+2)^{2k+2}$ ... is impossible by COR 1 and 2n+4 $(2n+2)^{2k}$ 2n ...  is impossible by Rhombus Rule C.

\noindent
QED.

\doublespace
\singlespace
\textbf{COR 3:} \textbf{(Rhombus Growth Rule 1)} This leads to the Growth Rule where if we let

\noindent $[S+4n+4]= 4n+10$ $(4n+8)^{2k-2}$ 4n+10 $(4n+12)^{2k-2}$ then for $k\geq2$, $n\geq-2$

\doublespace

... $(4n+12)^{2k-4}$ $[S+4n+4]$ ... (if it extends at least two spots to the left)

\noindent forces

... $(4n+12)^{2k-3}$ $[S+4n+4]$ ...

\noindent Proof: Replace n by 2n+5 and k by k-2 in Cor 2.

\noindent QED

\singlespace
\noindent \textbf{Convention:} We are using the notation that if S = 6 $4^{2k-2}$ 6 $8^{2k-2}$ then $[S+m]$= $m+6$ $(m+4)^{2k-2}$ m+6 $(m+8)^{2k-2}$. Similarly for other code sequences.

\doublespace
\singlespace
\noindent \textbf{Rhombus Rule M:}
The subcodes 2n+6 $(2n+4)^{2s+1}$ 2n+2 ... and $(2n)^{2s+3}$... with  $s\geq{0}$, $n\geq{1}$ never occur together in a  rhombus poolshot if they create a blue-black vector v=(a,b) and a  black-blue vector w=(c,d) with $ad=bc$.

Proof: By the Non-Rhombus Tower Test 6.

\noindent QED.

\noindent
As an example if the subcode 2n+6 $(2n+4)^{2s+1}$ 2n+2 ... produces the blue-black vector v=(a,b) between the starred points below with $a=(2s+3)sin2x+..+(2s+3)sin(2n+2)x+(s+1)sin(2n+4)x$ and $b=s+2+(2s+3)cos2x+...+(2s+3)cos(2n+2)x+(s+1)cos(2n+4)x$

\doublespace
\singlespace

\noindent
\underline {-2* 0 2 4 ...  2n+2}

\noindent
$\left.
\begin{aligned}
&\text{2n+4 2n+2 ... 4 2}\\
&\text{0 2 4 ... 2n 2n+2}\\
&\quad\vdots\quad\vdots\\
&\underline {\text{2n+4 2n+2 ... 4 2}}
\end{aligned} 
\;\right\}$ 2s+1 times

\noindent
\underline {0 2 4 ... 2n-2 2n}

\noindent
(2n+2)* ...

\doublespace
\singlespace 

\noindent
whereas $(2n)^{2s+3}$ ... produces the black-blue vector w=(c,d) between the double starred points below with c=(s+2)sin2x+(2s+3)sin4x+...+(2s+3)sin2nx+(s+1)sin(2n+2)x and d=(s+2)cos2x+(2s+3)cos4x+...+(2s+3)cos2nx+(s+1)cos(2n+2)x

\doublespace
\singlespace

\noindent
$\left.
\begin{aligned}
&\text{(2n+2)** 2n ... 6 4}\\
&\text{2 4 6 ... 2n-2 2n}\\
&\quad\vdots\quad\vdots\\
&\underline {\text{2n+2 2n ... 6 4}}
\end{aligned} 
\;\right\}$ 2s+3 times

\noindent 2** ...

\noindent then $ad=bc$.

\doublespace
\singlespace
\noindent Note a specific example in which this occurs is 2n+6 $(2n+4)^{2s+1}$ 2n+2 $(2n)^{2s+3}$ ... with $s\geq{0}$ so this subcode never occurs in the code of a rhombus poolshot.

\doublespace
\singlespace
\noindent
\textbf{Cor 1:} The subcodes 2n+6 $(2n+4)^{2r+1}$ 2n+2 ... and $(2n)^{2s+3}$ ... with  $s\geq{r}\geq0$, $n\geq{1}$ do not occur together in a rhombus poolshot if they create a  blue-black vector v=(a,b) and a  black-blue vector w=(c,d) with $ad=bc$.

Proof: By the Non-Rhombus Tower Test 6.

\noindent
QED.

\doublespace
\singlespace
\noindent Note a specific example in which this occurs is 2n+6 $(2n+4)^{2r+1}$ 2n+2 $(2n)^{2s+3}$ ... with $s\geq{r}\geq0$ so this subcode never occurs in the code of a rhombus poolshot.

\doublespace
\singlespace
\noindent
\textbf{Cor 2:} The subcode 2n+6 $(2n+4)^{2r+1}$ 2n+2 $(2n)^{2s+2}$ 2n+2 $(2n+4)^{2s+2}$ 2n+2 $(2n)^{2s+3}$ ... where $s\geq{r\geq0}$, $n\geq{1}$ never appears in a rhombus poolshot.

Proof: If s=r, it produces exactly the two vertical arrays in the first example above as seen below and hence is impossible. 

\doublespace
\singlespace

\noindent
\underline {-2* 0 2 4 ... 2n+2}

\noindent
$\left.
\begin{aligned}
&\text{2n+4 2n+2 ... 4 2}\\
&\text{0 2 4 ... 2n+2}\\
&\quad\vdots\quad\vdots\\
&\underline {\text{2n+4 2n+2 ... 4 2}}
\end{aligned} 
\;\right\}$ 2s+1 times

\noindent
\underline {0 2 4 ... 2n-2 2n}

\noindent
$\left.
\begin{aligned}
&\text{(2n+2)* 2n ... 6 4}\\
&\text{2 4 6 8 ... 2n-2 2n}\\
&\quad\vdots\quad\vdots\\
&\underline {\text{2 4 6 8 ... 2n-2 2n}}
\end{aligned} 
\;\right\}$ 2s+2 times

\noindent
\underline {(2n+2) 2n ... 6 4 2}

\noindent
$\left.
\begin{aligned}
&\text{0 2 4 6 8 ... 2n 2n+2}\\
&\text{2n+4 2n+2 ... 4 2}\\
&\quad\vdots\quad\vdots\\
&\underline {\text{2n+4 2n+2 ... 4 2}}
\end{aligned} 
\;\right\}$ 2s+2 times

\noindent
\underline {0 2 4 6 ... 2n-2 2n}

\noindent
$\left.
\begin{aligned}
&\text{(2n+2)** 2n ... 6 4}\\
&\text{2 4 6 8... 2n-2 2n}\\
&\quad\vdots\quad\vdots\\
&\underline {\text{2n+2 2n ... 8 6 4}}
\end{aligned} 
\;\right\}$ 2s+3 times

\noindent 2** ...

\noindent It follows that the subcode 2n+6 $(2n+4)^{2r+1}$ 2n+2 $(2n)^{2s+2}$ 2n+2 $(2n+4)^{2s+2}$ 2n+2 $(2n)^{2s+3}$ ... is impossible in a rhombus poolshot if $s\geq{r}$.

\noindent QED.

\noindent Note: It further follows that if $s\geq{r\geq0}$, $n\geq{1}$ and

1. if V=2n+2 $(2n)^{2s+2}$ 2n+2 $(2n+4)^{2s+2}$ then 2n+6 $(2n+4)^{2r+1}$ $V^t$ 2n+2 $(2n)^{2s+3}$ ... never appears in a rhombus poolshot for any $t\geq0$.

2. if W=$(2n)^{2s+2}$ 2n+2 $(2n+4)^{2s+2}$ 2n+2 then 2n+6 $(2n+4)^{2r+1}$ 2n+2 $W^t$  $(2n)^{2s+3}$ ... never appears in a rhombus poolshot for any $t\geq0$.

\doublespace
\singlespace
\noindent
\textbf{Cor 3:} \textbf{(Forcing Rule 4)} In a rhombus poolshot the subcode below with $k\geq{0}$, $n\geq{1}$ (if it extends at least two spots to the left) 

\doublespace
\noindent
... $(2n+4)^{2k+1}$ 2n+2 $(2n)^{2k+2}$ 2n+2 $(2n+4)^{2k+2}$ 2n+2 $(2n)^{2k+3}$ ... 

forces

\noindent
... $(2n+4)^{2k+2}$ 2n+2 $(2n)^{2k+2}$ 2n+2 $(2n+4)^{2k+2}$ 2n+2 $(2n)^{2k+3}$ ...

\singlespace
Proof: 2n+6 $(2n+4)^{2k+1}$ 2n+2 $(2n)^{2k+2}$ 2n+2 $(2n+4)^{2k+2}$ 2n+2 $(2n)^{2k+3}$ ... is impossible by COR 2 with $r=s=k$ and 
... 2n+2 $(2n+4)^{2k+1}$ 2n+2 ... is impossible by Rhombus Rule B.

\noindent QED.

\doublespace
\noindent Note: It further follows that in a rhombus poolshot if $s\geq{r\geq0}$, $n\geq{1}$ and

\doublespace
\singlespace
\noindent 1. if V=2n+2 $(2n)^{2s+2}$ 2n+2 $(2n+4)^{2s+2}$ then 

\doublespace
\singlespace
\noindent ... $(2n+4)^{2r+1}$ $V^t$ 2n+2 $(2n)^{2s+3}$ ...  (if it extends at least two spots to the left) 

 forces 

\noindent ... $(2n+4)^{2r+2}$ $V^t$ 2n+2 $(2n)^{2s+3}$ ...for any $t\geq0$.

\doublespace
\singlespace

\noindent 2. if W=$(2n)^{2s+2}$ 2n+2 $(2n+4)^{2s+2}$ 2n+2 then 

\doublespace
\singlespace

\noindent ...$(2n+4)^{2r+1}$ 2n+2 $W^t$  $(2n)^{2s+3}$... (if it extends at least two spots to the left) 

forces

\noindent ... $(2n+4)^{2r+2}$ 2n+2 $W^t$  $(2n)^{2s+3}$ ... for any $t\geq0$.

\doublespace
\noindent Note that if t=0, both of these also follow from Rhombus Rule F.

\doublespace
\noindent \textbf{Cor 4:} \textbf{(Rhombus Growth Rule 2)} This leads to the Growth Rule where if
\singlespace
\noindent $[S+4n+4]= 4n+10$ $(4n+8)^{2k-2}$ 4n+10 $(4n+12)^{2k-2}$ then for $k\geq2$, $t\geq0$, $n\geq{-1}$ (assuming the first subcode below extends at least two spots to the left) 

\doublespace

... $(4n+12)^{2k-3}$ $[S+4n+4]^t$ 4n+10 $(4n+8)^{2k-1}$...

\noindent forces

... $(4n+12)^{2k-2}$ $[S+4n+4]^t$ 4n+10 $(4n+8)^{2k-1}$ ...

\noindent Proof: In the first note above replace n by 2n+4 and r and s by k-2.

\noindent QED

\subsection*{13. More Rhombus Rules}

\singlespace
\textbf{Rhombus Rule N:}
The subcode  ... 4n+2 $(4n+4)^{2k}$ 4n+2 $(4n)^{2k}$ 4n+2 $(4n+4)^{2k+2}$ ... for $n\geq0$, $k\geq0$ never appears in the code sequence of any rhombus poolshot. Similarly for the subcode ... $(4n+4)^{2k+2}$ 4n+2 $(4n)^{2k}$ 4n+2 $(4n+4)^{2k}$ 4n+2 ... which is the above in reverse order.

Proof:  Let w=(c,d) be the black-blue vector between the first two starred points where  $c=(k+1)sin2x+(3k+3)sin4x+(4k+3)sin6x+...+(4k+3)sin(4n-4)x+(3k+3)sin(4n-2)x+(k+1)sin4nx$ and $d=4k+3+(7k+5)cos2x+(5k+3)cos4x+(4k+3)cos6x+...+(4k+3)cos(4n-4)x+(3k+3)cos(4n-2)x+(k+1)cos4nx$ and let v=(a,b) be the blue-black vector betweeen the last two starred points where $a=(k+1)sin4x+(2k+2)sin6x+...+(2k+2)sin(4n-2)x+(k+2)sin4nx$ and 

\noindent
$b=2k+2+(4k+4)cos2x+(3k+3)cos4x+(2k+2)cos6x+...+(2k+2)cos(4n-2)x+(k+2)cos4nx$, then $ad<bc$ and hence by the Non-Rhombus Tower Test 4 there is no poolshot since $ad<bc$
is equivalent to $sin4x+sin6x+...+sin(4n+2)x+sin(4n+4)x >0$ which holds since the code contains a 4n+4.

\doublespace
\singlespace

\noindent
\underline{...0*} black

\noindent
\underline{-2 0 2 4 ... 4n-2}

\noindent
$\left.
\begin{aligned}
&\text{4n 4n-2 ... 4 2 0 -2} \\
&\text{-4 -2 0 2 4 ... 4n-2} \\
&\vdots\quad\vdots\\
&\underline {\text{-4 -2 0 2 4 ... 4n-2}} \\
\end{aligned} 
\;\right\}$ 2k lines

\noindent
\underline {4n 4n-2 ... 4 2 0}

\noindent
$\left.
\begin{aligned}
&\text{-2 0 2 4 ... 4n-4} \\
&\text{4n-2 ... 6 4 2 0} \\
&\vdots\quad\vdots\\
&\underline {\text{4n-2 ... 6 4 2 0}} \\
\end{aligned} 
\;\right\}$ 2k lines

\noindent
\underline {-2 0 2 4 ... 4n-4 (4n-2)*} blue

\noindent
$\left.
\begin{aligned}
&\text{4n 4n-2 ... 6 4 2 0 -2} \\
&\text{-4 -2 0 2 ... 4n-4 4n-2} \\
&\vdots\quad\vdots\\
&\underline {\text{-4 -2 0 2 ... 4n-4 4n-2}} \\
\end{aligned} 
\;\right\}$ 2k+2 lines

\noindent
(4n)* black

In the special case of n=0 corresponding to the subcode ... 2 $4^{2k}$ 2 2 $4^{2k+2}$ ..., $w=(-(2k+2)sin2x-ksin4x,k+1+(2k+2)cos2x+kcos4x)$ and $v=(-(2k+2)sin2x-(k+1)sin4x,k+2+(2k+2)cos2x+(k+1)cos4x)$ and it follows that $ad<bc$ since this is equivalent to $sin4x>0$ which holds since 4 is one of the code numbers. Note if also k=0 this corresponds to the subcode ... 2 2 2 $4^{2}$ ... never appearing.

In the special case of n=1 corresponding to the subcode ... 6 $8^{2k}$ 6 $4^{2k}$ 6 $8^{2k+2}$ ..., $w=(sin2x+sin4x,4k+3+(6k+5)cos2x+(2k+1)cos4x)$ and $v=(sin4x,2k+2+(4k+4)cos2x+(2k+3)cos4x)$ and it follows that $ad<bc$ since this is equivalent to $sin4x+sin6x+sin8x>0$ which holds since 8 is one of the code numbers.

\noindent QED

\doublespace
\singlespace

\textbf{Cor 1:} Let S= 6 $4^{2k-2}$ 6 $8^{2k-2}$ and [S+4n]=4n+6 $(4n+4)^{2k-2}$ 4n+6 $(4n+8)^{2k-2}$ then in a rhombus pool shot  with $k\geq1$, $n\geq{-1}$ the subcode 

\doublespace
\noindent ... $(4n+8)^{2k-1}$ [S+4n] 4n+6 ...  (if it extends at least two spots to the left) 

forces

\noindent ... 4n+10 $(4n+8)^{2k-1}$ [S+4n] 4n+6 ...

\singlespace
Proof: Since ... 4n+6 $(4n+8)^{2k-1}$ 4n+6 ... is impossible by Rhombus Rule B and ... $(4n+8)^{2k}$ [S+4n] 4n+6 ... is impossible by Rhombus Rule N (reading it from right to left and replacing n by n+1 and k by k-1). This means the only remaining choice is ... 4n+10 $(4n+8)^{2k-1}$ [S+4n] 4n+6 ...

QED

\doublespace
\noindent
\textbf{Cor 2:} \textbf{(Rhombus Growth Rule 3)} It follows that  in a rhombus pool shot for $s\geq1$, $k\geq1$, $n\geq{-1}$

\noindent ... $(4n+8)^{2k-1}$ $[S+4n]^{s}$ 4n+6 ... (if it extends at least two spots to the left) 

 forces   

\noindent ... 4n+10 $(4n+8)^{2k-1}$ $[S+4n]^{s}$ 4n+6 ...

\singlespace

\noindent \textbf{Rhombus Rule O:} The subcode sequence $4^{2k+2}$ 6 $8^{2k}$ 6 $4^{2k}$ 6 for $k\geq0$ never appears in any rhombus poolshot. 

Note k=0 holds by Rhombus rule D (part a) using n=2 which means $4^2$ $6^3$ never appears.

Proof: The blue-black vector v=(a,b) between the first two starred points has $a=sin2x$ and $b=2k+2+(2k+1)cos2x$ whereas the black-blue vector w=(c,d) between the last two starred points has $c=sin2x+sin4x$ and $d=(4k+3)+(6k+5)cos2x+(2k+1)cos4x$.
But since $ad=bc$ and the blue-black-blue points are collinear then by the  Non-Rhombus Tower Collinear Test there is no poolshot.

\doublespace
\singlespace

\noindent
$\left.
\begin{aligned}
&\text{-2* 0} \\
&\text{2 0} \\
&\vdots\quad\vdots\\
&\underline {\text{2 0*}} \\
\end{aligned} 
\;\right\}$ 2k+2 lines

\noindent
\underline {-2 0 2}

\noindent
$\left.
\begin{aligned}
&\text{4 2 0 -2} \\
&\text{-4 -2 0 2} \\
&\vdots\quad\vdots\\
&\underline {\text{-4 -2 0 2}} \\
\end{aligned} 
\;\right\}$ 2k lines

\noindent
\underline {4 2 0}

\noindent
$\left.
\begin{aligned}
&\text{-2 0} \\
&\text{2 0} \\
&\vdots\quad\vdots\\
&\underline {\text{2 0}} \\
\end{aligned} 
\;\right\}$ 2k lines

\noindent
-2 0 2*

\doublespace
QED

\singlespace

\singlespace
\noindent \textbf{Rhombus Rule P:} The subcode sequences $(2n)^{2}$ 2n+2 $(2n+4)^{2}$ $(2n+2)^{2}$ or
$(2n+2)^{2}$ $(2n+4)^{2}$ 2n+2 $(2n)^{2}$ for $n\geq2$ never appear in any rhombus poolshot. Note for n=1, $2^2$ 4 $6^2$ $4^2$ can appear.

Proof: Using the first subcode above and the blue-black vector $v=(a,b)$ between the first two starred points where 
$a=2sin2x+4sin4x+5sin6x+...+5sin(2n-4)x+4sin(2n-2)x+2sin2nx$ and $b=5+8cos2x+6cos4x+5cos6x+...+5cos(2n-4)x+4cos(2n-2)x+2cos2nx$   and the black-blue vector $w=(c,d)$ between the last two starred points with $c=sin2x+2sin4x+2sin6x+...+2sin(2n-2)x$ and $d=2+3cos2x+2cos4x+...+2cos(2n-2)x$, we get $ad<bc$. This inequality is equivalent to $sin4x+sin6x+...+sin(2n-2)x>0$ which is true since $(2n+4)x<180$. Hence by the Non-Rhombus Tower Test 4, the subcode sequence never appears in a rhombus poolshot.

\doublespace
\singlespace

-2* 0 2 4 ... 2n-4 blue

\underline {2n-2 2n-4 ... 2 0}

\underline {-2 0 2 4 ... 2n-4 2n-2}

2n 2n-2 2n-4 ... 2 0 -2

\underline {-4 -2 0 2 4 ... 2n-4 2n-2}

2n* 2n-2 2n-4 ... 2 0 black

-2 0 2 4 ... 2n-4 (2n-2)* blue

\doublespace
\singlespace

Note in the special case $n=2$ the three starred blue-black-blue points are collinear where $v=(sin2x+sin4x,5+7cos2x+3cos4x)$ and $w=(sin2x,2+3cos2x)$ which is impossible by the Non-Rhombus Tower Collinear Test.
\noindent
QED

\doublespace
\singlespace
\noindent \textbf{Rhombus Rule Q:} The subcode ... $(2n+2)^{2}$ $(2n)^{2}$ 2n+2 $(2n+4)^{2}$ ... for $0\leq{n}$ never appears in any rhombus poolshot.

Proof: Consider the black-blue vector 

\noindent
$w=(c,d)=(2sin4x+2sin6x+...+2sin(2n-2)x+sin2nx,2+4cos2x+2cos4x+...+2cos(2n-2)x+cos2nx)$ between the first two starred points and the blue-black vector $v=(a,b)=(sin4x+2sin6x+...+2sin2nx,2+4cos2x+3cos4x+2cos6x+...+2cos2nx)$ between the last two starred points, then we would have ${ad}<{bc}$ which holds since this is equivalent to $sin(2n+4)x>0$ which would be true since 2n+4 is in the subcode. Hence by the Non-Rhombus Tower Test 5 the subcode is impossible.

\doublespace
\singlespace
\underline{... 0*} black

-2 0 2 4 ... 2n-2

\underline {2n 2n-2 ... 2 0}

-2* 0 2 ... 2n-4  blue

\underline {2n-2 2n-4 ... 2 0}

\underline {-2 0 2 ... 2n-4 (2n-2)*} blue

2n 2n-2 ... 2 0 -2

\underline {-4 -2 0 2 ... 2n-4 2n-2}

2n* ... black

\doublespace
\singlespace

\noindent Note the special cases.

1. n=0 which says that ... $2^{3}$ $4^{2}$ ... is impossible and this also follows from Rhombus Rule D part b. 

2. n=1 then  .. $4^{2}$ $2^{2}$ 4 $6^{2}$ ... is impossible since $w=(c,d)=(-sin2x,2+3cos2x)$ and $v=(a,b)=(-sin4x,2+4cos2x+cos4x)$ and $ad<bc$ since this is equivalent to $sin6x>0$ which holds since 6 is in the subcode.

3. n=2 then ... $6^2$ $4^2$ 6 $8^2$ ... is impossible since w=(c,d)=(sin4x,2+4cos2x+cos4x) and v=(a,b)=(sin4x,2+4cos2x+3cos4x) and $ad<bc$ since this is equivalent to $sin8x>0$ which holds since 8 is in the code.

QED

\doublespace
\singlespace

\noindent
\textbf{Forcing Rule 5:}

\doublespace
\singlespace

 \noindent Given a rhombus poolshot   and letting R=2n+2 $(2n)^{2k-2}$ 2n+2 $(2n+4)^{2k-2}$
then the subcode below (assuming that it extends to the left for at least 2k-1 spots)

\doublespace

\noindent 
... 2n+2 $(2n+4)^{2k-2}$ $R^{i}$ 2n+2 $(2n)^{2k-1}$ ...  with $k\geq2$, $i\geq0$ $n\geq2$ 

\doublespace
\noindent forces (Note: We will use $\downarrow$ to show that one subcode forces another.)

\noindent
... $(2n)^{2k-2}$ 2n+2 $(2n+4)^{2k-2}$ $R^{i}$ 2n+2 $(2n)^{2k-1}$ ...

\noindent Proof: We put the reasoning to the right of the forcing arrows.

\noindent
... 2n+2 $(2n+4)^{2k-2}$ $R^{i}$ 2n+2 $(2n)^{2k-1}$ ...

$\downarrow$ (2n+4 2n+2 2n+4 impossible by Rhombus Rule B )

 $\downarrow$  (For k=2, $n\geq2$,  $(2n+2)^{2}$ $(2n+4)^{2}$ 2n+2 $(2n)^{2}$ impossible by Rhombus Rule P)

$\downarrow$ ( For $k\geq3$, $(2n+2)^{2}$ $(2n+4)^{3}$  impossible by Rhombus Rule D )

\noindent ... 2n 2n+2 $(2n+4)^{2k-2}$ $R^{i}$ 2n+2 $(2n)^{2k-1}$ ...

$\downarrow$ Forcing Rule 2

\noindent ... $(2n)^2$ 2n+2 $(2n+4)^{2k-2}$ $R^{i}$ 2n+2 $(2n)^{2k-1}$ ...

$\downarrow$ Forcing Rule 1

\noindent ... $(2n)^3$ 2n+2 $(2n+4)^{2k-2}$ $R^{i}$ 2n+2 $(2n)^{2k-1}$ ...

$\downarrow$

\noindent
\vdots\quad\vdots\quad\vdots\quad\vdots

$\downarrow$ Forcing Rule 2

\noindent ... $(2n)^{2k-4}$ 2n+2 $(2n+4)^{2k-2}$ $R^{i}$ 2n+2 $(2n)^{2k-1}$ ...

$\downarrow$ Forcing Rule 1

\noindent ... $(2n)^{2k-3}$ 2n+2 $(2n+4)^{2k-2}$ $R^{i}$ 2n+2 $(2n)^{2k-1}$ ...

$\downarrow$ Forcing Rule 2

\noindent ... $(2n)^{2k-2}$ 2n+2 $(2n+4)^{2k-2}$ $R^{i}$ 2n+2 $(2n)^{2k-1}$ ...

QED

\singlespace
Note: 1. \textbf{(Rhombus Growth Rule 4)} This implies the growth rule that with i=0, $k\geq2$,  $n\geq{-1}$ and assuming that the first subcode extends to the left for at least 2k-1 spots.

\doublespace
\noindent 
... 4n+10 $(4n+12)^{2k-2}$ 4n+10 $(4n+8)^{2k-1}$ ...   

\doublespace
$\downarrow$ 

\noindent
... $(4n+8)^{2k-2}$ 4n+10 $(4n+12)^{2k-2}$ 4n+10 $(4n+8)^{2k-1}$ ... 

\singlespace
Note 2. \textbf{(Rhombus Growth Rule 5)} This also implies the growth rule that with $i=t\geq0$, $k\geq2$,  $n\geq{-1}$ where [S+4n+4]=4n+10 $(4n+8)^{2k-2}$ 4n+10 $(4n+12)^{2k-2}$ and assuming that the first subcode extends to the left for at least 2k-1 spots.

\doublespace
\noindent 
... 4n+10 $(4n+12)^{2k-2}$ $[S+4n+4]^t$ 4n+10 $(4n+8)^{2k-1}$ ...   

\doublespace
$\downarrow$ 

\noindent
... $(4n+8)^{2k-2}$ 4n+10 $(4n+12)^{2k-2}$ $[S+4n+4]^t$ 4n+10 $(4n+8)^{2k-1}$ ...   

\noindent
\textbf{Rhombus Rule R:}

The subcode $(2n+8)^{2k+1}$ 2n+6 $(2n+4)^{2k+1}$ 2n+2 $(2n)^{2k}$ 2n+2 $(2n+4)^{2k}$ 2n+2 ...  for $n\geq0$, $k\geq0$ never appears in a rhombus poolshot.

\singlespace
Proof: Considering the subcode from right to left and let $v=(a,b)=((k+1)sin4x+(3k+3)sin6x+(4k+3)sin8x+...+(4k+3)sin(2n-2)x+(3k+3)sin2nx+
(k+1)sin(2n+2)x,4k+3+(8k+6)cos2x+(7k+5)cos4x+(5k+3)cos6x+(4k+3)cos8x+...+(4k+3)cos(2n-2)x+
(3k+3)cos2nx+(k+1)cos(2n+2)x)$ be the blue-black vector between the second two starred points and $w=(c,d)=((k+1)sin2x+(3k+3)sin4x+(4k+3)sin6x+...+(4k+3)sin(2n-4)x+(3k+3)sin(2n-2)x+(k+1)sin2nx,4k+3+(7k+5)cos2x+(5k+3)cos4x+(4k+3)cos6x+...+(4k+3)cos(2n-4)x+(3k+3)cos(2n-2)x+(k+1)cos2nx)$ be the black-blue vector between the first two starred points. Then ad=bc and the points are collinear and hence by the Non-Rhombus Tower Collinear Test this subcode never appears in any rhombus pool shot.

\doublespace
\singlespace
\noindent
\underline {... 0*}

\noindent
\underline {-2 0 2 ... 2n-2}

\noindent
$\left.
\begin{aligned}
&\text{2n 2n-2 ... 2 0 -2}\\
&\quad\vdots\quad\vdots\\
&\underline {\text{-4 -2 0 2 ... 2n-4 2n-2}}
\end{aligned} 
\;\right\}$ 2k times

\noindent
\underline {2n 2n-2 ... 2 0}

\noindent
$\left.
\begin{aligned}
&\text{-2 0 2 ... 2n-4}\\
&\quad \vdots \quad \vdots\\
&\underline {\text{2n-2 2n-4 ... 2 0}}
\end{aligned} 
\;\right\}$ 2k times

\noindent
\underline {-2 0 2 ... 2n-4 (2n-2)*}

\noindent
$\left.
\begin{aligned}
&\text{2n 2n-2 ... 4 2 0 -2}\\
&\quad \vdots \quad \vdots\\
&\text{-4 -2 0 2 ... 2n-4 2n-2}\\
&\underline {\text{2n 2n-2 ... 4 2 0 -2}}
\end{aligned} 
\;\right\}$ 2k+1 times

\noindent
\underline {-4 -2 0 2 4 ... 2n-2 2n}

\noindent
$\left.
\begin{aligned}
&\text{2n+2 2n ... 4 2 0 -2 -4}\\
&\quad \vdots \quad \vdots\\
&\text{-6 -4 -2 0 2 4 ... 2n-2 2n}\\
&\underline {\text{2n+2 2n ... 4 2 0 -2 -4*}}
\end{aligned} 
\;\right\}$ 2k+1 times

\noindent
QED

\noindent Note the special cases in all of which $ad=bc$.

1. $n=0, w=(-(2k+2)sin2x-ksin4x,k+1+(2k+2)cos2x+kcos4x), v=(-(3k+2)sin2x-(3k+2)sin4x-ksin6x,3k+3+(5k+4)cos2x+(3k+2)cos4x+kcos6x)$

2. $n=1,  w=(-(2k+1)sin2x-ksin4x,3k+3+(4k+3)cos2x+kcos4x), v=(-ksin2x-(2k+1)sin4x-ksin6x,4k+3+(7k+6)cos2x+(4k+3)cos4x+kcos6x)$

3. $n=2,  w=(sin2x+sin4x,4k+3+(6k+5)cos2x+(2k+1)cos4x), v=(sin4x+sin6x,4k+3+(8k+6)cos2x+(6k+5)cos4x+(2k+1)cos6x)$

4. $n=3,  w=((k+1)sin2x+(2k+3)sin4x+(k+1)sin6x,4k+3+(7k+5)cos2x+(4k+3)cos4x+(k+1)cos6x), v=((k+1)sin4x+(2k+3)sin6x+(k+1)sin8x,4k+3+(8k+6)cos2x+(7k+5)cos4x+(4k+3)cos6x+(k+1)cos8x)$

\doublespace
\singlespace
\noindent \textbf{Cor 1: (Forcing Rule 6)} In a rhombus tower, for $n\geq0$, $k\geq0$, the subcode 

\doublespace
\noindent  ... $(2n+8)^{2k}$ 2n+6 $(2n+4)^{2k+1}$ 2n+2 $(2n)^{2k}$ 2n+2 $(2n+4)^{2k}$ 2n+2 ... 

$\downarrow$

\noindent  2n+6 $(2n+8)^{2k}$ 2n+6 $(2n+4)^{2k+1}$ 2n+2 $(2n)^{2k}$ 2n+2 $(2n+4)^{2k}$ 2n+2 ...

\singlespace \noindent Proof:  $(2n+8)^{2k+1}$ 2n+6 $(2n+4)^{2k+1}$ 2n+2 $(2n)^{2k}$ 2n+2 $(2n+4)^{2k}$ 2n+2 ... is impossible by Rhombus Rule R and if $k>0$, 2n+10 $(2n+8)^{2k}$ 2n+6 ... is impossible by Rhombus Rule C. If k=0, ... 2n+4 2n+6 2n+4 ... is impossible by Rhombus Rule B.

QED

\doublespace
\singlespace \noindent \textbf{Cor 2: (Rhombus Growth Rule 6)} In a rhombus tower, for $n\geq{-1}$, $k\geq1$, $s\geq1$ and S= 6 $4^{2k-2}$ 6 $8^{2k-2}$, the subcode

\doublespace \noindent ... $(4n+12)^{2k-2}$ 4n+10 $(4n+8)^{2k-1}$  $[S+4n]^s$ 4n+6 ...  

$\downarrow$

\noindent 4n+10 $(4n+12)^{2k-2}$ 4n+10 $(4n+8)^{2k-1}$  $[S+4n]^s$ 4n+6  ...

\noindent Proof: Replace n by 2n+2 and k by k-1 in Cor 1.

QED

\doublespace

\singlespace
\noindent
\textbf{Forcing Rule 7:} In a rhombus tower, for $n\geq{0}$, $k\geq0$, the subcode  (assuming it continues to the left for at least 2k+2 spots)

\doublespace
\noindent ... $(4n+4)^{2k+1}$ 4n+2 $(4n)^{2k}$ 4n+2 $(4n+4)^{2k}$ 4n+2 ... 

$\downarrow$  (... $(4n+4)^{2k+2}$ 4n+2 $(4n)^{2k}$ 4n+2 $(4n+4)^{2k}$ 4n+2 ... and
 
$\downarrow$    ... 4n+2  $(4n+4)^{2k+1}$ 4n+2... are impossible by Rhombus Rule N and B)

\noindent ... 4n+6 $(4n+4)^{2k+1}$ 4n+2 $(4n)^{2k}$ 4n+2 $(4n+4)^{2k}$ 4n+2 ...

$\downarrow$ Rhombus Rule F

\noindent ... $(4n+8)^{2k}$ 4n+6 $(4n+4)^{2k+1}$ 4n+2 $(4n)^{2k}$ 4n+2 $(4n+4)^{2k}$ 4n+2 ...

$\downarrow$ Forcing Rule 6 with n replaced by 2n

\noindent 4n+6 $(4n+8)^{2k}$ 4n+6 $(4n+4)^{2k+1}$ 4n+2 $(4n)^{2k}$ 4n+2 $(4n+4)^{2k}$ 4n+2 ...

\doublespace
\singlespace
\noindent
\textbf{Rhombus Rule S:} For $n\geq0$, $k\geq0$ the subcode

... $(2n+4)^{2k+1}$ 2n+6 $(2n+8)^{2k}$ 2n+6 $(2n+4)^{2k+1}$ 2n+2 $(2n)^{2k}$ 2n+2 $(2n+4)^{2k}$ 2n+2 $(2n)^{2k}$ 2n+2 $(2n+4)^{2k}$ 2n+2 ... never appears in a rhombus poolshot.

Proof:

\doublespace
\singlespace

\noindent
\underline{...0*}
\singlespace

\noindent
$\left.
\begin{aligned}
&\text{-2 0 2 4 ... 2n}\\
&\text{2n+2 2n ...  2 0}\\
&\quad\vdots\quad\vdots\\
&\underline {\text{-2 0 2 4 ... 2n}}
\end{aligned} 
\;\right\}$ 2k+1 times

\noindent
\underline {2n+2 2n ... 2 0 -2}

\noindent
$\left.
\begin{aligned}
&\text{-4 -2 0 2 4 ... 2n+2}\\
&\text{2n+4 2n+2 ... 2 0 -2}\\
&\quad\vdots\quad\vdots\\
&\underline {\text{2n+4 2n+2 ... 2 0 -2}}
\end{aligned} 
\;\right\}$ 2k times

\noindent
\underline {-4 -2 0 ... 2n-2 2n}

\noindent
$\left.
\begin{aligned}
&\text{2n+2 2n ... 2 0}\\
&\text{-2 0 2 ... 2n}\\
&\quad\vdots\quad\vdots\\
&\underline {\text{2n+2 2n ... 2 0}}
\end{aligned} 
\;\right\}$ 2k+1 times

\noindent
\underline {-2* 0 ... 2n-2}

\noindent
$\left.
\begin{aligned}
&\text{2n ... 4 2}\\
&\text{0 2 ... 2n-2}\\
&\quad\vdots\quad\vdots\\
&\underline {\text{0 2 ... 2n-2}}
\end{aligned} 
\;\right\}$ 2k times

\noindent
\underline {2n ... 6 4 2 0}

\noindent
$\left.
\begin{aligned}
&\text{-2 0 2 ... 2n}\\
&\text{2n+2 ... 2 0}\\
&\quad\vdots\quad\vdots\\
&\underline {\text{2n+2 ... 2 0}}
\end{aligned} 
\;\right\}$ 2k times

\noindent
\underline {-2 0 ... 2n-2}

\noindent
$\left.
\begin{aligned}
&\text{2n ... 4 2}\\
&\text{0 2 ... 2n-2}\\
&\quad\vdots\quad\vdots\\
&\underline {\text{0 2 ... 2n-2}}
\end{aligned} 
\;\right\}$ 2k times

\noindent
\underline {2n ... 6 4 2 0}

\noindent
$\left.
\begin{aligned}
&\text{-2 0 2 4 ... 2n}\\
&\text{2n+2 ... 4 2 0}\\
&\quad\vdots\quad\vdots\\
&\underline {\text{2n+2 ... 4  2 0}}
\end{aligned} 
\;\right\}$ 2k times

\noindent
\underline{-2 0 2 ...2n-2}

\noindent
2n*

\doublespace
\singlespace

Let $w=(c,d)=(2ksin2x+(5k+3)sin4x+(6k+4)sin6x+...+(6k+4)sin2nx+(4k+2)sin(2n+2)x+ksin(2n+4)x,
6k+4+(10k+8)cos2x+(7k+5)cos4x+(6k+4)cos6x+...+(6k+4)cos2nx+(4k+2)cos(2n+2)x+kcos(2n+4)x)$ be the black-blue vector from the first starred point  to the second starred point and $v=(a,b)=((6k+3)sin2x+(8k+5)sin4x+...+(8k+5)sin(2n-2)x+(6k+3)sin(2n)x+2ksin(2n+2)x, 6k+5+(10k+7)cos2x+(8k+5)cos4x+...+(8k+5)cos(2n-2)x+(6k+3)cos(2n)x+2kcos(2n+2)x)$ be the blue-black vector from the second starred point to the third starred point. Then $ad<bc$ which is equivalent to $sin4x+sin6x+...+sin(2n+4)x>0$  which holds since $(2n+8)x<180$.
 Hence by the Non-Rhombus Tower Test 4 the given subcode does not appear in a rhombus poolshot.

\noindent QED.

\noindent Note the special cases.

1. $n=0, w=(-2sin2x-sin4x,6k+4+(8k+6)cos2x+(2k+1)cos4x), v=(-2sin2x,4k+3+(4k+2)cos2x)$ and $ad<bc$ is equivalent to $sin4x>0$ which holds since $8x<180$.

2. $n=1, w=(2ksin2x+(3k+1)sin4x+ksin6x,6k+4+(10k+8)cos2x+(5k+3)cos4x+kcos6x), v=((4k+1)sin2x+2ksin4x,6k+5+(8k+5)cos2x+2kcos4x)$ and $ad<bc$ is equivalent to $sin4x+sin6x>0$ which holds since $10x<180$.

\doublespace

\singlespace
\noindent 
\textbf{Cor:} \textbf{(Rhombus Growth Rule 7)} Let  [S+4n] = 4n+6 $(4n+4)^{2k-2}$ 4n+6 $(4n+8)^{2k-2}$, then in a rhombus poolshot the subcode for $n\geq{-1}$, $k\geq1$ (assuming it continues to the left at least two spots)

\doublespace
\noindent
... $(4n+8)^{2k-2}$ 4n+10 $(4n+12)^{2k-2}$ 4n+10 $(4n+8)^{2k-1}$ $[S+4n]^2$ 4n+6 ... 

 $\downarrow$ 

\noindent
4n+10 $(4n+8)^{2k-2}$ 4n+10 $(4n+12)^{2k-2}$ 4n+10 $(4n+8)^{2k-1}$ $[S+4n]^2$ 4n+6 ...

\singlespace
Proof: ... $(4n+8)^{2k-1}$ 4n+10 $(4n+12)^{2k-2}$ 4n+10 $(4n+8)^{2k-1}$ $[S+4n]^2$ 4n+6 ...  is impossible by Rhombus Rule S  and ... 4n+6 $(4n+8)^{2k-2}$ 4n+10 is impossible by Rhombus Rule C.

\noindent QED.

\doublespace
\noindent
\textbf{Rhombus Rule T:}

\noindent
Let R = 2n+2 $(2n)^{2k-2}$ 2n+2 $(2n+4)^{2k-2}$ then for $n\geq0$, $k\geq1$, the subcode

 $(2n+8)^{2k-1}$ $[R+4]^{i-1}$ 2n+6 $(2n+4)^{2k-1}$ $R^i$ 2n+2 ...

\singlespace
\noindent
never appears in a rhombus poolshot for $i\geq1$ where [R+4]=2n+6 $(2n+4)^{2k-2}$ 2n+6 $(2n+8)^{2k-2}$

\doublespace
\singlespace
\noindent
Proof: Considering the subcode from right to left, let $w=(c,d)=(kisin2x+(3ki-i+1)sin4x+(4ki-2i+1)sin6x+ ...+(4ki-2i+1)sin(2n-4)x+(3ki-i+1)sin(2n-2)x+kisin2nx,4ki-2i+1+(7ki-4i+2)cos2x+(5ki-3i+1)cos4x+(4ki-2i+1)cos6x + ...+(4ki-2i+1)cos(2n-4)x+(3ki-i+1)cos(2n-2)x+kicos2nx)$ be a black-blue vector between the first two starred points and let $v=(a,b)=(kisin4x+(3ki-i+1)sin6x+(4ki-2i+1)sin8x+ ... +(4ki-2i+1)sin(2n-2)x+(3ki-i+1)sin2nx+kisin(2n+2)x,4ki-2i+1+(8ki-4i+2)cos2x+(7ki-4i+2)cos4x+(5ki-3i+1)cos6x+(4ki-2i+1)cos8x+ ... +(4ki-2i+1)cos(2n-2)x+(3ki-i+1)cos2nx+kicos(2n+2)x)$ be a blue-black vector between the second two starred points. The three starred C points are a collinear black-blue-black situation since ad=bc and hence by the Non-Rhombus Tower Collinear Test never appear  in a rhombus poolshot.

\doublespace

\singlespace
\noindent
\underline 0*

\noindent
$\left.
\parbox{0.5\linewidth}
{\underline {-2 0 2 ... 2n-2}

\noindent
$\left.
\begin{aligned}
&\text{2n 2n-2 ... 2 0 -2}\\
&\quad\vdots\quad\vdots\\
&\underline {\text{-4 -2 0 2 ... 2n-2}}
\end{aligned} 
\;\right\}$ $2k-2$ \text{ times}

\underline {2n 2n-2 ... 4 2 0}

\noindent
$\left.
\begin{aligned}
&\text{-2 0 2 ... 2n-4}\\
&\quad \vdots \quad \vdots\\
&\underline {\text{2n-2 ... 4 2 0}}
\end{aligned} 
\;\right\}$ $2k-2$ \text{ times}

\underline {-2 0 2 ... 2n-4 (2n-2)*}
}
\right\}$ i times

\noindent
$\left.
\begin{aligned}
&\text{2n 2n-2 ... 2 0 -2}\\
&\text{-4 -2 0 2 ... 2n-2}\\
&\quad\vdots\quad\vdots\\
&\underline {\text{2n 2n-2 ... 2 0 -2}}
\end{aligned} 
\;\right\}$ 2k-1 times

\noindent
$\left.
\parbox{0.5\linewidth}
{\underline {-4 -2 0 2 ...2n-2 2n}

\noindent
$\left.
\begin{aligned}
&\text{2n+2 2n ... 2 0 -2 -4}\\
&\text{-6 -4 -2 0  ... 2n-2 2n}\\
&\quad\vdots\quad\vdots\\
&\underline {\text{-6 -4 -2 0  ... 2n-2 2n}}
\end{aligned} 
\;\right\}$ $2k-2$ \text{ times}

\underline {2n+2 2n ... 2 0 -2}

\noindent
$\left.
\begin{aligned}
&\text{-4 -2 0 2 ... 2n-2}\\
&\text{2n 2n-2 ... 2 0 -2}\\
&\quad \vdots \quad \vdots\\
&\underline {\text{2n 2n-2 ... 2 0 -2}}
\end{aligned} 
\;\right\}$ $2k-2$ \text{ times}

\underline {-4 -2 0 2 ... 2n-2 2n}
}
\right\}$ i-1 times

\noindent
$\left.
\begin{aligned}
&\text{2n+2 2n 2n-2 ... 2 0 -2 -4}\\
&\text{-6 -4 -2 0 2 ... 2n-4 2n-2 2n}\\
&\quad\vdots\quad\vdots\\
&\underline {\text{2n+2 2n 2n-2 ... 2 0 -2 -4*}}
\end{aligned} 
\;\right\}$ 2k-1 times

\doublespace
\noindent 
QED.

\noindent Note the special cases.

\singlespace
1. $n=0, w=(-(2ki-i+1)sin2x-(k-1)isin4x,ki+(2ki-i+1)cos2x+(k-1)icos4x), v=(-(3ki-2i+1)sin2x-(3ki-2i+1)sin4x-(k-1)isin6x,3ki-i+1+(5ki-2i+1)cos2x+(3ki-2i+1)cos4x+(k-1)icos6x)$ and $ad=bc$.

2. $n=1, w=(-(2ki-1)sin2x-(k-1)isin4x,3ki-i+1+(4ki-2i+1)cos2x+(k-1)icos4x), v=(-(k-1)isin2x-(2ki-1)sin4x-(k-1)isin6x,4ki-2i+1+(7ki-3i+2)cos2x+(4ki-2i+1)cos4x+(k-1)icos6x)$ and $ad=bc$.

3. $n=2, w=(isin2x+isin4x,4ki-2i+1+(6ki-3i+2)cos2x+(2ki-i)cos4x), v=(isin4x+isin6x,4ki-2i+1+(8ki-4i+2)cos2x+(6ki-3i+2)cos4x+(2ki-i)cos6x)$ and $ad=bc$.

4. $n=3, w=(kisin2x+(2ki+1)sin4x+kisin6x,4ki-2i+1+(7ki-4i+2)cos2x+(4ki-2i+1)cos4x+kicos6x), v=(kisin4x+(2ki+1)sin6x+kisin8x,4ki-2i+1+(8ki-4i+2)cos2x+(7ki-4i+2)cos4x+(4ki-2i+1)cos6x+kicos8x)$ and $ad=bc$.

\doublespace
\noindent
\textbf{Cor 1:}

\noindent
Let R = 2n+2 $(2n)^{2k-2}$ 2n+2 $(2n+4)^{2k-2}$ then for $n\geq0$, $k\geq1$

 $(2n+8)^{2k-1}$ $[R+4]^{s}$ 2n+6 $(2n+4)^{2k-1}$ $R^i$ 2n+2 ...

\noindent
never appears in a rhombus poolshot where $0\leq{s}<i$.

\singlespace
Proof: Since

\noindent
 $(2n+8)^{2k-1}$ $[R+4]^{s}$ 2n+6 $(2n+4)^{2k-1}$ $R^{s+1}$ 2n+2 ...

\noindent
never appears in a rhombus poolshot by Rhombus Rule T.

QED

\doublespace
\noindent \textbf{COR 2: (Forcing Rule 8)} For $n\geq0$, $k\geq1$, $0\leq{s}<i$

\noindent ... $(2n+8)^{2k-2}$ $[R+4]^{s}$ 2n+6 $(2n+4)^{2k-1}$ $R^i$ 2n+2 ... 

 $\downarrow$

\noindent 2n+6 $(2n+8)^{2k-2}$ $[R+4]^{s}$ 2n+6 $(2n+4)^{2k-1}$ $R^i$ 2n+2 ...

\singlespace
\noindent Proof: Since  $(2n+8)^{2k-1}$ $[R+4]^{s}$ 2n+6 $(2n+4)^{2k-1}$ $R^i$ 2n+2 ... is impossible by corollary one
and 2n+10 $(2n+8)^{2k-2}$ 2n+6 ... is impossible by Rhombus Rule C.

QED

\doublespace
\noindent This leads to the Growth Rule where if S= 6 $4^{2k-2}$ 6 $8^{2k-2}$  then 

\singlespace
\noindent \textbf{Rhombus Growth Rule 8:} In a rhombus tower, for $k\geq1$, $0\leq{s}<i$, $n\geq-1$ the subcode

\doublespace

\noindent ... $(4n+12)^{2k-2}$ $[S+4n+4]^s$ 4n+10 $(4n+8)^{2k-1}$ $[S+4n]^i$  4n+6 ...

 $\downarrow$ 

\noindent 4n+10 $(4n+12)^{2k-2}$ $[S+4n+4]^s$ 4n+10 $(4n+8)^{2k-1}$ $[S+4n]^i$  4n+6 ...

\singlespace
\noindent Proof: In Forcing Rule 8, replace n by 2n+2.

QED

\doublespace
\noindent
\textbf{Rhombus Rule U:}

\doublespace

\noindent Let R = 2n+2 $(2n)^{2k-2}$ 2n+2 $(2n+4)^{2k-2}$ then for $n\geq0$, $k\geq1$ the subcode

\noindent
... $(2n+4)^{2k-1}$ 2n+6 $(2n+8)^{2k-2}$ $[R+4]^j$ 2n+6 $(2n+4)^{2k-1}$ $R^i$ 2n+2 ...

\noindent
is impossible for $j\geq0$ and $i\geq(j+2)$.

\singlespace
\noindent
Proof: Considering the subcode from right to left, it is enough to show it for $i=j+2$. Let $w=(c,d)=((kj+2k)sin2x+(3kj-j+6k-1)sin4x+(4kj-2j+8k-3)sin6x+...+(4kj-2j+8k-3)sin(2n-4)x+(3kj-j+6k-1)sin(2n-2)x+(kj+2k)sin2nx,4kj-2j+8k-3+(7kj-4j+14k-6)cos2x+(5kj-3j+10k-5)cos4x+(4kj-2j+8k-3)cos6x+...+(4kj-2j+8k-3)cos(2n-4)x+(3kj-j+6k-1)cos(2n-2)x+(kj+2k)cos2nx)$ be a black-blue vector between the first two starred points and let $v=(a,b)=((kj+2k)sin4x+(3kj-j+5k-1)sin6x+(4kj-2j+6k-2)sin8x+...+(4kj-2j+6k-2)sin(2n-2)x+(3kj-j+4k)sin(2n)x+(kj+k)sin(2n+2)x,4kj-2j+6k-2+(8kj-4j+12k-4)cos2x+(7kj-4j+10k-4)cos4x+(5kj-3j+7k-3)cos6x+(4kj-2j+6k-2)cos8x+...+(4kj-2j+6k-2)cos(2n-2)x+(3kj-j+4k)cos(2n)x+(kj+k)cos(2n+2)x)$ be a blue-black vector between the second two starred points. Then $ad<bc$ is equivalent to $sin4x+sin6x+...sin(2n+4)x>0$ and this holds since $(2n+8)x>0$. Hence by the Non-Rhombus Tower Test 4 this subcode is not part of a rhombus poolshot.

\doublespace

\singlespace
\noindent
\underline {... 0*}

\noindent
$\left.
\parbox{0.5\linewidth}
{\underline {-2 0 2 ... 2n-2}

\noindent
$\left.
\begin{aligned}
&\text{2n 2n-2 ... 2 0 -2}\\
&\quad\vdots\quad\vdots\\
&\underline {\text{-4 -2 0 2 ... 2n-2}}
\end{aligned} 
\;\right\}$ $2k-2$ \text{ times}

\underline {2n 2n-2 ... 4 2 0}

\noindent
$\left.
\begin{aligned}
&\text{-2 0 2 ... 2n-4}\\
&\quad \vdots \quad \vdots\\
&\underline {\text{2n-2 ... 4 2 0}}
\end{aligned} 
\;\right\}$ $2k-2$ \text{ times}

\underline {-2 0 2 ... 2n-4 (2n-2)*}
}
\right\}$ i times

\noindent
$\left.
\begin{aligned}
&\text{2n 2n-2 ... 2 0 -2}\\
&\text{-4 -2 0 2 ... 2n-2}\\
&\quad\vdots\quad\vdots\\
&\underline {\text{2n 2n-2 ... 2 0 -2}}
\end{aligned} 
\;\right\}$ 2k-1 times

\noindent
$\left.
\parbox{0.5\linewidth}
{\underline {-4 -2 0 2 ...2n-2 2n}

\noindent
$\left.
\begin{aligned}
&\text{2n+2 2n ... 2 0 -2 -4}\\
&\text{-6 -4 -2 0  ... 2n-2 2n}\\
&\quad\vdots\quad\vdots\\
&\underline {\text{-6 -4 -2 0  ... 2n-2 2n}}
\end{aligned} 
\;\right\}$ $2k-2$ \text{ times}

\underline {2n+2 2n ... 2 0 -2}

\noindent
$\left.
\begin{aligned}
&\text{-4 -2 0 2 ... 2n-2}\\
&\text{2n 2n-2 ... 2 0 -2}\\
&\quad \vdots \quad \vdots\\
&\underline {\text{2n 2n-2 ... 2 0 -2}}
\end{aligned} 
\;\right\}$ $2k-2$ \text{ times}

\underline {-4 -2 0 2 ... 2n-2 2n}
}
\right\}$ j times

\noindent
$\left.
\begin{aligned}
&\text{2n+2 2n 2n-2 ... 0 -2 -4}\\
&\text{-6 -4 -2 0 2 ... 2n-2 2n}\\
&\quad\vdots\quad\vdots\\
&\underline {\text{-6 -4 -2 0 2 ... 2n-2 2n}}
\end{aligned} 
\;\right\}$ $2k-2$ times

\noindent
\underline {2n+2 2n ... 2 0 -2}

\noindent
$\left.
\begin{aligned}
&\text{-4 -2 0 2 ... 2n-2}\\
&\text{2n 2n-2 ... 2 0 -2}\\
&\quad \vdots \quad \vdots\\
&\underline {\text{-4 -2 0 2 ... 2n-2}}
\end{aligned} 
\;\right\}$ $2k-1$ \text{ times}

\noindent
2n* ...

\noindent
QED.

\noindent Note the special cases using $i=j+2$.

1. $n=0, w=((-2kj+j-4k+1)sin2x+(-kj+j-2k+2)sin4x,kj+2k+(2kj-j+4k-1)cos2x+(kj-j+2k-2)cos4x), v=((-3kj+2j-5k+2)sin2x+(-3kj+2j-4k+2)sin4x+(-kj+j-k+1)sin6x,3kj-j+4k+(5kj-2j+7k-2)cos2x+(3kj-2j+4k-2)cos4x+(kj-j+k-1)cos6x)$ with $ad<bc$ equivalent to $sin4x>0$.

2. $n=1, w=((-2kj+2j-4k+3)sin2x+(-kj+j-2k+2)sin4x,3kj-j+6k-1+(4kj-2j+8k-3)cos2x+(kj-j+2k-2)cos4x), v=((-kj+j-2k+2)sin2x+(-2kj+2j-3k+2)sin4x+(-kj+j-k+1)sin6x,4kj-2j+6k-2+(7kj-3j+10k-2)cos2x+(4kj-2j+5k-2)cos4x+(kj-j+k-1)cos6x)$ with $ad<bc$ equivalent to $sin4x+sin6x>0$.

3. $n=2, w=((j+2)sin2x+(j+2)sin4x,4kj-2j+8k-3+(6kj-3j+12k-4)cos2x+(2kj-j+4k-2)cos4x), v=((j+2)sin4x+(j+1)sin6x,4kj-2j+6k-2+(8kj-4j+12k-4)cos2x+(6kj-3j+8k-2)cos4x+(2kj-j+2k-1)cos6x)$ with $ad<bc$ equivalent to $sin4x+sin6x+sin8x>0$.

4. $n=3, w=((kj+2k)sin2x+(2kj+4k+1)sin4x+(kj+2k)sin6x,4kj-2j+8k-3+(7kj-4j+14k-6)cos2x+(4kj-2j+8k-3)cos4x+(kj+2k)cos6x), v=((kj+2k)sin4x+(2kj+3k+1)sin6x+(kj+k)sin8x,4kj-2j+6k-2+(8kj-4j+12k-4)cos2x+(7kj-4j+10k-4)cos4x+(4kj-2j+5k-1)cos6x+(kj+k)cos8x)$ with $ad<bc$ equivalent to $sin4x+sin6x+sin8x+sin10x>0$.

\doublespace
\singlespace
\noindent \textbf{COR 1: (Forcing Rule 9)}
\noindent Let R = 2n+2 $(2n)^{2k-2}$ 2n+2 $(2n+4)^{2k-2}$ then for $n\geq0$, $k\geq1$, $j\geq0$ and $i\geq(j+2)$ the subcode (assuming it extends at least two spots to the left)

\doublespace
\noindent
... $(2n+4)^{2k-2}$ 2n+6 $(2n+8)^{2k-2}$ $[R+4]^j$ 2n+6 $(2n+4)^{2k-1}$ $R^i$ 2n+2 ...

$\downarrow$ 

\noindent
...2n+6 $(2n+4)^{2k-2}$ 2n+6 $(2n+8)^{2k-2}$ $[R+4]^j$ 2n+6 $(2n+4)^{2k-1}$ $R^i$ 2n+2...

\singlespace
\noindent Proof: Since ... $(2n+4)^{2k-1}$ 2n+6 $(2n+8)^{2k-2}$ $[R+4]^j$ 2n+6 $(2n+4)^{2k-1}$ $R^i$ 2n+2 ... is impossible by Rhombus Rule U and if $k>1$ ... 2n+2 $(2n+4)^{2k-2}$ 2n+6 is impossible by Rhombus Rule C. If $k=1$, 2n+8  $(2n+6)^{2j+2} 2n+4$ ... is also impossible by Rhombus Rule C.

QED

\doublespace
\singlespace
\noindent This leads to the Growth Rule where if S= 6 $4^{2k-2}$ 6 $4^{2k-2}$ then

\singlespace
\noindent \textbf{COR 2:} \textbf{(Rhombus Growth Rule 9)} In a rhombus tower, for $k\geq2$, $3\leq{t+2}\leq s$, the subcode (assuming it extends at least two spots to the left)

\doublespace
\singlespace
\noindent ... $(4n+8)^{2k-2}$ 4n+10 $(4n+12)^{2k-2}$ $[S+4n+4]^t$ 4n+10 $(4n+8)^{2k-1}$ $[S+4n]^s$  4n+6 ...

\doublespace
$\downarrow$ 

\singlespace
\noindent ... 4n+10 $(4n+8)^{2k-2}$ 4n+10 $(4n+12)^{2k-2}$ $[S+4n+4]^t$ 4n+10 $(4n+8)^{2k-1}$ $[S+4n]^s$  4n+6 ...

\doublespace
\singlespace
\noindent Proof: In Forcing Rule 9, replace n by 2n+2.

QED

\subsection*{14. More Corridor Results}

\singlespace
\textbf{Corridor Lemma 5:} The code sequence of a corridor rhombus tower never ends in the subcode $4^{2s+2}$ 6 $8^{2s}$ 6 $4^{2s+1}$ 2 for $s\geq0$.

Proof: The three starred points represent a blue-black-black collinear situation since if $v=(a,b)$ is a vector between the first two starred points with $a=sinx$ and $b=2s+2+(2s+1)cos2x$ and if $w=(c,d)$ is a vector between the last two starred points with $c=sin2x+sin4x$ and  $d= (4s+3)+(6s+5)cos2x+(2s+1)cos4x$ then $ad=bc$. Hence this is not the code sequence of a corridor by the Non-Corridor Test 2. Note that the last starred point is $C_n$.

\doublespace
\singlespace

\noindent
$\left.
\begin{aligned}
&\text{-2* 0} \\
&\text{2 0} \\
&\vdots\quad\vdots\\
&\underline {\text{2 0*}} \\
\end{aligned} 
\;\right\}$ 2s+2 lines

\noindent
\underline {-2 0 2}

\noindent
$\left.
\begin{aligned}
&\text{4 2 0 -2} \\
&\text{-4 -2 0 2} \\
&\vdots\quad\vdots\\
&\underline {\text{-4 -2 0 2}} \\
\end{aligned} 
\;\right\}$ 2s lines

\noindent
\underline {4 2 0}

\noindent
$\left.
\begin{aligned}
&\text{-2 0} \\
&\text{2 0} \\
&\vdots\quad\vdots\\
&\underline {\text{-2 0}} \\
\end{aligned} 
\;\right\}$ 2s+1 lines

\noindent
2*

QED

\doublespace\noindent
\textbf{Corridor Lemma 6:}

\singlespace
Let S= 6 $4^{2k-2}$ 6 $8^{2k-2}$ and let ... 6 $4^{2k-1}$ 2 be the end of the code of a corridor rhombus tower, then for $i\geq0$, $n\geq1$, $k\geq1$ the corridor tower never ends with the subcode

\doublespace
\singlespace

... $(4n+4)^{2k-1}$ 4n+6 $(4n+8)^{2k-2}$ $[S+4n]^{i}$ 4n+6

$(4n+4)^{2k-1}$ $[S+4n-4]^{i+1}$ 4n+2

\quad \vdots \quad \vdots

$12^{2k-1}$ $[S+4]^{i+1}$ 10

$8^{2k-1}$ $[S]^{i+1}$ 6 $4^{2k-1}$ 2

\doublespace
\singlespace
\noindent
Proof: Note, for convenience we  have written the corresponding vertical array in reverse of the usual order so that 0* corresponds to the last C point $C_{n}$ and 0** corresponds to the special point $C_{m}$. \textbf{Backwards standard position} would be to have the base $A_{n}$$B_{n}$ horizontal with A to the right of B and $C_{n}$ above the base and with the corridor extending upwards. Observe that a corridor in backwards standard position will have all B points and  black C points on the left side and all A points and blue C points on the right side of the corridor and will still lean to the right.  Further for a corridor in backwards standard position, the integers in the blue lines now decrease while the integers in the black lines now increase. So in the array below, the * points are black and the ** points are blue. 

\doublespace
Figure 29

\doublespace
\singlespace
\noindent
\underline{0*}

\singlespace
\noindent
$\left.
\begin{aligned}
&\text{2 0}\\
&\text{-2 0}\\
&\quad \vdots \quad \vdots\\
&\underline {\text{2 0**}}
\end{aligned} 
\;\right\}$ $2k-1$ \text{ times}

\noindent
$\left.
\parbox{0.5\linewidth}
{\underline {-2 0 2}

\noindent
$\left.
\begin{aligned}
&\text{4 2 0 -2}\\
&\quad\vdots\quad\vdots\\
&\underline {\text{-4 -2 0 2}}
\end{aligned} 
\;\right\}$ $2k-2$ \text{ times}

\underline {4 2 0}

\noindent
$\left.
\begin{aligned}
&\text{-2 0}\\
&\quad \vdots \quad \vdots\\
&\underline {\text{2 0}}
\end{aligned} 
\;\right\}$ $2k-2$ \text{ times}

\underline {-2 0 2}
}
\right\}$ i+1 times

\noindent
$\left.
\begin{aligned}
&\text{4 2 0 -2}\\
&\text{-4 -2 0 2}\\
&\quad \vdots \quad \vdots\\
&\underline {\text{4 2 0 -2}}
\end{aligned} 
\;\right\}$ $2k-1$ \text{ times}

\noindent
$\left.
\parbox{0.5\linewidth}
{\underline {-4 -2 0 2 4}

\noindent
$\left.
\begin{aligned}
&\text{6 4 2 0 -2 -4}\\
&\quad\vdots\quad\vdots\\
&\underline {\text{-6 -4 -2 0 2 4}}
\end{aligned} 
\;\right\}$ $2k-2$ \text{ times}

\underline {6 4 2 0 -2}

\noindent
$\left.
\begin{aligned}
&\text{-4 -2 0 2}\\
&\quad \vdots \quad \vdots\\
&\underline {\text{4 2 0 -2}}
\end{aligned} 
\;\right\}$ $2k-2$ \text{ times}

\underline {-4 -2 0 2 4}
}
\right\}$ i+1 times

\quad\vdots\quad\vdots

\quad\vdots\quad\vdots

\quad\vdots\quad\vdots

\noindent
$\left.
\parbox{0.5\linewidth}
{\underline {-2n -(2n-2) ...-2 0 2 ... 2n-2 2n}

\noindent
$\left.
\begin{aligned}
&\text{2n+2 2n ... 0 -2 ... -2n}\\
&\quad\vdots\quad\vdots\\
&\underline {\text{-(2n+2) -2n ... 0 2 ... 2n}}
\end{aligned} 
\;\right\}$ 2k-2 times

\underline {2n+2 2n ... 2 0 -2 ... -(2n-2)}

\noindent
$\left.
\begin{aligned}
&\text{-2n -(2n-2) ... 0 2 ... 2n-2}\\
&\quad \vdots \quad \vdots\\
&\underline {\text{2n 2n-2 ... 0 -2 ... -(2n-2)}}
\end{aligned} 
\;\right\}$2k-2 times

\underline {-2n -(2n-2) ... -2 0 2 ... 2n*}
}
\right\}$ i+1 times

\noindent
$\left.
\begin{aligned}
&\text{2n+2 2n ... 2 0 -2 ... -(2n)}\\
&\text{-(2n+2) -2n ... 0 2 ... 2n}\\
&\quad \vdots \quad \vdots\\
&\underline {\text{2n+2 2n ... 2 0 -2 ... -(2n)}}
\end{aligned} 
\;\right\}$ $2k-1$ \text{ times}

\noindent
$\left.
\parbox{0.7\linewidth}
{\underline {-(2n+2) -2n ...-2 0 2 ... 2n 2n+2}

\noindent
$\left.
\begin{aligned}
&\text{2n+4 2n+2 ... 0 -2 ... -2n -(2n+2)}\\
&\quad\vdots\quad\vdots\\
&\underline {\text{-(2n+4) -(2n+2) ... 0 2 ... 2n+2}}
\end{aligned} 
\;\right\}$2k-2 times

\underline {2n+4 2n+2 2n ... 2 0 -2 ... -(2n)}

\noindent
$\left.
\begin{aligned}
&\text{-(2n+2) -2n ... 0 2 ... 2n}\\
&\quad \vdots \quad \vdots\\
&\underline {\text{2n+2 2n ... 0 -2 ... -(2n)}}
\end{aligned} 
\;\right\}$2k-2 times

\underline {-(2n+2) -2n ... -2 0 2 ... 2n+2}
}
\right\}$ i times

\noindent
$\left.
\begin{aligned}
&\text{2n+4 2n+2 ... 2 0 -2 ... -(2n+2)}\\
&\text{-(2n+4) -(2n+2) ...-2  0 2 ... 2n+2}\\
&\quad \vdots \quad \vdots\\
&\underline {\text{-(2n+4) -(2n+2) ...-2  0 2 ... 2n+2}}
\end{aligned} 
\;\right\}$ $2k-2$ \text{ times}

\noindent
\underline {2n+4 2n+2 ... 2 0 -2 ... -(2n)}

\noindent
$\left.
\begin{aligned}
&\text{-(2n+2) -2n ...-2  0 2 ... 2n}\\
&\text{2n+2 2n ... 2 0 -2 ... -(2n)}\\
&\quad \vdots \quad \vdots\\
&\underline {\text{-(2n+2) -2n ...-2  0 2 ... 2n}}
\end{aligned} 
\;\right\}$ $2k-1$ \text{ times}

\noindent
(2n+2)** ...

\doublespace
\singlespace
\noindent Let  w=(c,d) be the black-black vector between the * points where the first  starred point is $C_{n}$ and

\doublespace
\singlespace
\noindent $d=6nk-2n+4nki-2ni+(12nk-4n-4k+2+8nki-4ni-2ki+i)cos2x+(12nk-4n-14k+5+8nki-4ni-8ki+4i)cos4x
+(12nk-4n-26k+9+8nki-4ni-16ki+8i)cos6x+(12nk-4n-38k+13+8nki-4ni-24ki+12i)cos8x+ ... +(22k-7+16ki-8i)cos(2n-2)x
+(10k-3+8ki-4i)cos2nx+(2k-1+2ki-i)cos(2n+2)x$

\doublespace
\noindent $c=(2+i)sin2x+(3+2i)sin4x+ ... +(3+2i)sin2nx+(1+i)sin(2n+2)x$

\singlespace
\noindent and let v=(a,b) be the blue-blue vector between the ** points where the first double starred point is $C_{m}$ and

\doublespace
\singlespace
\noindent $b=6nk-2n+4k-1+4nki-2ni+4ki-2i+(12nk-4n+6k-1+8nki-4ni+6ki-3i)cos2x+(12nk-4n-2k+1+8nki-4ni)cos4x
+(12nk-4n-14k+5+8nki-4ni-8ki+4i)cos6x+ ... +(22k-7+16ki-8i)cos(2n)x
+(10k-3+8ki-4i)cos(2n+2)x+(2k-1+2ki-i)cos(2n+4)x$

\doublespace
\noindent $a=(1+i)sin2x+(3+2i)sin4x+ ... +(3+2i)sin(2n+2)x+(1+i)sin(2n+4)x$

\singlespace
\noindent then $ad<bc$ which is equivalent to $sin4x+sin6x+... +sin(2n+2)x>0$ which is true since $(4n+6)x<180$. Hence by the Non-Corridor Test 8, there is no corridor.

QED

\doublespace
\singlespace
\noindent
\textbf{Cor:} \textbf{(Corridor Growth Rule 1)} In a corridor rhombus tower ending ... 6 $4^{2k-1}$ 2, then for $i\geq0$, $n\geq 1$, $k\geq1$ the subcode (where S= 6 $4^{2k-2}$ 6 $8^{2k-2}$)

\doublespace
\singlespace

... $(4n+4)^{2k-2}$ 4n+6 $(4n+8)^{2k-2}$ $[S+4n]^{i}$ 4n+6

$(4n+4)^{2k-1}$ $[S+4n-4]^{i+1}$ 4n+2

\quad \vdots \quad \vdots

$12^{2k-1}$ $[S+4]^{i+1}$ 10

$8^{2k-1}$ $[S]^{i+1}$ 6 $4^{2k-1}$ 2

\doublespace
\singlespace
\noindent forces (noting that by the symmetry of the code numbers of the first level this subcode must extend at least two spots to the left)

\doublespace
\singlespace
... 4n+6 $(4n+4)^{2k-2}$ 4n+6 $(4n+8)^{2k-2}$ $[S+4n]^{i}$ 4n+6

$(4n+4)^{2k-1}$ $[S+4n-4]^{i+1}$ 4n+2

\quad \vdots \quad \vdots

$12^{2k-1}$ $[S+4]^{i+1}$ 10

$8^{2k-1}$ $[S]^{i+1}$ 6 $4^{2k-1}$ 2

\doublespace
\singlespace
\noindent Proof:  Since if $k>1$

... $(4n+4)^{2k-1}$ 4n+6 $(4n+8)^{2k-2}$ $[S+4n]^{i}$ 4n+6

$(4n+4)^{2k-1}$ $[S+4n-4]^{i+1}$ 4n+2

\quad \vdots \quad \vdots

$12^{2k-1}$ $[S+4]^{i+1}$ 10

$8^{2k-1}$ $[S]^{i+1}$ 6 $4^{2k-1}$ 2 is impossible by Corridor Lemma 6 and ... 4n+2 $(4n+4)^{2k-2}$ 4n+6 is impossible by
Rhombus Rule C. If k=1, use Corridor Lemma 6 again and the fact that 4n+8 $(4n+6)^{2i+2}$ 4n+4 ... is impossible by Rhombus Rule C.

QED

\doublespace
\singlespace
\noindent
\textbf{Corridor Lemma 7:} Let S= 6 $4^{2k-2}$ 6 $8^{2k-2}$ and let ... 6 $4^{2k-1}$ 2 be the end of a corridor rhombus tower then for $i\geq-1$, $n\geq 1$, $k\geq1$ the corridor tower never ends with the subcode

\doublespace
\singlespace

... 4n+6 $(4n+8)^{2k-2}$ $( S+4n)^{i+1}$ 4n+6

\quad \quad \quad $(4n+4)^{2k-1}$ $( S+4n-4)^{i+1}$ 4n+2

\quad \quad \quad \quad \vdots \quad \vdots \quad \vdots \quad \vdots

 \quad \quad \quad $12^{2k-1}$ $( S+4)^{i+1}$ 10

\quad \quad \quad $8^{2k-1}$ $ S^{i+1}$ 6 $4^{2k-1}$ 2

\doublespace
\singlespace
\noindent Note 1: If k=1, n=1, $i=s-1$, then S=$6^2$,  S+4=$10^2$ and this rule says that no corridor tower ends in ... $10^{2s+2}$ 8 $6^{2s+1}$ 4 2 for $s\geq0$.

\noindent Note 2: If k=1, n=2, $i=s-1$, then this rule says that no corridor tower ends in ... $14^{2s+2}$ 12 $10^{2s+1}$ 8 $6^{2s+1}$ 4 2 for $s\geq0$.

\noindent Note 3: If k=1, $i=-1$, then this rule says that no corridor tower ends in ... $(4n+6)^2$ 4n+4 ... 8 6 4 2 for $n\geq1$. If we replace n by n-1, then no corridor tower ends in ... $(4n+2)^2$ 4n ... 8 6 4 2 for $n\geq2$.

\doublespace
Proof: Using backwards standard position

\doublespace
\singlespace
\noindent
\underline{0*}

\singlespace
\noindent
$\left.
\begin{aligned}
&\text{2 0}\\
&\text{-2 0}\\
&\quad \vdots \quad \vdots\\
&\underline {\text{2 0**}}
\end{aligned} 
\;\right\}$ $2k-1$ \text{ times}

\noindent
$\left.
\parbox{0.5\linewidth}
{\underline {-2 0 2}

\noindent
$\left.
\begin{aligned}
&\text{4 2 0 -2}\\
&\quad\vdots\quad\vdots\\
&\underline {\text{-4 -2 0 2}}
\end{aligned} 
\;\right\}$ $2k-2$ \text{ times}

\underline {4 2 0}

\noindent
$\left.
\begin{aligned}
&\text{-2 0}\\
&\quad \vdots \quad \vdots\\
&\underline {\text{2 0}}
\end{aligned} 
\;\right\}$ $2k-2$ \text{ times}

\underline {-2 0 2}
}
\right\}$ i+1 times

\noindent
$\left.
\begin{aligned}
&\text{4 2 0 -2}\\
&\text{-4 -2 0 2}\\
&\quad \vdots \quad \vdots\\
&\underline {\text{4 2 0 -2}}
\end{aligned} 
\;\right\}$ $2k-1$ \text{ times}

\noindent
$\left.
\parbox{0.5\linewidth}
{\underline {-4 -2 0 2 4}

\noindent
$\left.
\begin{aligned}
&\text{6 4 2 0 -2 -4}\\
&\quad\vdots\quad\vdots\\
&\underline {\text{-6 -4 -2 0 2 4}}
\end{aligned} 
\;\right\}$ $2k-2$ \text{ times}

\underline {6 4 2 0 -2}

\noindent
$\left.
\begin{aligned}
&\text{-4 -2 0 2}\\
&\quad \vdots \quad \vdots\\
&\underline {\text{4 2 0 -2}}
\end{aligned} 
\;\right\}$ $2k-2$ \text{ times}

\underline {-4 -2 0 2 4}
}
\right\}$ i+1 times

\quad\vdots\quad\vdots

\quad\vdots\quad\vdots

\quad\vdots\quad\vdots

\noindent
$\left.
\parbox{0.5\linewidth}
{\underline {-2n -(2n-2) ...-2 0 2 ... 2n}

\noindent
$\left.
\begin{aligned}
&\text{2n+2 2n ... 0 -2 ... -2n}\\
&\quad\vdots\quad\vdots\\
&\underline {\text{-(2n+2) -2n ... 0 2 ... 2n}}
\end{aligned} 
\;\right\}$ 2k-2 times

\underline {2n+2 2n ... 2 0 -2 ... -(2n-2)}

\noindent
$\left.
\begin{aligned}
&\text{-2n -(2n-2) ... 0 2 ... 2n-2}\\
&\quad \vdots \quad \vdots\\
&\underline {\text{2n 2n-2 ... 0 -2 ... -(2n-2)}}
\end{aligned} 
\;\right\}$2k-2 times

\underline {-2n -(2n-2) ... -2 0 2 ... 2n}
}
 \right\}$ i+1 times

\noindent
$\left.
\begin{aligned}
&\text{2n+2 2n ... 2 0 -2 ... -2n}\\
&\text{-(2n+2) -2n ... 0 2 ... 2n}\\
&\quad \vdots \quad \vdots\\
&\underline {\text{2n+2 2n ... 2 0 -2 ... -2n**}}
\end{aligned} 
\;\right\}$ $2k-1$ \text{ times}

\noindent
$\left.
\parbox{0.7\linewidth}
{\underline {-(2n+2) -2n ...-2 0 2 ... 2n 2n+2}

\noindent
$\left.
\begin{aligned}
&\text{2n+4 2n+2 ... 0 -2 ... -2n -(2n+2)}\\
&\quad\vdots\quad\vdots\\
&\underline {\text{-(2n+4) -(2n+2) ... 0 2 ... 2n+2}}
\end{aligned} 
\;\right\}$2k-2 times

\underline {2n+4 2n+2 2n ... 2 0 -2 ... -2n}

\noindent
$\left.
\begin{aligned}
&\text{-(2n+2) -2n ... 0 2 ... 2n}\\
&\quad \vdots \quad \vdots\\
&\underline {\text{2n+2 2n ... 0 -2 ... -2n}}
\end{aligned} 
\;\right\}$2k-2 times

\underline {-(2n+2) -2n ... -2 0 2 ... 2n+2}
}
\right\}$ i+1 times

\noindent
$\left.
\begin{aligned}
&\text{2n+4 2n+2 ... 2 0 -2 ... -(2n+2)}\\
&\text{-(2n+4) -(2n+2) ...-2  0 2 ... 2n+2}\\
&\quad \vdots \quad \vdots\\
&\underline {\text{-(2n+4) -(2n+2) ...-2  0 2 ... 2n+2}}
\end{aligned} 
\;\right\}$ $2k-2$ \text{ times}

\noindent
\underline {2n+4 2n+2 ... 2 0 -2 ... -2n}

\noindent
-(2n+2)* ...

\doublespace
\singlespace
\noindent Let  v=(a,b) be the blue-blue vector between the ** points where the first  starred point is $C_{m}$ and

\singlespace
$b=6nk-2n+4nki-2ni+(12nk-4n-2k+1+8nki-4ni-2ki+i)cos2x+(12nk-4n-10k+3+8nki-4ni-8ki+4i)cos4x
+(12nk-4n-22k+7+8nki-4ni-16ki+8i)cos6x+ ... +(26k-9+16ki-8i)cos(2n-2)x
+(14k-5+8ki-4i)cos2nx+(4k-2+2ki-i)cos(2n+2)x$

\doublespace
$a=(1+i)sin2x+(3+2i)sin4x+ ... +(3+2i)sin2nx+(2+i)sin(2n+2)x$

\singlespace
\noindent Let  w=(c,d) be the black-black vector between the * points where the first  starred point is $C_{n}$ and

\singlespace
$d=6nk-2n+8k-3+4nki-2ni+4ki-2i+(12nk-4n+12k-4+8nki-4ni+6ki-3i)cos2x+(12nk-4n+2k-1+8nki-4ni)cos4x
+(12nk-4n-10k+3+8nki-4ni-8ki+4i)cos6x+ ... +(26k-9+16ki-8i)cos2nx
+(14k-5+8ki-4i)cos(2n+2)x+(4k-2+2ki-i)cos(2n+4)x$

\doublespace
$c=(2+i)sin2x+(3+2i)sin4x+ ... +(3+2i)sin(2n+2)x+(2+i)sin(2n+4)x$

\singlespace
\noindent then $ad<bc$ which is equivalent to $sin4x+sin6x+... +sin(2n+2)x>0$ which is true since $(4n+6)x<180$. Hence by the Non-Corridor Test 8, there is no corridor.

QED

\doublespace
\singlespace
\noindent
\textbf{Cor:} \textbf{(Corridor Growth Rule 2)} In a corridor rhombus tower ending ... 6 $4^{2k-1}$ 2, then for $i\geq-1$, $n\geq 1$, $k\geq1$ the subcode (where S= 6 $4^{2k-2}$ 6 $8^{2k-2}$)

\doublespace
\singlespace

\quad \quad ... $(4n+8)^{2k-2}$ $( S+4n)^{i+1}$ 4n+6

\quad \quad \quad $(4n+4)^{2k-1}$ $( S+4n-4)^{i+1}$ 4n+2

\quad \quad \quad \quad \vdots \quad \vdots \quad \vdots \quad \vdots

 \quad \quad \quad $12^{2k-1}$ $( S+4)^{i+1}$ 10

\quad \quad \quad $8^{2k-1}$ $ S^{i+1}$ 6 $4^{2k-1}$ 2

\doublespace
\singlespace
\noindent forces  (noting that by the symmetry of the code numbers of the first level this subcode must extend at least two spots to the left)

\doublespace
\singlespace
\quad \quad ... $(4n+8)^{2k-1}$ $( S+4n)^{i+1}$ 4n+6

\quad \quad \quad $(4n+4)^{2k-1}$ $( S+4n-4)^{i+1}$ 4n+2

\quad \quad \quad \quad \vdots \quad \vdots \quad \vdots \quad \vdots

 \quad \quad \quad $12^{2k-1}$ $( S+4)^{i+1}$ 10

\quad \quad \quad $8^{2k-1}$ $ S^{i+1}$ 6 $4^{2k-1}$ 2

\doublespace
\singlespace
\noindent Proof: Since 

\quad \quad ... 4n+6 $(4n+8)^{2k-2}$ $( S+4n)^{i+1}$ 4n+6

\quad \quad \quad $(4n+4)^{2k-1}$ $( S+4n-4)^{i+1}$ 4n+2

\quad \quad \quad \quad \vdots \quad \vdots \quad \vdots \quad \vdots

 \quad \quad \quad $12^{2k-1}$ $( S+4)^{i+1}$ 10

\quad \quad \quad $8^{2k-1}$ $ S^{i+1}$ 6 $4^{2k-1}$ 2 is impossible by Corridor Lemma 7 and 4n+10 $(4n+8)^{2k-2}$ 4n+6 ... is impossible by Rhombus Rule C.

QED

\noindent Note: 1. If If k=1, n=1, $i=s-1$, then S=$6^2$,  S+4=$10^2$ and this rule says that 

... $10^{2s+1}$ 8 $6^{2s+1}$ 4 2 forces ... 12 $10^{2s+1}$ 8 $6^{2s+1}$ 4 2 for $s\geq0$

\noindent Note 2. If If k=1, n=2, $i=s-1$, then S=$6^2$,  S+4=$10^2$, S+8=$14^2$ and this rule says that 

... $14^{2s+1}$ 12 $10^{2s+1}$ 8 $6^{2s+1}$ 4 2 forces ... 16 $14^{2s+1}$ 12 $10^{2s+1}$ 8 $6^{2s+1}$ 4 2 for $s\geq0$

\doublespace
\singlespace
\noindent
\textbf{Corridor Lemma 8:}  A corridor rhombus tower cannot end in   

... 4n $(4n+2)^{2k}$ 4n $(4n-2)^{2k+1}$ ... $10^{2k+1}$ 8  $6^{2k+1}$ 4 2 with $k\geq0$, $n\geq2$ 

\singlespace
Proof: Using backwards standard position and the array below,  0* corresponds to the last C point $C_{n}$ and 0** corresponds to the special point $C_{m}$. The two single starred points are  black points and the two double starred points are blue points. 

\singlespace

\underline{0*}

\underline{2 0**}

$\left.
\begin{aligned}
&\text{-2 0 2}\\
&\text{4 2 0}\\
&\vdots \quad \vdots\\
&\underline {\text{-2 0 2}}
\end{aligned} 
\;\right\}$ $2k+1$ \text{ times}

\underline{4 2 0 -2}

$\left.
\begin{aligned}
&\text{-4 -2 0 2 4}\\
&\text{6 4 2 0 -2}\\
&\quad \vdots \quad \vdots\\
&\underline {\text{-4 -2 0 2 4}}
\end{aligned} 
\;\right\}$ $2k+1$ \text{ times}

\underline{6 4 2 0 -2 -4}

\quad \vdots \quad \vdots

\quad \vdots \quad \vdots

\underline{2n-2 2n-4 ... 2 0 -2 ... -(2n-4)}

$\left.
\begin{aligned}
&\text{-(2n-2) -(2n-4) ... -2 0 2 ... 2n-2}\\
&\text{2n 2n-2 ... 2 0 -2 ... -(2n-4)}\\
&\quad \vdots \quad \vdots\\
&\underline {\text{-(2n-2) -(2n-4) ... -2 0 2 ... (2n-2)*}}
\end{aligned} 
\;\right\}$ $2k+1$ \text{ times}

\underline{2n 2n-2 ... 2 0 -2 ... -(2n-4) -(2n-2)}

$\left.
\begin{aligned}
&\text{-2n -(2n-2)  ... -2 0 2 ... 2n-4 2n-2 2n}\\
&\text{2n+2 2n  ... 2 0 -2 ... -(2n-4) -(2n-2)}\\
&\quad \vdots \quad \vdots\\
&\underline {\text{2n+2 2n  ... 2 0 -2 ... -(2n-4) -(2n-2)}}
\end{aligned} 
\;\right\}$ $2k$ \text{ times}

\underline{-2n -(2n-2)  ... -2 0 2 ... 2n-2}

2n**

Let $w=(c,d)=((k+1)sin2x+(2k+1)sin4x+...+(2k+1)sin(2n-2)x+ksin2nx,2nk-2k+2n-2+(4nk-5k+4n-5)cos2x+(4nk-8k+4n-9)cos4x+(4nk-12k+4n-13)cos6x+...+(8k+7)cos(2n-4)x+(4k+3)cos(2n-2)x+kcos2nx)$  be the vector between the two black single starred points and $v=(a,b)=(ksin2x+(2k+1)sin4x+...+(2k+1)sin2nx+ksin(2n+2)x,2nk+2n-1+(4nk-k+4n-2)cos2x+(4nk-4k+4n-5)cos4x+(4nk-8k+4n-9)cos6x+...+(8k+7)cos(2n-2)x+(4k+3)cos2nx+kcos(2n+2)x)$ be the vector between the two blue double starred points then $ad<bc$ is equivalent to $sin4x+sin6x+...+sin2nx >0$  which holds since $4nx<180$. But by the Non-Corridor Test 8, this code sequence is impossible.

QED

\doublespace
\singlespace
\noindent
\textbf{Cor:} \textbf{(Corridor Growth Rule 3)} In a corridor rhombus tower ending ... 8  $6^{2k+1}$ 4 2, then for $k\geq0$, $n\geq2$

\doublespace
\noindent ...  $(4n+2)^{2k}$ 4n $(4n-2)^{2k+1}$ ... $10^{2k+1}$ 8  $6^{2k+1}$ 4 2 

\doublespace
\singlespace
 forces  (noting that by the symmetry of the code numbers of the first level this subcode must extend at least two spots to the left)

\doublespace
\noindent ...  $(4n+2)^{2k+1}$ 4n $(4n-2)^{2k+1}$ ... $10^{2k+1}$ 8  $6^{2k+1}$ 4 2

\singlespace
\noindent Proof: Since ...  4n $(4n+2)^{2k}$ 4n $(4n-2)^{2k+1}$ ... $10^{2k+1}$ 8  $6^{2k+1}$ 4 2 is impossible by Corridor Lemma 8 and 4n+4 $(4n+2)^{2k}$ 4n ... is impossible by Rhombus Rule C.

QED

\subsection*{15. More Growth Rules}

\singlespace
\textbf{Corridor Growth Rule 4:} Let  ... 6 $4^{2k-1}$  2 represent the end of a corridor tower with $k\geq1$ then it cannot end in the subcode

... $(4n)^{2k}$ 4n-2 ... $12^{2k-1}$ 10 $8^{2k-1}$ 6 $4^{2k-1}$  2 for $n\geq2$. 

\noindent Hence this means that if a corridor rhombus tower ends (and extends to the left for at least 2 spots) then

\doublespace
\noindent
... $(4n)^{2k-1}$ 4n-2 ... $12^{2k-1}$ 10 $8^{2k-1}$ 6 $4^{2k-1}$  2 

$\downarrow$

\noindent2
... 4n+2 $(4n)^{2k-1}$ 4n-2 ... $12^{2k-1}$ 10 $8^{2k-1}$ 6 $4^{2k-1}$  2 for $n\geq2$, $k\geq1$.

\singlespace
\noindent (since the only other choice ... 4n-2 $(4n)^{2k-1}$ 4n-2 ... is impossible by Rhombus Rule B.)

\doublespace
\singlespace
\noindent Proof:  Using backwards standard position

\doublespace
\singlespace
\noindent
\underline {0*}

\noindent
$\left.
\begin{aligned}
&\text{2 0}\\
&\text{-2 0}\\
& \vdots \quad \vdots\\
&\underline {\text{2 0**}}
\end{aligned} 
\;\right\}$ $2k-1$ \text{ times}

\noindent
\underline {-2 0 2}

\noindent
$\left.
\begin{aligned}
&\text{4 2 0 -2}\\
&\text{-4 -2 0 2}\\
&\vdots\quad\vdots\\
&\underline {\text{4 2 0 -2}}
\end{aligned} 
\;\right\}$ $2k-1$ \text{ times}

\noindent
\underline {-4 -2 0 2 4}

\noindent
\vdots \quad \vdots

\noindent
\vdots \quad \vdots

\noindent
\underline {-(2n-2) -(2n-4) ... -4 -2 0 2 ... (2n-2)*}

\noindent
$\left.
\begin{aligned}
&\text{2n 2n-2 ... 2 0 -2 ... -(2n-2)}\\
&\text{-2n -(2n-2) ... -2 0 2 ... (2n-2)}\\
&\vdots\quad\vdots\\
&\underline {\text{-2n -(2n-2) ... -2 0 2 ... (2n-2)}}
\end{aligned} 
\;\right\}$ $2k$ \text{ times}

\noindent 2n**

Let $w=(c,d)=(sin2x+sin4x+...+sin(2n-2)x,2nk-2k+(4nk-6k+1)cos2x+(4nk-10k+1)cos4x+...+(6k+1)cos(2n-4)x+(2k+1)cos(2n-2)x)$  be the vector between the two black single starred points and $v=(a,b)=(sin4x+sin6x+...+sin2nx,2nk-2k+1+(4nk-4k+2)cos2x+(4nk-6k+1)cos4x+ (4nk-10k+1)cos6x+...+(6k+1)cos(2n-2)x+(2k+1)cos2nx)$ be the vector between the two blue double starred points then $ad<bc$ is equivalent to $sin4x+sin6x+...+sin2nx >0$ which holds since $4nx<180$. Hence by the Non-Corridor Test 8, this code sequence is impossible.

QED

\doublespace
\singlespace
\textbf{Corridor Growth Rule 5:} Let  ... 6 $4^{2k-1}$  2 represent the end of a corridor rhombus tower with $k\geq1$ then it cannot end in the subcode

... 4n+2 $(4n+4)^{2k-2}$ 4n+2 ... $12^{2k-1}$ 10 $8^{2k-1}$ 6 $4^{2k-1}$  2 for $n\geq2$. 

\noindent Hence this means that if a corridor rhombus tower ends (and extends to the left for at least 2 spots) then

\doublespace
\singlespace
\noindent
... $(4n+4)^{2k-2}$ 4n+2 ... $12^{2k-1}$ 10 $8^{2k-1}$ 6 $4^{2k-1}$  2

\doublespace
$\downarrow$

\doublespace
\noindent
... $(4n+4)^{2k-1}$ 4n+2 ... $12^{2k-1}$ 10 $8^{2k-1}$ 6 $4^{2k-1}$  2 for $n\geq2$, $k\geq1$.

\singlespace
\noindent (since if $k>1$, the only other choice 4n+6 $(4n+4)^{2k-2}$ 4n+2 ... is impossible by Rhombus Rule C and if k=1 the other two choices namely .. 4n 4n+2 4n ... and ... $(4n+2)^{2}$ 4n ... $12$ 10 $8$ 6 $4$  2 are impossible by Rhombus Rule B and Corridor Lemma 7.)

\doublespace
\singlespace
Proof: Using backwards standard position so that 0* corresponds to the last C point $C_{n}$ and 0** corresponds to the special point $C_{m}$.

\doublespace
\singlespace

\noindent
\underline {0*}
\singlespace
\noindent
$\left.
\begin{aligned}
&\text{2 0}\\
&\text{-2 0}\\
&\vdots \quad \vdots\\
&\underline {\text{2 0**}}
\end{aligned} 
\;\right\}$ $2k-1$ \text{ times}

\noindent
\underline {-2 0 2}

\noindent
$\left.
\begin{aligned}
&\text{4 2 0 -2}\\
&\text{-4 -2 0 2}\\
&\vdots \quad \vdots\\
&\underline {\text{4 2 0 -2}}
\end{aligned} 
\;\right\}$ $2k-1$ \text{ times}

\noindent
\underline {-4 -2 0 2 4}

\noindent
\vdots \quad \vdots

\noindent
\vdots \quad \vdots

\noindent
\underline {-(2n-2) -(2n-4) ... -4 -2 0 2 .. (2n-2)}

\noindent
$\left.
\begin{aligned}
&\text{2n 2n-2 ... 2 0 -2 ... -(2n-2)}\\
&\text{-2n -(2n-2) ... -2 0 2 ... (2n-2)}\\
&\vdots \quad \vdots\\
&\underline {\text{2n 2n-2 ... 2 0 -2 ... -(2n-2)**}}
\end{aligned} 
\;\right\}$ $2k-1$ \text{ times}

\noindent
\underline{-2n -(2n-2) ... -2 0 2 ... 2n-2 2n}

\noindent
$\left.
\begin{aligned}
&\text{2n+2 2n 2n-2 ... 2 0 -2 ... -2n}\\
&\text{-(2n+2) -2n ... -2 0 2 ... 2n-2 2n}\\
&\vdots \quad \vdots\\
&\underline {\text{-(2n+2) -2n ... -2 0 2 ... 2n-2 2n}}
\end{aligned} 
\;\right\}$ $2k-2$ \text{ times}

\noindent
\underline{2n+2 2n  ... 2 0 -2 ... -(2n-2)}

\noindent
-2n*

\doublespace

\singlespace

Let $w=(c,d)=(sin2x+sin4x+...+sin(2n+2)x,2nk+2k-1+(4nk+2k-1)cos2x+(4nk-2k-1)cos4x+...+(6k-1)cos2nx+(2k-1)cos(2n+2)x)$  be the vector between the two black single starred points and $v=(a,b)=(sin4x+sin6x+...+sin2nx,2nk-2k+(4nk-4k)cos2x+(4nk-6k-1)cos4x+...+(6k-1)cos(2n-2)x+(2k-1)cos2nx)$ be the vector between the two blue double starred points, then $ad<bc$ is equivalent to $sin4x+sin6x+...+sin2nx >0$ which holds since $(4n+2)x<180$. Hence by the Non-Corridor Test 8, this code sequence is impossible.

QED

\noindent COR: \textbf{(Corridor Growth Rule 6)} Let $n\geq2$, $k\geq1$ and  ... 6 $4^{2k-1}$  2 represent the end of a corridor tower then the subcode

\doublespace
\singlespace
\noindent ... 4n+2 $(4n)^{2k-1}$ 4n-2 ... $12^{2k-1}$ 10 $8^{2k-1}$ 6 $4^{2k-1}$  2 (assuming it extends to the left at least 2k spots)

\doublespace
 $\downarrow$

\noindent ... $(4n+4)^{2k-1}$ 4n+2 $(4n)^{2k-1}$ 4n-2 ... $12^{2k-1}$ 10 $8^{2k-1}$ 6 $4^{2k-1}$  2

\noindent Proof:

\noindent ... 4n+2 $(4n)^{2k-1}$ 4n-2 ... $12^{2k-1}$ 10 $8^{2k-1}$ 6 $4^{2k-1}$  2 

 $\downarrow$ Rhombus Rule F

\noindent ... $(4n+4)^{2k-2}$ 4n+2 $(4n)^{2k-1}$ 4n-2 ... $12^{2k-1}$ 10 $8^{2k-1}$ 6 $4^{2k-1}$  2 

  $\downarrow$  Corridor Growth Rule 5 

\noindent ... $(4n+4)^{2k-1}$ 4n+2 $(4n)^{2k-1}$ 4n-2 ... $12^{2k-1}$ 10 $8^{2k-1}$ 6 $4^{2k-1}$  2

QED

\noindent
\textbf{Consequence 4} 

\singlespace
\noindent It follows that if  ... 6 $4^{2k-1}$  2 with $k\geq1$ represents the end of a corridor tower and if it expands to  ... $8^{2k-1}$ 6 $4^{2k-1}$  2 (noting that by the symmetry of the code numbers of the first level that the subcode must extend to the left at least two spots) then 

\doublespace
\noindent ... $8^{2k-1}$ 6 $4^{2k-1}$  2

  $\downarrow$ Corridor Growth Rule 4 with n=2

\noindent ... 10 $8^{2k-1}$ 6 $4^{2k-1}$  2 

\singlespace
  $\downarrow$ Corridor Growth Rule 6 (by symmetry, it extends to the left at least 2k spots)

\doublespace
\noindent ... $12^{2k-1}$ 10 $8^{2k-1}$ 6 $4^{2k-1}$  2

  $\downarrow$ Corridor Growth Rule 4  with n=3

\noindent ... 14  $12^{2k-1}$ 10 $8^{2k-1}$ 6 $4^{2k-1}$  2

\singlespace
  $\downarrow$ Corridor Growth Rule 6  (by symmetry, it extends to the left at least 2k spots)

\doublespace
\noindent ... $16^{2k-1}$ 14  $12^{2k-1}$ 10 $8^{2k-1}$ 6 $4^{2k-1}$  2

$\downarrow$

\noindent
\vdots \quad \vdots\quad\vdots\\

\singlespace
Evidently this means that the code sequence would grow without bound (noting that a corridor tower must start with a 2) and hence no corridor tower can end
... $8^{2k-1}$ 6 $4^{2k-1}$  2 with $k\geq1$.

\end{document}


\title  {On non-periodic and non-dense billiard trajectories Part 2}

\author{George William Tokarsky\\Math Department, University of Alberta}\maketitle
\singlespace\

\setcounter{page}{70}

\subsection*{16. Growth Trees}

\doublespace
\singlespace
\noindent
\textbf{The Growth tree:} Note by the symmetry of the code numbers of the first level, each subcode extends at least the required number of spots to the left.

\doublespace
\singlespace
\noindent Stage one: Let  S= 6 $4^{2k-2}$ 6 $8^{2k-2}$ and ... 6 $4^{2k-1}$  2 be the end of a corridor tower, then for $n\geq0$, $s\geq1$, $k\geq2$

\doublespace
\singlespace
... $(4n+8)^{2k-1}$ $( S+4n)^{s}$ 4n+6

 \quad \vdots \quad \vdots \quad \vdots \quad \vdots

$12^{2k-1}$ $( S+4)^{s}$ 10

$8^{2k-1}$ $ S^{s}$ 6 $4^{2k-1}$ 2

\doublespace
$\downarrow$ Rhombus Growth Rule 3

\singlespace

... 4n+10 $(4n+8)^{2k-1}$ $( S+4n)^{s}$ 4n+6

 \quad \vdots \quad \vdots \quad \vdots \quad \vdots

$12^{2k-1}$ $( S+4)^{s}$ 10

$8^{2k-1}$ $ S^{s}$ 6 $4^{2k-1}$ 2

\doublespace
$\downarrow$  Rhombus Rule F 

\singlespace

 ... $(4n+12)^{2k-2}$ 4n+10 $(4n+8)^{2k-1}$ $( S+4n)^{s}$ 4n+6

 \quad \vdots \quad \vdots \quad \vdots \quad \vdots

$12^{2k-1}$ $( S+4)^{s}$ 10

$8^{2k-1}$ $ S^{s}$ 6 $4^{2k-1}$ 2

\doublespace
$\downarrow$  Rhombus Growth Rule 6

\singlespace

 ... 4n+10 $(4n+12)^{2k-2}$ 4n+10 $(4n+8)^{2k-1}$ $( S+4n)^{s}$ 4n+6

 \quad \vdots \quad \vdots \quad \vdots \quad \vdots

$12^{2k-1}$ $( S+4)^{s}$ 10

$8^{2k-1}$ $ S^{s}$ 6 $4^{2k-1}$ 2

\doublespace
$\downarrow$ Rhombus Growth Rule 4 

\singlespace

... $(4n+8)^{2k-2}$ 4n+10 $(4n+12)^{2k-2}$ 4n+10 $(4n+8)^{2k-1}$ $( S+4n)^{s}$ 4n+6

 \quad \vdots \quad \vdots \quad \vdots \quad \vdots

$12^{2k-1}$ $( S+4)^{s}$ 10

$8^{2k-1}$ $ S^{s}$ 6 $4^{2k-1}$ 2

\doublespace
\singlespace
$\downarrow$  If $s\geq2$ use Rhombus Growth Rule 7, if s=1 use Corridor Growth Rule 1 with i=0 and n replaced by n+1.

\doublespace
\singlespace

... 4n+10 $(4n+8)^{2k-2}$ 4n+10 $(4n+12)^{2k-2}$ 4n+10 $(4n+8)^{2k-1}$ $( S+4n)^{s}$ 4n+6

 \quad \vdots \quad \vdots \quad \vdots \quad \vdots

$12^{2k-1}$ $( S+4)^{s}$ 10

$8^{2k-1}$ $ S^{s}$ 6 $4^{2k-1}$ 2

which is the same as (where  S+4n+4=4n+10 $(4n+8)^{2k-2}$ 4n+10 $(4n+12)^{2k-2}$)

\doublespace

\singlespace

\quad \quad \quad ... $( S+4n+4)^{1}$ 4n+10

$(4n+8)^{2k-1}$ $( S+4n)^{s}$ 4n+6

 \quad \vdots \quad \vdots \quad \vdots \quad \vdots

$12^{2k-1}$ $( S+4)^{s}$ 10

$8^{2k-1}$ $ S^{s}$ 6 $4^{2k-1}$ 2

\doublespace

\noindent This is the end of Stage one.

\singlespace
\noindent Stage two: Now suppose it has grown to the below where  $1\leq{t}\leq{s}$. If t=s, we jump to stage four and if t=s-1, we jump to the stage three. Hence we can assume that $1\leq{t}\leq{s-2}$.

\doublespace
\singlespace

\quad \quad \quad ... $( S+4n+4)^{t}$ 4n+10

$(4n+8)^{2k-1}$ $( S+4n)^{s}$ 4n+6

 \quad \vdots \quad \vdots \quad \vdots \quad \vdots

$12^{2k-1}$ $( S+4)^{s}$ 10

$8^{2k-1}$ $ S^{s}$ 6 $4^{2k-1}$ 2

\doublespace

$\downarrow$  Rhombus Rule F, $k\geq2$ Note for k=2, there is no change

\singlespace

... $(4n+12)^{2k-4}$ $( S+4n+4)^{t}$ 4n+10

$(4n+8)^{2k-1}$ $( S+4n)^{s}$ 4n+6

 \quad \vdots \quad \vdots \quad \vdots \quad \vdots

$12^{2k-1}$ $( S+4)^{s}$ 10

$8^{2k-1}$ $ S^{s}$ 6 $4^{2k-1}$ 2

\doublespace
\singlespace
$\downarrow$  Rhombus Growth Rule 1 (... $(4n+12)^{2k-4}$ [S+4n+4] $\longrightarrow$ $(4n+12)^{2k-3}$ [S+4n+4])

\doublespace
\singlespace

... $(4n+12)^{2k-3}$ $( S+4n+4)^{t}$ 4n+10

$(4n+8)^{2k-1}$ $( S+4n)^{s}$ 4n+6

 \quad \vdots \quad \vdots \quad \vdots \quad \vdots

$12^{2k-1}$ $( S+4)^{s}$ 10

$8^{2k-1}$ $ S^{s}$ 6 $4^{2k-1}$ 2

\doublespace
\singlespace
$\downarrow$  Rhombus Growth Rule 2  (... $(4n+12)^{2k-3}$ $( S+4n+4)^{t}$ 4n+10 $(4n+8)^{2k-1}$ ... $\longrightarrow$ $(4n+12)^{2k-2}$ 
$( S+4n+4)^{t}$ 4n+10 $(4n+8)^{2k-1}$ ...)

\doublespace

\singlespace

... $(4n+12)^{2k-2}$ $( S+4n+4)^{t}$ 4n+10

$(4n+8)^{2k-1}$ $( S+4n)^{s}$ 4n+6

 \quad \vdots \quad \vdots \quad \vdots \quad \vdots

$12^{2k-1}$ $( S+4)^{s}$ 10

$8^{2k-1}$ $ S^{s}$ 6 $4^{2k-1}$ 2

\doublespace

$\downarrow$  Rhombus Growth Rule 8

\singlespace

... 4n+10 $(4n+12)^{2k-2}$ $( S+4n+4)^{t}$ 4n+10

$(4n+8)^{2k-1}$ $( S+4n)^{s}$ 4n+6

 \quad \vdots \quad \vdots \quad \vdots \quad \vdots

$12^{2k-1}$ $( S+4)^{s}$ 10

$8^{2k-1}$ $ S^{s}$ 6 $4^{2k-1}$ 2

\doublespace
\singlespace
$\downarrow$ Rhombus Growth Rule 5  (... 4n+10 $(4n+12)^{2k-2}$ $[S+4n+4]^t$ 4n+10 $(4n+8)^{2k-1}$ ...   $\longrightarrow$ $(4n+8)^{2k-2}$ 4n+10 $(4n+12)^{2k-2}$ $[S+4n+4]^t$ 4n+10 $(4n+8)^{2k-1}$ ...) 

\doublespace
\singlespace

... $(4n+8)^{2k-2}$ 4n+10 $(4n+12)^{2k-2}$ $( S+4n+4)^{t}$ 4n+10

$(4n+8)^{2k-1}$ $( S+4n)^{s}$ 4n+6

 \quad \vdots \quad \vdots \quad \vdots \quad \vdots

$12^{2k-1}$ $( S+4)^{s}$ 10

$8^{2k-1}$ $ S^{s}$ 6 $4^{2k-1}$ 2

\doublespace
$\downarrow$  Rhombus Growth Rule 9

\singlespace

... 4n+10 $(4n+8)^{2k-2}$ 4n+10 $(4n+12)^{2k-2}$ $( S+4n+4)^{t}$ 4n+10

$(4n+8)^{2k-1}$ $( S+4n)^{s}$ 4n+6

 \quad \vdots \quad \vdots \quad \vdots \quad \vdots

$12^{2k-1}$ $( S+4)^{s}$ 10

$8^{2k-1}$ $ S^{s}$ 6 $4^{2k-1}$ 2

\doublespace
which is the same as

\singlespace

... $( S+4n+4)^{t+1}$ 4n+10

$(4n+8)^{2k-1}$ $( S+4n)^{s}$ 4n+6

 \quad \vdots \quad \vdots \quad \vdots \quad \vdots

$12^{2k-1}$ $( S+4)^{s}$ 10

$8^{2k-1}$ $ S^{s}$ 6 $4^{2k-1}$ 2

\doublespace
\singlespace
\noindent
We now go back to the beginning of stage two and keep repeating this process until we get

\doublespace

\singlespace

... $( S+4n+4)^{s-1}$ 4n+10

$(4n+8)^{2k-1}$ $( S+4n)^{s}$ 4n+6

 \quad \vdots \quad \vdots \quad \vdots \quad \vdots

$12^{2k-1}$ $( S+4)^{s}$ 10

$8^{2k-1}$ $ S^{s}$ 6 $4^{2k-1}$ 2

\doublespace
which we will designate as ... K

\noindent Stage three.

\singlespace
... K

$\downarrow$ Rhombus Rule F (Note if k=2 then there is no change)

... $(4n+12)^{2k-4}$ K

$\downarrow$ Rhombus Growth Rule 1 ... $(4n+12)^{2k-4}$ [S+4n+4] $\longrightarrow$ $(4n+12)^{2k-3}$ [S+4n+4]

... $(4n+12)^{2k-3}$ K

$\downarrow$ Rhombus Growth Rule 2

... $(4n+12)^{2k-2}$ K

$\downarrow$ Rhombus Growth Rule 8

... 4n+10 $(4n+12)^{2k-2}$ K

$\downarrow$ Rhombus Growth Rule 5 

... $(4n+8)^{2k-2}$ 4n+10 $(4n+12)^{2k-2}$ K

$\downarrow$  Corridor Growth Rule 1 with i=s-1 and n replaced by n+1.

... 4n+10 $(4n+8)^{2k-2}$ 4n+10 $(4n+12)^{2k-2}$ K 

\doublespace
which is the same as 

\doublespace
\singlespace
... $( S+4n+4)^{s}$ 4n+10

$(4n+8)^{2k-1}$ $( S+4n)^{s}$ 4n+6

 \quad \vdots \quad \vdots \quad \vdots \quad \vdots

$12^{2k-1}$ $( S+4)^{s}$ 10

$8^{2k-1}$ $ S^{s}$ 6 $4^{2k-1}$ 2

\doublespace
which we designate as ... L say then

\noindent Stage four

\singlespace
... L

$\downarrow$ Rhombus Rule F (Note if k=2 then there is no change)

... $(4n+12)^{2k-4}$ L

$\downarrow$ Rhombus Growth Rule 1

... $(4n+12)^{2k-3}$ L

$\downarrow$ Rhombus Growth Rule 2

... $(4n+12)^{2k-2}$ L

$\downarrow$ Corridor Growth Rule 2 with n replaced by n+1 and i+1=s

... $(4n+12)^{2k-1}$ L

\doublespace
which is the same as

\doublespace
\singlespace

... $(4n+12)^{2k-1}$ $( S+4n+4)^{s}$ 4n+10

$(4n+8)^{2k-1}$ $( S+4n)^{s}$ 4n+6

 \quad \vdots \quad \vdots \quad \vdots \quad \vdots

$12^{2k-1}$ $( S+4)^{s}$ 10

$8^{2k-1}$ $ S^{s}$ 6 $4^{2k-1}$ 2

\doublespace 
and we can go back to the start and repeat this process forever

QED

\singlespace
\noindent Note this means that if a corridor rhombus tower were to end in ... $8^{2k-1}$ $ S^{s}$ 6 $4^{2k-1}$ 2 for $s\geq1$, $k\geq2$, then the code numbers would grow without bound and hence no such corridor tower can exist.

\doublespace
\singlespace
\noindent
\textbf{Special cases of the Growth Tree}: Note by the symmetry of the code numbers of the first level, each subcode extends at least the required number of spots to the left.

\noindent 
\textbf{Growth Tree for  $k\geq1$, s=0}

\noindent
... $8^{2k-1}$ 6 $4^{2k-1}$ 2

\noindent
$\downarrow$ Corridor Growth Rule 4 (using n=2)

\noindent
... 10 $8^{2k-1}$ 6 $4^{2k-1}$ 2

\noindent
$\downarrow$ Corridor Growth Rule 6  (using n=2)

\noindent
...  $12^{2k-1}$ 10 $8^{2k-1}$ 6 $4^{2k-1}$ 2

\noindent
$\downarrow$ Corridor Growth Rule 4 (using n=3)

\noindent
... 14  $12^{2k-1}$ 10 $8^{2k-1}$ 6 $4^{2k-1}$ 2

\noindent
$\downarrow$ Corridor Growth Rule 6 (using n=3)

\noindent
... $16^{2k-1}$ 14  $12^{2k-1}$ 10 $8^{2k-1}$ 6 $4^{2k-1}$ 2

\noindent
\vdots\quad\vdots

\singlespace
\noindent Note 1. if  k=1,s=0 this becomes

\noindent
... 8 6 4 2 $\longrightarrow$ 
... 10 8 6 4 2
$\longrightarrow$
... 12 10 8 6 4 2
$\longrightarrow$
... 14 12 10 8 6 4 2
$\longrightarrow$
... 16 14 12 10 8 6 4 2
$\longrightarrow$ ...

\noindent Note 2. if  k=2,s=0 this becomes

\noindent
... $8^3$ 6 $4^3$  2 $\longrightarrow$ 
... 10 $8^3$  6 $4^3$ 2
$\longrightarrow$
... $12^3$ 10 $8^3$  6 $4^3$ 2
$\longrightarrow$
... 14 $12^3$ 10 $8^3$  6 $4^3$ 2
$\longrightarrow$
... $16^3$ 14 $12^3$ 10 $8^3$  6 $4^3$ 2
$\longrightarrow$ ...

\singlespace
\noindent Note 3. this means that if a corridor rhombus tower were to end in ... $8^{2k-1}$  6 $4^{2k-1}$ 2 for $s=0$, $k\geq1$, then the code numbers would grow without bound and hence no such corridor tower can exist.

\doublespace
\singlespace
\noindent
\textbf{Growth tree for  k=1, $s\geq1$}  where S = 6 6  since S= 6 $4^{2k-2}$ 6 $8^{2k-2}$

\doublespace
\singlespace
\noindent
... 8 $6^{2s+1}$ 4 2

\noindent
$\downarrow$ Rhombus Rule F 

\noindent
... $10^{2s}$ 8 $6^{2s+1}$ 4 2

\noindent
$\downarrow$ Corridor Growth Rule 3 (using n=2)

\noindent
... $10^{2s+1}$ 8 $6^{2s+1}$ 4 2

\noindent
$\downarrow$ Corridor Growth Rule 2 (using k=1, n=1, i=s-1)

\noindent
... 12 $10^{2s+1}$ 8 $6^{2s+1}$ 4 2

\noindent
$\downarrow$ Rhombus Rule F

\noindent
... $14^{2s}$ 12 $10^{2s+1}$ 8 $6^{2s+1}$ 4 2

\noindent
$\downarrow$ Corridor Growth Rule 3 (using n=3)

\noindent
... $14^{2s+1}$ 12 $10^{2s+1}$ 8 $6^{2s+1}$ 4 2

\noindent
$\downarrow$ Corridor Growth Rule 2 (using k=1, n=2, i=s-1)

\noindent
... 16 $14^{2s+1}$ 12 $10^{2s+1}$ 8 $6^{2s+1}$ 4 2

\noindent
\vdots

\singlespace
\noindent Note this means that if a corridor rhombus tower were to end in ... 8 $6^{2s+1}$ 4 2 for $s\geq1$, $k=1$, then the code numbers would grow without bound and hence no such corridor tower can exist.

\doublespace
\singlespace
\noindent \textbf{Unbounded Conclusion:}  If a corridor rhombus tower were to end in ... $8^{2k-1}$ $ S^{s}$ 6 $4^{2k-1}$ 2 for $s\geq0$, $k\geq1$ where S= 6 $4^{2k-2}$ 6 $8^{2k-2}$, then the code numbers would grow without bound and hence no such corridor tower can exist.

\subsection*{17. Locked in Rules}

\noindent
\singlespace
\textbf{Locked in Rule 1:} Let $k\geq1$, $N\geq1$ and  J= 6 $8^{2k-2}$ 6 $4^{2k-2}$ then if the code of a  corridor rhombus tower ends 
in ... 2 $4^{2k-1}$ $J^{N}$ 4 2 where 2 $4^{2k-2}$ 2 is its second level, then its boundary lines are locked in and are the same as  for the Corridor Rhombus tower 2 $J^{N}$ 4 2 and hence produce exactly the same ratio.

Proof: Using backwards standard position. Let $w=(c,d)=(Nsin2x+Nsin4x,4kN-2N+1+(6kN-3N+2)cos2x+(2kN-N)cos4x)$ be the black-black vector between $C_n$ and a lower black point (the two single starred points) and let $v=(a,b)=(Nsin2x+Nsin4x,4kN-2N+1+(6kN-3N+2)cos2x+(2kN-N)cos4x)$ be the blue-blue vector between the special blue point $C_m$ and a lower blue point (the two double starred points), then these vectors are parallel since they are exactly the same. Hence by the Locked In Test since ad=bc the two corridor rhombus towers have the same corridor boundary lines and produce exactly the same ratio.

\doublespace
\singlespace

\noindent
\underline {0*}

\noindent
$\left.
\begin{aligned}
&\text{2 0}\\
&\text{-2 0}\\
&\quad \vdots \quad \vdots\\
&\underline {\text{2 0**}}
\end{aligned} 
\;\right\}$ $2k-1$ \text{ lines}

\noindent
\underline {-2 0 2}

\noindent
$\left.
\begin{aligned}
&\text{4 2 0 -2}\\
&\quad\vdots\quad \vdots\\
&\underline {\text{-4 -2 0 2}}
\end{aligned} 
\;\right\}$2k-2 lines

\noindent
\underline {4 2 0}

\noindent
$\left.
\parbox{0.4\linewidth}
{ 
\noindent
$\left.
\begin{aligned}
&\text{-2 0}\\
&\quad\vdots\quad\vdots\\
&\underline {\text{2 0}}
\end{aligned} 
\;\right\}$2k-2 times

\underline {-2 0 2}

\noindent
$\left.
\begin{aligned}
&\text{4 2 0 -2}\\
&\quad\vdots\quad \vdots\\
&\underline {\text{-4 -2 0 2}}\\
\end{aligned} 
\;\right\}$2k-2 times

\underline {4 2 0}
}
\right\}$ N-1 times

\noindent  
$\left.
\begin{aligned}
&\text{-2* 0}\\
&\text{2 0}\\
&\quad \vdots \quad \vdots\\
&\underline {\text{-2 0}}
\end{aligned} 
\;\right\}$ $2k-1$ \text{ times}

\noindent 2**

\noindent Observing that all the integers between the two single starred points inclusive is exactly the vertical array of 2 $J^{N}$ 4 2 in backwards standard position with 0* corresponding to $C_n$, 0** to $C_m$ and -2* to $C_1$. Also observe that right boundary line of the first corridor and the location of  $C_m$ completely determines the other boundary lines.

QED

Note the special case when k=1 in which case J=$6^2$ and if a corridor rhombus tower ends in ... 2 4 $6^{2k}$ 4 2
then its boundary lines are locked in and produce exactly the same ratio as the corridor rhombus tower 2 $6^{2k}$ 4 2.

\singlespace
\textbf{Locked in Rule 2:} Let  $N\geq1$ and then if the code of a  corridor rhombus tower ends 
in ... 4 $6^{2N}$ 4 2 where 2  2 is its second level, then its boundary lines are locked in and are the same as  for the Corridor Rhombus tower 2 $6^{2N}$ 4 2 and hence produce exactly the same ratio.

Proof: Using backwards standard position. Let $w=(c,d)=(Nsin2x+Nsin4x,2N+1+(3N+2)cos2x+Ncos4x)$ be the black-black vector between $C_n$ and a lower black point (the two single starred points) and let $v=(a,b)=(Nsin2x+Nsin4x,2N+1+(3N+2)cos2x+Ncos4x)$ be the blue-blue vector between the special blue point $C_m$ and a lower blue point (the two double starred points), then these vectors are parallel since they are exactly the same. Hence by the Locked In Test since ad=bc the two corridor rhombus towers have the same corridor boundary lines and produce exactly the same ratio.

\doublespace
\singlespace

\noindent
\underline {0*}

\noindent
\underline {2 0**}

\noindent
$\left.
\begin{aligned}
&\text{-2 0 2}\\
&\text{4 2 0}\\
&\quad \vdots \quad \vdots\\
&\underline {\text{4 2 0}}
\end{aligned} 
\;\right\}$ $2N$ \text{ lines}

\noindent
\underline {-2* 0}

\noindent 2** ...

\noindent Observing that all the integers between the two starred points inclusive is exactly the vertical array of 2 $6^{2N}$ 4 2 in backwards standard position with 0* corresponding to $C_n$, 0** to $C_m$ and -2* to $C_1$.
QED

\textbf{Consequence 1:} If the code of a  corridor rhombus tower which we can consider to be in standard position ends in ... 2 $4^{2k-1}$ $J^{N}$ 4 2 with  J= 6 $8^{2k-2}$ 6 $4^{2k-2}$ and for $k\geq1$, $N\geq1$, then since its boundary lines are locked in and since the top portion of the tower contains a copy of the corridor rhombus tower 2 $J^{N}$ 4 2 with the same boundary lines, it must be one of its periodic corridor extensions and hence is exactly of the form 2 $K^s$ $ J^N$ 4 2 for some $s\geq1$ where K= $ J^N$ 4 $L^N$ 4 and L = $2^2$ $4^{2k-2}$.

\textbf{Consequence 2:} If the code of a  corridor rhombus tower which we can consider to be in standard position ends in ... 4 $6^{2N}$ 4 2 for $N\geq1$ then since its boundary lines are locked in and since the top portion of the tower contains a copy of the corridor rhombus tower 2 $6^{2N}$ 4 2 with the same boundary lines, it must be one of its periodic corridor extensions and hence is exactly of the form 2 $K^s$ $ 6^{2N}$ 4 2 for some $s\geq1$ where K= $ 6^{2N}$ 4 $2^{2N}$ 4.

\subsection*{18. Flow Charts}

\doublespace

\singlespace
We can now develop  a flow chart tree that produces \textbf{potential} corridor rhombus towers given that we know that the first level as it passes to the second level is of the form  ... 6 $4^{2k-1}$ 2 which we now will write as ... 6 $4^{2k-2}$ 4 2 with $k\geq1$. Actual existance of these potential corridor towers is then established later.

 We start with the case $k\geq2$. The special cases k=1 is discussed separately.  Note that for the code of a corridor tower \textbf{not to continue on the left}, it must have a symmetric first level with an even number of code numbers and start 2 6 ... . In other words that branch of the flow chart tree would have to end since otherwise it would violate Rhombus Rule A. 

Keep in mind that a corridor rhombus tower could also be written in the form 2 6 ... 6 2 plus 2 $4^{2k-2}$ 2 if we wanted to separate the codes of the first and second levels.

\doublespace
\singlespace

\noindent
\textbf{Flow Chart that produces the corridor towers for $k\geq2$.} Given a potential corridor rhombus tower ending in the subcode ... 6 $4^{2k-1}$ 2, then the following tree of subcodes is forced. Observe that each rule that requires the subcode to extend a certain number of spots to the left is satisfied because of the symmetry of the first level.

\doublespace

... 6 $4^{2k-2}$ 4 2 (same as ... 6 $4^{2k-1}$ 2)

$\downarrow$ \quad (Rhombus Rule F with n=2)

\singlespace
... $8^{2k-2}$ 6 $4^{2k-2}$ 4 2 $\longrightarrow $ ... $8^{2k-1}$ 6 $4^{2k-2}$ 4 2 grows without bound by the Unbounded Conclusion with s=0.

\doublespace
$\downarrow$ (10 $8^{2k-2}$ 6 ... impossible, Rhombus Rule C, n=3)

\singlespace
... 6 $8^{2k-2}$ 6 $4^{2k-2}$ 4 2 $\longrightarrow $  potential corridor tower of form \fbox{2 J 4 2} (ends since it can't extend to the left any further) where J = 6 $8^{2k-2}$ 6 $4^{2k-2}$ .

\doublespace
\singlespace
$\downarrow$ (Forcing Rule 5 with i=0, n=2 since the subcode keeps extending to the left)

\singlespace
... $4^{2k-2}$ J 4 2 $\longrightarrow$ ... $4^{2k-1}$J 4 2 $\longrightarrow$ ... 2 $4^{2k-1}$J 4 2 by Locked in Rule 1 (Noting that $4^{2k}$ J 4 2 Corridor Lemma 5 with s=k-1 and 6 $4^{2k-1}$ 6 Rhombus Rule B(b) with n=3 are impossible.) $\longrightarrow$ \fbox{periodic corridor extensions}

\doublespace
$\downarrow$ ( ... 2 $4^{2k-2}$ 6 impossible, Rhombus Rule C, n=1.)

.... 6 $4^{2k-2}$ J 4 2

$\downarrow$ (Rhombus Rule F with n=2)

\singlespace
... $8^{2k-4}$ 6 $4^{2k-2}$ J 4 2 $\longrightarrow$ ... 6 $8^{2k-4}$ 6 $4^{2k-2}$ 6 $8^{2k-2}$ ... impossible, Rhombus Rule H, example 2, s=k

\noindent (Note if k=2 this reads ... 6 $4^2$ J 4 2 $\longrightarrow$ ... 6 6 $4^2$ 6 $8^2$ ... impossible, Rhombus Rule H, example 1)

\doublespace
\singlespace
$\downarrow$  (If $k\geq3$, 10 $8^{2k-4}$ 6 ... is impossible by Rhombus Rule C using n=3 or if k=2, ... 4 6 4 ... is impossible by Rhombus Rule B(a) using n=2)

\doublespace
\singlespace
... $8^{2k-3}$ 6 $4^{2k-2}$ J 4 2 $\longrightarrow$ 10 $8^{2k-3}$ 6 $4^{2k-2}$ J 4 ... impossible, Rhombus Rule I, example 5 with $i=1$, $s=k$

\doublespace
$\downarrow$ ( ... 6 $8^{2k-3}$ 6 ... impossible, Rhombus Rule B(a), n=3.)

\singlespace
... $8^{2k-2}$ 6 $4^{2k-2}$ J 4 2 $\longrightarrow$ ... $8^{2k-1}$ 6 $4^{2k-2}$ J 4 2 grows without bound by the Unbounded Conclusion, (Since this is the same as ... $8^{2k-1}$ S 6 $4^{2k-1}$ 2 where S=6 $4^{2k-2}$ 6 $8^{2k-2}$.)

\doublespace
$\downarrow$ (10 $8^{2k-2}$ 6 ... impossible, Rhombus Rule C, n=3)

  ... $ J^2$ 4 2 $\longrightarrow$ potential corridor tower of form \fbox{2 $ J^2$ 4 2}

\singlespace
$\downarrow$ (Forcing Rule 5 with i=1, n=2 since the subcode keeps extending to the left)

\doublespace
\singlespace
... $4^{2k-2}$ $ J^2$ 4 2 $\longrightarrow$ ... $4^{2k-1}$ $ J^2$ 4 2 $\longrightarrow$ ... 2 $4^{2k-1}$ $ J^2$ 4 2 by Locked in Rule 1 (Noting that  $4^{2k}$ J 6 =  $4^{2k}$ 6 $8^{2k-2}$ 6 $4^{2k-2}$ 6 Rhombus Rule O and 6 $4^{2k-1}$ 6 Rhombus Rule B(b) with n=3 are impossible.)
$\longrightarrow$ \fbox{periodic corridor extensions}

\doublespace
$\downarrow$ ( ... 2 $4^{2k-2}$ 6 impossible, Rhombus Rule C, n=1.)

... 6 $4^{2k-2}$ $ J^2$ 4 2

$\downarrow$ (Rhombus Rule F with n=2)

\singlespace
... $8^{2k-4}$ 6 $4^{2k-2}$ $ J^2$ 4 2 $\longrightarrow$ ... 6 $8^{2k-4}$ 6 $4^{2k-2}$ 6 $8^{2k-2}$ ... impossible, Rhombus Rule H, example 2, s=k

\noindent (Note if k=2 this reads ... 6 $4^2$ $ J^2$ 4 2 $\longrightarrow$ ... 6 6 $4^2$ 6 $8^2$ ... impossible, Rhombus Rule H, example 1)

\doublespace
\singlespace
$\downarrow$ (If $k\geq3$, 10 $8^{2k-4}$ 6 ... is impossible by Rhombus Rule C using n=3 or if k=2,  ... 4 6 4 ... is impossible by Rhombus Rule B(a) using n=2)

\doublespace
\singlespace
... $8^{2k-3}$ 6 $4^{2k-2}$ $ J^2$ 4 2 $\longrightarrow$ 10 $8^{2k-3}$ 6 $4^{2k-2}$ $ J^2$ 4 ... impossible, Rhombus Rule I, example 5 with $i=2$, $s=k$

\doublespace
$\downarrow$  ( ... 6 $8^{2k-3}$ 6 ... impossible Rhombus Rule B(a), n=3.)

\singlespace
... $8^{2k-2}$ 6 $4^{2k-2}$ $ J^2$ 4 2 $\longrightarrow$ ... $8^{2k-1}$ 6 $4^{2k-2}$ $ J^2$ 4 2 grows without bound by the Unbounded Conclusion (since this is the same as ... $8^{2k-1}$ $ S^2$ 6 $4^{2k-1}$ 2 where S=6 $4^{2k-2}$ 6 $8^{2k-2}$.)

\doublespace
$\downarrow$ (10 $8^{2k-2}$ 6 ... impossible, Rhombus Rule C, n=3)

... $ J^3$ 4 2  $\longrightarrow$ potential corridor tower of the form \fbox{2 $ J^3$ 4 2}

\singlespace
$\downarrow$ (Forcing Rule 5 with i=2, n=2 if the subcode keeps extending to the left)

\doublespace
\singlespace
... $4^{2k-2}$ $ J^3$ 4 2 $\longrightarrow$ ... $4^{2k-1}$ $ J^3$ 4 2 $\longrightarrow$ ... 2 $4^{2k-1}$ $ J^3$ 4 2 by Locked in Rule 1 (Noting that  $4^{2k}$ J 6 =  $4^{2k}$ 6 $8^{2k-2}$ 6 $4^{2k-2}$ 6 Rhombus Rule O and 6 $4^{2k-1}$ 6 Rhombus Rule B(b) with n=3 are impossible.)
$\longrightarrow$ \fbox{periodic corridor extensions}

\noindent
\quad\quad\vdots\quad\vdots\quad\vdots

\noindent and so forth.

\doublespace
\noindent \textbf{Flow chart that produces the corridor towers for k=1}

\singlespace
... 6 4 2 $\longrightarrow$ ... 8 6 4 2 grows without bound by the Unbounded Conclusion with s=0, k=1.

\doublespace
$\downarrow$ (... 4 6 4 ... impossible, Rhombus Rule B, n=2).

... $6^{2}$ 4 2  $\longrightarrow$ potential corridor  tower of form \fbox{2 $6^{2}$ 4 2} 

$\downarrow$ (8 $6^{2}$ 4 ... impossible, Rhombus by Rule C, n=2). 

$\downarrow$  (... 4 $6^{2}$ 4 2 by Locked in Rule 2  $\longrightarrow$ \fbox{periodic corridor extensions})

\singlespace
... $6^{3}$ 4 2 $\longrightarrow$ ... 8 $6^{3}$ 4 2 grows without bound by the Unbounded Conclusion with s=1, k=1.

\doublespace
$\downarrow$  (... 4 $6^{3}$ 4 ... impossible, Rhombus Rule B, n=2).

... $6^{4}$ 4 2 $\longrightarrow$  potential corridor  tower of form  \fbox{2 $6^{4}$ 4 2}

$\downarrow$ (8 $6^{4}$ 4 ... impossible, Rhombus Rule C, n=2).

$\downarrow$ (... 4 $6^{4}$ 4 2 by Locked in Rule 2 $\longrightarrow$ \fbox{periodic corridor extensions}). 

\singlespace
... $6^{5}$ 4 2 $\longrightarrow$ ... 8 $6^{5}$ 4 2 grows without bound by the Unbounded Conclusion with s=2, k=1

\doublespace
$\downarrow$  (... 4 $6^{5}$ 4 ... impossible, Rhombus Rule B, n=2)

... $6^{6}$ 4 2  $\longrightarrow$ potential corridor  tower of form  \fbox{2 $6^{6}$ 4 2}

\noindent
\quad\quad\vdots\quad\vdots\quad\vdots

\noindent
and so forth

\singlespace
\noindent
\textbf{Conclusion:} it follows that the only potential corridor rhombus towers are of the form 2 $J^{N}$ 4 2 where J= 6 $8^{2k-2}$ 6 $4^{2k-2}$ for $k\geq1$, $N\geq1$ or one of its periodic corridor extensions.

\subsection*{19. Existance}

\doublespace
\singlespace
It is enough to show that corridor rhombus towers with code sequence 2 $J^{N}$ 4 2 exist since if they do, then so too do their periodic corridor extensions.

\doublespace
\singlespace
\textbf{Theorem:} Rhombus towers of the code sequence form 2 $J^{N}$ 4 2 where J=6 $8^{2k-2}$ 6 $4^{2k-2}$ and $N\geq1$, $k\geq1$ satisfy all three corridor conditions for all x such that $0<x<\alpha$ for some sufficently small $\alpha$ and hence are corridor rhombus towers for those x's.

\doublespace
\singlespace
\textbf{Proof:} We will assume the rhombus tower is in standard position and that the first reflection is in side AC.

Condition 1 is satisfied by starting the first reflection of the rhombus tower in side $A_{0}$$C_{1}$ and since the given code sequence has an odd number of code numbers in it, the last reflection must also be in side $A_{n}$$C_{n}$. Noting that if a code sequence of a rhombus tower has an odd number of even code numbers in it then the first and last reflection are in the same side.

Condition 2 is satisfied since the first level has an even length symmetric code sequence 2 6 $8^{2k-2}$ 6 $4^{2k-2}$ $J^{N-2}$ 6 $8^{2k-2}$ 6 2 and the second level also has an even length, symmetric code sequence 2 $4^{2k-2}$ 2. Noting that if a code sequence of a rhombus tower has an even number of even code numbers in it and its first reflection is in side AC then its last reflection is in side BC and vice versa. Further if the number of blue points in a rhombus tower equals the number of black points then the first and last appearance of side AB are parallel. This is always the case if its code sequence is symmetric and of even length.

It remains to show that condition 3 is satisfied.

Now if $k>1$, consider the vertical array corresponding to  2 $J^{N}$ 4 2 below in standard position and observe that all the 2N+1 single starred points are collinear black C points.

\doublespace
\singlespace

\noindent
\underline {0*}

\noindent
$\left.
\parbox{0.4\linewidth}
{

\underline {-2 0 2}

\noindent
$\left.
\begin{aligned}
&\text{4 2 0 -2}\\
&\quad\vdots\quad \vdots\\
&\underline {\text{-4 -2 0 2}}\\
\end{aligned} 
\;\right\}$2k-2 times

\underline {4* 2 0}

\noindent
$\left.
\begin{aligned}
&\text{-2 0}\\
&\quad\vdots\quad\vdots\\
&\underline{\text{2 0*}}
\end{aligned} 
\;\right\}$2k-2 times
}
\right\}$ N times

\noindent
\underline{ -2 0}

\noindent 2

This is easily seen by observing that if $k>1$, then $v=(a,b)=(sin4x,2k-1+(4k-2)cos2x+(2k-1)cos4x)$ is a vector between the first pair (and as well as between every succesive odd numbered pair) of starred points and $w=(c,d)=(sin2x,2k-1+(2k-1)cos2x)$ is a vector between the second pair (and as well as between every succesive even numbered pair) of starred points and since ad=bc all these starred points are all collinear black C points on a line L.  This line L has positive slope (2k-1+(4k-2)cos2x+(2k-1)cos4x)/sin4x since $0<8x<180$ and this is the same as $(N(4k-2)+N(6k-3)cos2x+N(2k-1)cos4x)/(Nsin2x+Nsin4x)$ which is the slope between the first and last starred points. 

If k=1, $N\geq1$, then if $v=(a,b)=(sin4x,1+2cos2x+cos4x)$ is a vector between the first pair of starred points and $w=(c,d)=(sin2x,1+cos2x)$ is a vector between the second pair of starred points then observe that all other vectors between successive starred points alternate between v and w. Now since ad=bc all these starred points are all collinear black C points on a line L with positive slope since $0<6x<180$.

\doublespace
\singlespace

\noindent
\underline {0*}

\noindent
$\left.
\begin{aligned}
&\text{-2 0 2}\\
&\underline{\text{4* 2 0*}}
\end{aligned} 
\;\right\}$N times

\noindent
\underline{ -2 0}

\noindent 2

\noindent Now further observe if $k>1$

1. That the Y coordinates of \textbf{all} C points strictly increase as we move down and from left to right along the vertical array as long as $8x<180$ since the Y coordinate of the first C point is 1 and this increases by either 1, cos2x or cos4x  as we go to the next successive C point. If $k=1$ then the Y coordinates of successive \textbf{black} C points increase since this increases by either 1, cos2x or $1+2cos2x+cos4x=4cos2xcos^{2}x$  as long as $6x<180$.

2. That the X coordinate of the (2i+1)st starred point is $isin2x+isin4x$ for $i\geq0$ and that the X coordinate of the (2i+2)th starred point is $isin2x+(i+1)sin4x$ for $i\geq0$.

3. That the X coordinate of every black C point between the (2i+1)st and the (2i+2)th starred points is either $(i+1)sin2x+(i+1)sin4x$ or $isin2x+(i+1)sin4x$ not including the (2i+1)st starred point.

4. That the X coordinate of every black C point between the (2i+2)th and the (2i+3)th starred points is $(i+1)sin2x+(i+1)sin4x$ not including the (2i+2)th starred point.

5. The slope of line L is strictly less than the slope of the right boundary line $C_{1}$$C_{n}$ of the first corridor which is given by $(N(4k-2)+1+(N(6k-3)+2)cos2x+N(2k-1)cos4x)/(Nsin2x+Nsin4x)$ for $8x<180$ since this inequality is equivalent to $sin2x+sin4x+sin6x>0$. If k=1, then the same result holds as long as $sin2x+sin4x>0$ which is the case since we must have $6x<180$.

\doublespace
\singlespace
\noindent If $k>1$, it follows from this that the slope from $C_{1}$ to any other black C point different from $C_{n}$ is $\leq$ slope of L and hence strictly less than the slope of $C_{1}$$C_{n}$ which means that their are no black C points in the strict interior of the two corridors if $8x<180$. If $k=1$, then since the only black C points not on line L (corresponding to the code number 2 located between 4* and 0*) have the same x coordinate as that of the adjacent 0* and since points 1 and 5 still hold with $6x<180$, we get the same result.

Interestingly if $k>1$ and we consider all the 2N starred blue C points below, they are also collinear on a line $L'$  which has exactly the same positive slope as the line L above with the last starred point being the special blue point $C_m$. Further this line $L'$ also goes through $C_n$ the last black C point.

\doublespace

Figure 30

\doublespace
\singlespace

\noindent
\underline {0}

\noindent
$\left.
\parbox{0.4\linewidth}
{

\underline {-2 0 2*}

\noindent
$\left.
\begin{aligned}
&\text{4 2 0 -2}\\
&\quad\vdots\quad \vdots\\
&\underline {\text{-4 -2 0 2}}\\
\end{aligned} 
\;\right\}$2k-2 times

\underline {4 2 0}

\noindent
$\left.
\begin{aligned}
&\text{-2* 0}\\
&\quad\vdots\quad\vdots\\
&\text{2 0}
\end{aligned} 
\;\right\}$2k-2 times
}
\right\}$ N times

\noindent
\underline{ -2 0}

\noindent 2

If k=1, we get the same result using the 2N+1 starred blue C points below.

\doublespace
\singlespace

\noindent
\underline {0}

\noindent
$\left.
\begin{aligned}
&\text{-2* 0 2*}\\
&\underline{\text{4 2 0}}
\end{aligned} 
\;\right\}$N times

\noindent
\underline{ -2* 0}

\noindent 2

\doublespace
Figure 31

\singlespace
Now by antisymmetry of the first level the line between any blue point $C_{m-i+1}$ with ${m-1}\geq{i}>1$ and the blue point $C_m$ is parallel to the line through the first black point $C_1$ and the black point $C_i$ and hence has slope $\leq$ that of the line $L'$. This means there can be no blue C points in the interior of the first level of the right corridor.

Finally all blue points $C_i$ with $i\geq m$ have the same X coordinate as $C_m$ and hence none of these blue points can lie in the interior of the second corridor.

Thus the interiors of the corridors are free of C points for $x<22.5$ if $k>1$ and for $x<30$ if $k=1$. Now since the coordinates of the A,B and C points in the tower are continuous functions of x, if we let $x\rightarrow0$ the A, B and C points will converge to the base triangle ABC. This means that for sufficiently small x the A and B points will also lie outside the two corridors, noting that the limit of the slope of $C_1$$C_n$ is infinity as x $\rightarrow$ 0. Hence condition 3 is satisfied for sufficiently small x.

QED

\doublespace
Figure 32

\noindent We give some more precise bounds on $\alpha$ below.
\singlespace

\subsection*{20. Existance Bounds}

By studying the behavior of the A and B points we can find bounds for $\alpha$. We will assume the corridor rhombus towers are in standard position and introduce a second coordinate system with the first B point being the origin and the distance between successive B points being one unit. The X-axis is the line through AB with the positive Y axis pointing upwards. We will call this the B coordinate system. Observe that the two coordinate systems have a different unit but the slope of any fixed line remains the same in both systems.

\doublespace
\singlespace
\noindent
\textbf{Conversion between the two systems:} If a point has coordinates (u,v) in the B system, then it has coordinates $(\frac{1}{2}cotx+2ucosx, \frac{1}{2}+2vcosx)$ in the original coordinate system.

Proof: Referencing Fig. 33, the origin of the B system has coordinates $(\frac{1}{2}cotx,\frac{1}{2})$ in the original system and since the distance between successive B points in the original system is $cscxcos(90-2x)=2cosx$, this represents the unit in the B system. It follows that (u,v) in the B system converts to $(\frac{1}{2}cotx+2ucosx, \frac{1}{2}+2vcosx)$ in the original coordinate system.

QED

\doublespace
Figure 33

\singlespace
\noindent
\textbf{Algorithm Three:} Given a corridor rhombus tower in standard position, we can calculate the coordinates of the B points in the B coordinate system as follows.

1. Introduce the corresponding vertical array startng with the black line 0. Note the only integers in the array will then be of the type -2, 0, 2 or 4.

2. Ignore all blue lines to get the \textbf{black vertical array}.

3. There is a one to one correspondence between the black C points and the B points in the tower omitting the first B point $B_0$. Note because of this correspondence, we can think of each integer in the black vertical array as representing its corresponding B point.

\noindent
a. The first B point $B_0$ has coordinates (0,0).

\noindent
b. For each integer of the following type add the following to the previous B coordinates

0: Add (-sinx,cosx)

2: Add (sinx,cosx)

4: Add (sin3x,cos3x)

-2: Add (-sin3x,cos3x)

\doublespace
\singlespace
Figure 34

\doublespace
\singlespace
\noindent Example: Given the tower 2 $6^{2}$ 4 2

\noindent
Step 1: The corresponding vertical array.
\doublespace
\singlespace

\noindent
\underline {0}

\noindent
-2 0 2

\noindent
\underline {4 2 0}

\noindent
\underline{-2 0}

\noindent
2

\noindent
Step 2: Ignore the blue lines to form the black vertical array.
\doublespace
\singlespace

\noindent
\underline {0}

\noindent
\underline {4 2 0}

\noindent
2

\doublespace
\singlespace

\noindent
Step 3: The successive coordinates of the B points in the B coordinate system after $B_0$(0,0) are then

\doublespace
\singlespace

\noindent
0: (-sinx,cosx)

\noindent 4: (-sinx+sin3x,cosx+cos3x)

\noindent 2: (sin3x,2cosx+cos3x)

\noindent 0: (-sinx+sin3x,3cosx+cos3x)

\noindent 2: (sin3x,4cosx+cos3x)

\doublespace
\singlespace

\noindent
\textbf{Fact:} Given the tower 2 $J^{N}$ 4 2 where J=6 $8^{2k-2}$ 6 $4^{2k-2}$ and $N\geq2$, $k\geq2$ and its black vertical array below the starred B points are collinear on a line $L_1$.

\doublespace
\noindent Proof:

\singlespace
 \noindent The starred points are collinear since the vector between successive starred points is the same.

\doublespace
\singlespace

\noindent
\underline {0}

\noindent
$\left.
\begin{aligned}
&\text{4 2 0 -2}\\
&\quad\vdots\quad \vdots\\
&\underline {\text{4 2 0 -2*}}\\
\end{aligned} 
\;\right\}$k-1 lines

\noindent
\underline {4 2 0}

\noindent
$\left.
\begin{aligned}
&\text{2 0}\\
&\quad\vdots\quad\vdots\\
&\underline{\text{2 0}}
\end{aligned} 
\;\right\}$k-1 lines

\noindent
$\left.
\begin{aligned}
&\text{4 2 0 -2}\\
&\quad\vdots\quad \vdots\\
&\underline {\text{4 2 0 -2*}}\\
\end{aligned} 
\;\right\}$k-1 lines

\noindent
\underline {4 2 0}

\noindent
$\left.
\begin{aligned}
&\text{2 0}\\
&\quad\vdots\quad\vdots\\
&\underline{\text{2 0}}
\end{aligned} 
\;\right\}$k-1 lines

\noindent
$\left.
\parbox{0.4\linewidth}
{ 

\noindent
$\left.
\begin{aligned}
&\text{4 2 0 -2}\\
&\quad\vdots\quad \vdots\\
&\underline {\text{4 2 0 -2*}}\\
\end{aligned} 
\;\right\}$k-1 lines

\underline {4 2 0}

\noindent
$\left.
\begin{aligned}
&\text{2 0}\\
&\quad\vdots\quad\vdots\\
&\underline{\text{2 0}}
\end{aligned} 
\;\right\}$k-1 lines
}
\right\}$ N-2 times

\noindent 
2

\doublespace
\noindent QED

\doublespace
\singlespace
\noindent Now working in the B coordinate system, noting that the B point corresponding to the first 0 has coordinates $(-sinx,cosx)$, observe

\doublespace
\singlespace

\noindent
1. that the Y coordinate of the B points strictly increases by either cosx or cos3x as we move from one B point to the next successive B point since $0<8x<180$. This means that successive B points are higher up.

\noindent
2. that the X coordinate of the B point corresponding to the (i+1)st appearance of -2* is $-sinx+isin3x$ for $i\geq0$.

\noindent
3. that the slope of $L_1$ in the B coordinate system is $((4k-2)cosx+(2k-1)cos3x)/sin3x$ which is positive since $0<x<22.5$.

\noindent
4. that the X coordinates of the B points corresponding to the integers between the (i+1)st (for $i\geq0$) appearance of -2* and the (i+2)th appearance of -2* (not including the  (i+1)st appearance of -2*) are of the form -sinx+(i+1)sin3x, (i+1)sin3x, (i+2)sin3x or -sinx+(i+2)sin3x. Note that $-sinx+(i+1)sin3x$ is the smallest number in this list since $0<x<22.5$.
It further follows since successive Y coordinates are increasing that this particular set of B points must lie to the right of the line $L_1$.

\noindent
5. that the X coordinates of the B points corresponding to the integers before the first appearance of -2* are of the form $-sinx$, $-sinx+sin3x$, $sin3x$ or 0. Note that -sinx which corresponds to the first -2* is the smallest number in this list. It also follows that this particular set of B points must lie to the right of the line $L_1$.

\noindent
6. that the X coordinates of the B points corresponding to the integers strictly after the last appearance of -2* are of the form -sinx+Nsin3x or Nsin3x. Note that -sinx+Nsin3x corresponds to the second last integer 0 and that Nsin3x corresponds to the  last integer 2 and that the slopes of the lines between the B point corrsponding to the last starred integer and the B point corresponding

a.  to the last integer 2 is ((2k+1)cosx+cos3x)/(sinx+sin3x) 

b. to the second last integer 0 is (2kcosx+cos3x)/(sin3x)

\noindent
Note that these slopes are $\leq$ the slope of $L_1$ since $k\geq2$ and $0<x<22.5$. It then follows that this particular set of B points must lie to the right of the line $L_1$. Finally note that the slope from a. is less than the slope from b. since that is equivalent to $(k-1)sin2x>0$.

\noindent
\textbf{Conclusion 1:} All the B points lie on the line $L_1$ or lie to the right of the line $L_1$.

\doublespace
Figure 35

\singlespace
\noindent
\textbf{Conclusion 2:} The slope of $L_1$ is less than the slope m=b/a of the boundary line of the right corridor  where $a= Nsin2x+Nsin4x$ and $b=(4k-2)N+1+((6k-3)N+2)cos2x+(2k-1)Ncos4x$.

Proof: Since $b>((4k-2)N+(6k-3)Ncos2x+(2k-1)Ncos4x)$ for $0<x<22.5$, it is enough to show that slope of  $L_1$ $\leq$ $((4k-2)N+(6k-3)Ncos2x+(2k-1)Ncos4x)/(Nsin2x+Nsin4x)$ which  is equivalent to

\noindent
(2cosx+cos3x)/sin3x $\leq$ (2+3cos2x+cos4x)/(sin2x+sin4x) which is true since both sides equal cotx.

QED

\noindent
\textbf{Conclusion 3:} The only B point that can lie on right boundary line of the right corridor is the B point corresponding to the first starred code number -2* with coordinates (-sinx,(2k-1)cosx+(2k-2)cos3x) in the B system.

\doublespace
\singlespace

\noindent
\textbf{Special k=1 Case:} Given the corridor rhombus tower 2 $6^{2N}$ 4 2 in standard position where $N\geq1$, $k=1$, the starred B points are collinear on a line $L_1$.

Proof: In the B coordinate system using the black vertical array below, the starred B points are collinear on a line $L_1$ since the vector between successive 0*'s is v=(a,b) where a=sin3x  and  b=2cosx+cos3x and the vector between the last two starred points is w=(c,d) where c=sinx and d=cosx and ad=bc. Note the slope of $L_1$ is cotx which is positive for $0<x<30$.

\doublespace
\singlespace

\noindent
\underline {0*}

\noindent
$\left.
\begin{aligned}
&\text{4 2 0*}\\
&\quad\vdots\quad \vdots\\
&\underline {\text{4 2 0*}}\\
\end{aligned} 
\;\right\}$N times

\noindent
 {2*}

QED

\doublespace
Figure 36

\noindent
Further observe 

\doublespace
\singlespace

\noindent
1. that the Y coordinates of the B points increase by either cos x or cos 3x as we move from one B point to the next successive B point. (Keep in mind that since $0<6x<180$ successive B points are higher up)

\noindent
2. that the X coordinate of the B point corresponding to the ith appearance of 0* is -sinx+(i-1)sin3x for $i\geq1$.

\noindent
3. that the X coordinates of the B points corresonding to the code numbers between the ith appearance of 0* and the (i+1)th appearance of 0* (not including the  ith appearance of 0*) are successively

4: -sinx+isin3x

2: isin3x

0*: -sinx+isin3x

\noindent 4. The last B point $B_n$ is on $L_1$ and the first B point $B_0$ is to the right of $L_1$.

\doublespace
\singlespace
\noindent
It follows that all B points are on the line $L_1$ or lie to the right of the line $L_1$ since $0<x<30$.

\doublespace
\singlespace
\noindent
5. that the slope of $L_1$ is less than the slope m=b/a of the boundary line of the right corridor  where a= Nsin2x+Nsin4x and b=2N+1+(3N+2)cos2x+Ncos4x since $acosx<bsinx$ is equivalent to $0<sin3x$ which is the case.

This means that the only B point that can lie on the right boundary line of  the right corridor is the B point corresponding to the first 0* and with coordinates (-sinx,cosx).

\doublespace
\singlespace

\noindent
\textbf{Special N=1, $k\geq2$ Case:} Given the corridor rhombus tower 2 $J$ 4 2  in standard position where J=6 $8^{2k-2}$ 6 $4^{2k-2}$  and  where  $k\geq2$, the starred  B points are collinear on a vertical line $L_1$.

Proof: In the B coordinate system using the black vertical array below, the starred points all have the same X coordinate namely -sinx.

\doublespace
\singlespace

\noindent
\underline {0*}

\noindent
$\left.
\begin{aligned}
&\text{4 2 0 -2*}\\
&\quad\vdots\quad \vdots\\
&\underline {\text{4 2 0 -2*}}\\
\end{aligned} 
\;\right\}$k-1 times

\noindent
\underline {4 2 0}

\noindent
$\left.
\begin{aligned}
&\text{2 0}\\
&\quad\vdots\quad\vdots\\
&\underline{\text{2 0}}
\end{aligned} 
\;\right\}$k-1 times

\noindent
2

QED

\noindent
Observe that in the B coordinate system

\noindent
1. that the Y coordinates of the B points increase by either cos x or cos 3x as we move from one B point to the next successive B point. (Keep in mind that since $0<8x<180$ successive B points are higher up)

\noindent
2. The X coordinates of all B points (except for $B_0$ which is 0) are -sinx, -sinx+sin3x or sin3x.

\noindent
3. The slope of $C_1$$C_n$ is (4k-1+(6k-1)cos2x+(2k-1)cos4x)/(sin2x+sin4x) which is positive for all $x\leq22.5$.

\noindent
4. The X coordinate of the last B point $B_n$ is sin3x and it follows that all B points with the same X coordinate lie beneath it on a vertical line. This means that the only B point on this vertical line that could possibly lie on the right boundary line of the right corridor is $B_n$ but this can't happen since ($<$$C_{n-3}$$C_n$$B_n$) = $90+x$ $\leq$ ($<$$C_1$$C_n$$B_n$) $\leq$ $90+3x$ = ($<$$C_{n-1}$$C_n$$B_n$) where $0<x<22.5$ .

\doublespace
Figure 37

\singlespace
\noindent
5. The X coordinate of the second last B point $B_{n-1}$ is -sinx+ sin3x and it follows that all B points with the same X coordinate lie beneath it on a vertical line. This means that the only B point on this vertical line that could possibly lie on the right boundary line of the right corridor is $B_{n-1}$ but this can't happen since ($<$$C_{n-3}$$C_n$$B_{n-1}$) = $90-3x$ $\leq$ ($<$$C_1$$C_n$$B_{n-1}$) $\leq$ $90-x$ = ($<$$C_{n-1}$$C_n$$B_{n-1}$) where $0<x<22.5$.

This means that the only B point that can lie on the right boundary line of  the right corridor is the B point with X coordinate $-sinx$ and which has the highest Y coordinate namely (2k-1)cosx+(2k-2)cos3x corresponding to the last -2*.

\doublespace
\singlespace
\noindent
\textbf{Conclusion:} In all cases the only B point that can lie on the right boundary line of  the right corridor of the corridor rhombus tower 2 $J^{N}$ 4 2 where J=6 $8^{2k-2}$ 6 $4^{2k-2}$ and $N\geq1$, $k\geq1$ is the B point with B-coordinates $(-sinx,(2k-1)cosx+(2k-2)cos3x)$ which corresponds to the coordinates $((cotx)/2-sin2x,(4k-1)/2+(4k-3)cos2x+(2k-2)cos4x)$ in the original coordinate system. Note that this \textbf{key B point} is in the first level.

\doublespace
\singlespace
\noindent
\textbf{Fact:} The key B point lies on the right boundary line of the right corridor of 2 $J^{N}$ 4 2 where J=6 $8^{2k-2}$ 6 $4^{2k-2}$ and $N\geq1$, $k\geq1$ if $(kN+1)cosx+[(3k-1)N+1]cos3x+[(4k-2)N+1)cos5x+(3k-2)Ncos7x+(k-1)Ncos9x=0$.

Proof: The only B point that can lie on the right boundary line of the right corridor would occur when the slope $m$ of $C_1$$C_n$ equals the slope $m_1$ from $C_1$ to the B point $((cotx)/2-sin2x,(4k-1)/2+(4k-3)cos2x+(2k-2)cos4x)$.

Now $m_1$ equals $((4k-3)/2+(4k-3)cos2x+(2k-2)cos4x)/((cotx)/2-sin2x)$

\noindent
=  $((4k-3)+(8k-6)cos2x+(4k-4)cos4x)/(cotx-2sin2x)$

\noindent
= $sinx((4k-3)+(8k-6)cos2x+(4k-4)cos4x)/cos3x$ using $cosx-2sinxsin2x=cos3x$

Setting this equal to the slope $m=(N(4k-2)+1+[N(6k-3)+2]cos2x+N(2k-1)cos4x)/(Nsin2x+Nsin4x)$
cross multiplying, simplifying and using the identities $sinx(sin2x+sin4x)=(cosx-cos5x)/2$ and $2cosnxcosmx=cos(n-m)x+cos(n+m)x$, we get the desired condition.

QED

\noindent Note that if $(kN+1)cosx+[(3k-1)N+1]cos3x+[(4k-2)N+1)cos5x+(3k-2)Ncos7x+(k-1)Ncos9x<0$ then the key B point lies to the left of the right boundary line and the corridor rhombus tower would not exist for those x's. 

\doublespace
\singlespace
\noindent
\textbf{Important Comment 1:} We don't have to worry about the A points since by antisymmetry an A point below the second level will lie on the left boundary line of the right corridor exactly when the antisymmetric B point lies on the right boundary line of the right corridor. Similarly an A point above the second level will lie on the left boundary line of the upper corridor exactly when the antisymmetric B point lies on the right boundary line of the upper corridor which as it turns out can never happen in a corridor rhombus tower because of the above conclusion.

\doublespace
Figure 38

\singlespace
It is worth noting that in a corridor rhombus tower, since $A_0$$B_0$ is parallel to $A_m$$B_m$ which is parallel to $A_n$$B_n$ and since successive B points (from which it follows that successive A points) have increasing Y coordinates, that any B point in the first level that lies on $C_1$$C_n$ must do so in the first level. Similarly for second level B points.

\doublespace
\singlespace
\noindent
\textbf{Important Comment 2:} Given the corridor rhombus tower 2 $J^{N}$ 4 2 with $N\geq1$, $k\geq1$, then if

1. k=1, the corridor tower only exists for those x's in the interval (0,30) for which $g_N(x)=(N+1)cosx+(2N+1)cos3x+(2N+1)cos5x+Ncos7x\geq0$.

2. $k>1$, the corridor tower only exists for those x's in the interval (0,22.5) for which $f_{k,N}(x)=(kN+1)cosx+[(3k-1)N+1]cos3x+[(4k-2)N+1)cos5x+(3k-2)Ncos7x+(k-1)Ncos9x\geq0$.

\doublespace
\singlespace
\noindent
\textbf{Consequences:}

\noindent
1. Observe that $g'_N(x)=-(N+1)sinx-(6N+3)sin3x-7N(sin5x+sin7x)-(3N+5)sin5x<0$ on (0,30) since sinx, sin3x, sin5x and sin5x+sin7x=2sin6xcosx are all positive on (0,30). It follows that $g_N(x)$ is strictly decreasing on (0,30) and since $g_N(0)=6N+3$ and $g_N(30)=-\sqrt3N$, $g_N(x)$ has an unique root on (0,30) given by

\noindent
 arctan$\sqrt{(6N+1-2\sqrt{8N^2+4N+1})/(2N-3)}$. Since $g_N(22.5)=cos22.5>0$, this unique root occurs in (22.5,30).

\doublespace
\singlespace
\noindent
\textbf{Conclusion 1:} This means that the corridor rhombus tower 2 $6^{2N}$ 4 2 for $N\geq1$ exists  if and only if $0<x\leq$arctan$\sqrt{(6N+1-2\sqrt{8N^2+4N+1})/(2N-3)}$.

\doublespace
\singlespace
\noindent
2. As a special case if k=1, N=1 the corridor rhombus tower 2 $6^2$ 4 2 exists for $2cosx+3cos3x+3cos5x+cos7x\geq0$ and x in (0,30) which only holds for $0<x\leq$arctan$\sqrt{2\sqrt{13}-7}$ or approximately $0<x\leq24.6768$. 

\doublespace
\singlespace
\noindent
3. As a further special case if k=1, N=2 the corridor rhombus tower 2 $6^4$ 4 2 exists for $3cosx+5cos3x+5cos5x+2cos7x\geq0$ and x in( 0,30) which only holds for $0<x\leq$arctan$\sqrt{13-2\sqrt{41}}$ or approximately $0<x\leq23.7578$.

\doublespace
\singlespace
Now $f(N)=(6N+1-2\sqrt{8N^2+4N+1})/(2N-3)$ is strictly decreasing on $[2,\infty)$ since its derivative $f'(N)<0$ is equivalent to $0<(2N-3)^2$ which holds on $[2,\infty)$. This means the unique root of $g_N(x)$ on (0,30) decreases as N increases on $[2,\infty)$. 

\doublespace
\singlespace
\noindent
\textbf{Conclusion 2:} Further since $\lim_{N \to +\infty}$arctan$\sqrt{(6N+1-2\sqrt{8N^2+4N+1})/(2N-3)}$=

\noindent arctan$\sqrt{3-\sqrt{8}}$ which equals 22.5 the corridor rhombus tower 2 $6^{2N}$ 4 2 with $N\geq1$ exists for all $0<x\leq22.5$ at least.

\doublespace
\singlespace
\noindent
4. For $N\geq1$ and $k\geq2$, $f_{k,N}(x)$ has an unique root in (18,22.5).

\textbf{Proof:} Consider the derivative $f'_{k,N}(x)=-(kN+1)sinx-3[(3k-1)N+1]sin3x-5[(4k-2)N+1)sin5x-7(3k-2)Nsin7x-9(k-1)Nsin9x$ which we can write in the form $f'_{k,N}(x)=-(kN+1)sinx-3[(3k-1)N+1]sin3x-5[(4k-2)N+1)sin5x-9(k-1)N(sin7x+sin9x)-(12Nk-5N)sin7x$ by rearranging its terms.

Now observe that since $sinx$, $sin3x$, $sin5x$, $sin7x$ and $sin7x+sin9x=2sin8xcosx$ are all positive on (0,22.5), $f'_{k,N}(x)$ is negative on (0,22.5) which means that $f_{k,N}(x)$ is strictly decreasing on (0,22.5).

Further since $f_{k,N}(18)=(N+1)cos18+(N+1)cos54>0$ and $f_{k,N}(22.5)=(3N(1-k)+1)cos22.5+N(1-k)cos(67.5)<0$ for $k\geq2$, $f_{k,N}(x)$ has an unique root in (0,22.5) which occurs in (18,22.5). 

QED

\doublespace
\singlespace
\noindent
\textbf{Conclusion 3:}  From the proof since $f_{k,N}(x)$ is strictly decreasing on (0,22.5) and has a root in (18,22.5), the corridor rhombus tower
2 $J^{N}$ 4 2 with $N\geq1$, $k\geq2$ exists for all $0<x\leq18$ at least.

\doublespace
\singlespace
\noindent
5. $f_{k+1,N}(x)<f_{k,N}(x)$ on (18,22.5) for $N\geq1$ and $k\geq2$.

\textbf{Proof:}
Note that  $f_{k+1,N}(x)-f_{k,N}(x)=Ncosx+3Ncos3x+4Ncos5x+3Ncos7x+Ncos9x$ and hence $f'_{k+1,N}(x)-f'_{k,N}(x)=-Nsinx-9Nsin3x-20Nsin5x-21Nsin7x-9Nsin9x=-Nsinx-9Nsin3x-20Nsin5x-9N(sin7x+sin9x)-12Nsin7x$. This derivative is $<0$ for all x in (0,22.5) since sinx, sin3x, sin5x, sin7x and sin7x+sin9x=2sin8xcosx are all $>0$ on (0,22.5). 

This means that $f_{k+1,N}(x)-f_{k,N}(x)$ is strictly decreasing on (0,22.5) and since $f_{k+1,N}(18)-f_{k,N}(18)=0$ and $f_{k+1,N}(22.5)-f_{k,N}(22.5)=-Ncos67.5-3Ncos22.5<0$, it follows that $f_{k+1,N}(x)<f_{k,N}(x)$ on (18,22.5). 

QED

\doublespace
\singlespace
\noindent
6. It follows that the unique root of $f_{k+1,N}(x)$ is strictly less than the unique root of $f_{k,N}(x)$ in (0,22.5) for $N\geq1$ and $k\geq2$ and they both occur in (18,22.5).

\doublespace
\singlespace
\noindent
7. For fixed $N\geq1$, the $\lim_{k \to +\infty}$ unique root of $f_{k,N}(x)$ is 18 ($\pi/10$ in radians).

\textbf{Proof:} We can assume that $k\geq6$. We first write $f_{k,N}(x)$ in the form $(kN-N)(cosx+cos9x)+(3k-2)N(cos3x+cos7x)+((4k-2)N+1)cos5x+(N+1)cosx+(N+1)cos3x$ and observe that cosx+cos9x=2cos5xcos4x, 

\noindent cos3x+cos7x=2cos5xcos2x and cos 5x are all $<0$ on (18,22.5). It follows that $f_{k,N}(x)\leq{g_{k,N}(x)}$ on (18,22.5) where $g_{k,N}(x)$=$2kNcos5x+4N$ since $(N+1)cosx+(N+1)cos3x\leq{4N}$ and $(4k-2)N+1\geq2kN$.

Now $g_{k,N}(18)=4N>0$, $g_{k,N}(22.5)={2kNcos112.5+4N}<0$ since $k\geq6$ and $g'_{k,N}(x)=-10kNsin5x<0$ on (18,22.5). This means that 18 $<$ the unique root of $f_{k,N}(x)$ $\leq$ the unique root of $g_{k,N}(x)$ on (18,22.5). But $g_{k,N}(x)=0$ when $x=.2arccos(-2/k)$. This means that the $\lim_{k \to +\infty}$ unique root of $f_{k,N}(x)$ is $\pi/10$.

QED

Also note that $f_{k,N+1}(x)-f_{k,N}(x)=kcosx+(3k-1)cos3x+(4k-2)cos5x+(3k-2)cos7x+(k-1)cos9x$.

Now if $f_{k,N}(\beta)=0$ where $\beta$ is in (18,22.5), then $f_{k,N+1}(\beta)=kcos\beta+(3k-1)cos3\beta+(4k-2)cos5\beta+(3k-2)cos7\beta+(k-1)cos9\beta$ but $f_{k,N}(\beta)= (kN+1)cos\beta+[(3k-1)N+1]cos3\beta+[(4k-2)N+1)cos5\beta+(3k-2)Ncos7\beta+(k-1)Ncos9\beta=N(kcos\beta+(3k-1)cos3\beta+(4k-2)cos5\beta+(3k-2)cos7\beta+(k-1)cos9\beta)+cos\beta+cos3\beta+cos5\beta=Nf_{k,N+1}(\beta)+cos\beta+cos3\beta+cos5\beta=0$.
\noindent And since $cos3x+cos5x=2cos4xcosx>0$ on (18,22.5), it follows that $f_{k,N+1}(\beta)=-(cos\beta+cos3\beta+cos5\beta)/N<0$.

\doublespace
\singlespace
\noindent
8. Noting that $f_{k,N+1}(18)=(N+2)cos18+(N+2)cos54>0$ this means that the unique root of $f_{k,N+1}(x)$ is strictly less than the unique root $\beta$ of $f_{k,N}(x)$ on (0,22.5) for $N\geq1$ and $k\geq2$ and they both occur in (18,22.5).

\doublespace
\singlespace
\noindent
Caution: It is not true that $\lim_{N \to +\infty}$ unique root of $f_{k,N}(x)$ on (18,22.5) is 18. For example if k=2, $f_{2,N}(19)$ is approximately $.3751N+1.4029>0$ whereas as previously noted $f_{2,N}(22.5)<0$ which means the unique root of $f_{2,N}(x)=0$ occurs in (19,22.5) for all N.

\doublespace
\singlespace
\noindent \textbf{Acknowledgements:}
Thanks to Arno Berger from the University of Alberta for listening to my thoughts and Boyan Marinov and David Szepesvari from the University of Alberta for developing computer programs that were useful in doublechecking the results of this paper.

\doublespace
\singlespace
\noindent \textbf{References:} \singlespace
\noindent [1] Galperin, G.A., Non-periodic and not everywhere dense billiard trajectories in convex polygons and polyhedrons, Commun. Math. Phys. 91, 187-211 (1983)

\doublespace
\singlespace
\noindent [2] Tokarsky, G.W, Galperin's Triangle Example, Commun. Math. Phys. 335, 1211-1213 (2015)

\doublespace
\singlespace
\noindent [3] Gutkin, Eugene, Billiards in polygons, Physica 19D, 311-333 (1986)

\doublespace
\singlespace
\noindent [4] Klee V., Wagon S., Old and New Unsolved Problems in Plane Geometry and Number Theory, (1991)

\doublespace
\singlespace
\noindent [5] Boshernitzan, M.D., Billiards and rational periodic directions in polygons, The American Mathematical Monthly 99(6), 522-529, (1992)

\doublespace
\singlespace
\noindent [6] Masur, H., Tabachnikov S., Rational billiards and flat structures, In: B. Hasselblatt, A. Katok (ed) Handbook of Dynamical Systems Volume 1A Elsevier Science B.V. (2002)

\doublespace
\singlespace
\noindent [7] Troubetzkoy, S., Periodic billiard orbits in right triangles, Annales de I'institut Fourier, 55:4, 1195-1217, (2005)

\doublespace
\singlespace
\noindent
[8] Bianca, C., On the mathematical transport theory in microporous media: the billiard approach, Nonlinear Analysis Hybrid Systems, 4, 699-735, (2010)

\doublespace
\singlespace
\noindent [9] Berger, M., Geometry Revealed A Jacob's Ladder to Modern Higher Geometry, Springer (2010)